\documentclass[a4paper]{report}
\usepackage{lmodern}
\usepackage[mathjax]{lwarp}
\usepackage{stmaryrd,amsmath,amssymb,ifthen,mathdots}
\usepackage{amsthm}
\usepackage{mathrsfs}
\usepackage{hyperref}
\usepackage{tikz}
\usetikzlibrary{positioning}
\usepackage{wrapfig}

\title{\Huge \bf Infinite Permutation Groups}
\author{\LARGE \bf Peter M.\ Neumann}
\date{\Large Michaelmas Term 1988--Trinity Term 1989\\Version: \today}

\newcommand{\mathmacro}[1]{#1\CustomizeMathJax{#1}}

\mathmacro{
\newcommand{\ms}[1]{\mathscr{#1}}
\newcommand{\mc}[1]{\mathcal{#1}}
}

\mathmacro{
\newcommand{\N}{\mathbb{N}}
\newcommand{\Q}{\mathbb{Q}}
\newcommand{\R}{\mathbb{R}}
\newcommand{\Z}{\mathbb{Z}}

\newcommand{\ep}{\varepsilon}
\newcommand{\normal}{\trianglelefteqslant}
\def\bdots{\mathinner{%
  \mkern1mu\raise1pt\hbox{.}%
  \mkern2mu\raise4pt\hbox{.}%
  \mkern2mu\raise7pt\vbox{\kern7pt\hbox{.}}\mkern1mu}}
}

\newtheorem{prop}{Proposition}
\newtheorem{thm}[prop]{Theorem}
\newtheorem{claim}[prop]{Claim}
\newtheorem{cor}[prop]{Corollary}
\newtheorem{lem}[prop]{Lemma}
\newtheorem{fact}[prop]{Fact}
\newtheorem{variant}[prop]{Variant}

\newtheorem*{coro}{Corollary}

\numberwithin{prop}{chapter}

\theoremstyle{definition}

\newtheorem{defn}[prop]{Definition}
\newtheorem{example}[prop]{Example}

\newtheorem*{definition}{Definition}
\newtheorem*{terminology}{Terminology}

\theoremstyle{remark}

\newtheorem{observation}[prop]{Observation}

\newtheorem{prob}{Problem}
\newtheorem{idleq}{Idle Question}
\newtheorem*{remark}{Remark}
\newtheorem*{comm}{Comment}
\newtheorem*{ex}{Example}
\newtheorem*{exs}{Examples}

\newtheorem*{question}{Question}

\newtheorem{note}{Note}

\usepackage{etoolbox}
\makeatletter
\patchcmd{\chapter}{\if@openright\cleardoublepage\else\clearpage\fi}{}{}{}
\makeatother

\counterwithout*{footnote}{chapter}

\newcommand\blankpage{
    \null
    \thispagestyle{empty}
    \newpage
    }

\mathmacro{
\newcommand{\nin}{\,\not\in\,}
\newcommand{\gen}[1]{\langle#1\rangle}
\newcommand{\Sym}{\mathrm{Sym}}
\newcommand{\Alt}{\mathrm{Alt}}
\newcommand{\FS}{\mathrm{FS}}
\newcommand{\BS}{\mathrm{BS}}
\newcommand{\AG}{\mathrm{AG}}
\newcommand{\AGL}{\mathrm{AGL}}
\newcommand{\ASL}{\mathrm{ASL}}
\newcommand{\AGammaL}{\mathrm{A\Gamma L}}
\newcommand{\PGL}{\mathrm{PGL}}
\newcommand{\PGammaL}{\mathrm{P\Gamma L}}
\newcommand{\GL}{\mathrm{GL}}
\newcommand{\PG}{\mathrm{PG}}
\newcommand{\PSL}{\mathrm{PSL}}
\newcommand{\Aut}{\mathrm{Aut}}
\newcommand{\Stab}{\mathrm{Stab}}
\newcommand{\supp}{\mathrm{supp}}
\newcommand{\Fix}{\mathrm{Fix}}
\newcommand{\im}{\mathrm{im}}
\newcommand{\bdisun}{\mathop{\dot{\bigcup}}}
\newcommand{\Wrr}{\mathrel{\mathrm{Wr}}}
\newcommand{\wrr}{\mathrel{\mathrm{wr}}}
\newcommand{\ha}{\upharpoonright}
\newcommand{\cof}{\mathrm{cof}}
\newcommand{\Orb}{\mathrm{Orb}}
\newcommand{\Hom}{\mathrm{Hom}}
\newcommand{\Homeo}{\mathrm{Homeo}}
\newcommand{\spa}{\mathrm{span}}
\newcommand{\charac}{\mathrm{char}}
}
\begin{document}

\maketitle

\blankpage
\pagenumbering{roman}

\chapter*{Preface}

What you have in front of you, dear reader, is a transcript of a lecture course on Infinite Permutation Groups given by Peter M.\ Neumann (1940--2020) in Oxford during the academic year 1988--1989. For information about Peter, his life and work, the reader is referred to the obituary that appeared in the \emph{Bulletin of the London Mathematical Society}\footnote{Martin W.\ Liebeck and Cheryl E.\ Praeger, \lq Peter Michael Neumann, 1940–2020\lq, \emph{Bull.\ London Math.\ Soc.} 54(2022), 1487--1514.}.  

The field of Group Theory was born out of the study of Finite Permutation Groups. But the study of infinite groups did not start with Infinite Permutation Groups.  Admittedly there were isolated results, such as the classification of normal subgroups of infinite symmetric groups from 1933 and 1934\footnote{J.\ Schreier and S.\ Ulam, \lq Ueber die Permutationsgruppe der nat\"urlichen Zahlenfolge\rq, \emph{Studia Math.} 4 (1933), 134-141 and R.\ Baer, \lq Die Kompositionreihe der Gruppe aller einein deutigigen Abbildungen einer unendlicher Menge auf sich\rq,  \emph{Studia Math.} 5 (1934), 15--17.}, but the study of Infinite Permutation Groups emerged as an independent field around 1980.   Peter's involvement with Infinite Permutation Groups started with two papers on finitary groups published in 1975 and 1976 (see Lectures 29, 30 and 31 in these notes) and from then on he was one of the leaders in this new field.   Most of the results described in these notes were at the time of the lectures brand new and had either just recently appeared in print or had not appeared formally.  A large part of the results described is either due to Peter himself or heavily influenced by him. These notes offer Peter's personal take on a field that he was instrumental in creating and in many cases ideas and questions that can not be found in the published literature.  

At the time these lectures were delivered the only available expository work on Infinite Permutation Groups were notes from Helmut Wielandt's lecture course \emph{Unendliche Permutationsgruppen} in T\"ubingen in 1959--1960\footnote{An English translation by P.\ V.\ Bruyns is printed in Helmut Wielandt, \emph{Mathematische Werke/Mathematical works. Vol.\ 1. Group theory.}  With essays on some of Wielandt's works by G.\ Betsch, B.\ Hartley, I.\ M.\ Isaacs, O.\ H.\ Kegel and P.\ M.\ Neumann. Edited and with a preface by Bertram Huppert and Hans Schneider. Walter de Gruyter \& Co., Berlin, 1994. xx+802~pp.}.   There are now other books treating Infinite Permutation Groups.  One can mention
Peter J.\ Cameron's book \emph{Oligomorphic Permutation groups}\footnote{Peter J.\ Camerons, \emph{Oligomorphic Permutation Groups}, London Mathematical Society Lecture Note Series, 152. Cambridge University Press, Cambridge, 1990. viii+160 pp.} and the book \emph{Notes on Infinite Permutation Groups}  that Peter himself had a part in writing\footnote{Meenaxi Bhattacharjee, Dugald Macpherson, Rögnvaldur G.\ Möller and Peter M.\ Neumann, \emph{Notes on infinite permutation groups}, Texts and Readings in Mathematics, 12, Hindustan Book Agency, New Delhi and co-published as volume 1698 of the Lecture Notes in Mathematics series by Springer-Verlag, Berlin, 1997. xii+202 pp.}.  Furthermore, textbooks on permutation groups by Brian Mortimer and John D.\ Dixon\footnote{John D.\ Dixon and Brian Mortimer, \emph{Permutation Groups}, Graduate Texts in Mathematics, 163. Springer-Verlag, New York, 1996. xii+346 pp.}, and by Peter J.\ Cameron\footnote{Peter J.\ Cameron, \emph{Permutation groups}, London Mathematical Society Student Texts, 45. Cambridge University Press, Cambridge, 1999. x+220 pp.} contain chapters on Infinite Permutation Groups.   Surprisingly the overlap of the present lecture notes with these books is small.  Thus the notes complement these books and treat in a unified way many topics that are not at all discussed in them.

The Oxford academic year is divided into three terms:  Michaelmas term, Hilary term and Trinity term, each spanning 8 weeks.  Peter's lecture course stretched over all three terms with two lectures per week (Mondays and Tuesdays).   As was Peter's habit the lectures are individually numbered and each lecture has a title.  Our sources are on one hand copies of Peter's own notes that he gave to David A.\ Craven and on the other hand notes that Röggi Möller took at the actual lectures.  There were some lectures missing from Peter's own notes and in these cases we rely solely on Röggi's notes.  The two sources are almost in 100\% agreement.   In preparing the notes the policy was to preserve the flavour of the original text with only the absolute minimum of editing.  The details of those references that had not appeared at the time have now been added and some additional references are mentioned in footnotes; these are mostly papers alluded to in the notes or papers relating to questions or comments in the text.  No attempt has been made to provide complete and up-to-date references for the topics treated in the notes. 

For us that took part in preparing these notes they are a way to pay tribute to Peter and his legacy.  Peter was a unique man; brimming with goodwill and enthusiasm towards everyone and everything.  His Mathematics have been highly influential but he was equally influential through his students and through the encouragement and support he offered to everyone.  We hope that by making these notes available more people will be smitten with Peter's enthusiasm for the subject and encouraged to explore the fascinating subject of Infinite Permutation Groups.  

\medskip
\noindent
{\bf Acknowledgements.}  We are very grateful to Sylvia Neumann for giving her approval to us making these notes available.  We are also grateful to Cheryl E.\ Praeger for her enthusiasm and support for this project. 

\bigskip\bigskip

\hspace{10cm}
David A.\ Craven (Birmingham)

\hspace{10cm}
Dugald Macpherson (Leeds)

\hspace{10cm}
Rögnvaldur G.\ Möller (Reykjavík)

\newpage

\tableofcontents
\newpage \blankpage

\pagenumbering{arabic}

\part{Michaelmas Term 1988}

\blankpage

\chapter{Introduction}

{\bf Sources.}

H.\ Wielandt:  Unendliche Permutationsgruppen, Lecture Notes, T\"ubingen 1960.

H.\ Wielandt:  Finite permutation groups, Academic Press, 1964.

$\Pi$MN et al:  Groups and Geometry.

\bigskip

\noindent
{\bf Prerequisites.}

(1)  Some knowledge of and sympathy for Group Theory.

(2)  The same for Set Theory.  Comfortable set theory + Axiom of choice.  Cardinal numbers are well-ordered.  If $n$ is a cardinal then $n^+$ is the next one.  If $m, n$ are infinite then
$$m+n=m\cdot n=\max\{m,n\}.$$

 \bigskip

\noindent
{\bf Notation.}

$\Omega$ is a set.  

$n=|\Omega|$, cardinal number of $\Omega$.

$\Sym(\Omega)=\{f:\Omega\to\Omega\mid f \mbox{ bijective}\}$.

When $n$ is infinite, $|\Sym(\Omega)|=2^n$.

If $f\in\Sym(\Omega)$ then 
\begin{align*}
\Fix(f)&=\{\omega\in\Omega\mid \omega f=\omega\},\\
\supp(f)&=\{\omega\in\Omega\mid \omega f\neq \omega\},\\
\deg(f)&=|\supp(f)|.
\end{align*}

If $G\leq \Sym(\Omega)$ then $G$ is said to be a \emph{permutation group}  and $n$ is known as its degree.

\begin{lem} Every permutation is uniquely (up to order) a product of disjoint cycles.
\end{lem}

Product may be infinite ($x\mapsto -x$ on $\Z$ or $x\mapsto 2x$ on $\Q$).

\bigskip

Let $G$ be a group.  A \emph{$G$-space} is a set $\Omega$ with an action
$$\Omega\times G\to\Omega; (\omega,g)\mapsto \omega g$$
such that (1) $\omega 1=\omega$ for all $\omega\in \Omega$ and (2) $(\omega g)h=\omega(gh)$ for all $\omega\in\Omega, g,h\in G$.

A $G$-space yields a \emph{permutation representation} of $G$, namely a homomorphism $G\to \Sym(\Omega)$.  Define $\rho_g:\Omega\to\Omega; \omega\mapsto \omega g$.  Then $\rho:G\to \Sym(\Omega); g\mapsto\rho_g$ is a homomorphism.

Conversely, given a permutation representation $G\to \Sym(\Omega)$ we get a natural $G$-space structure on $\Omega$.

Define an equivalence relation $\tau$ 
$$\alpha \equiv \beta \mod \tau\Leftrightarrow \big(\exists g\big) (\alpha g=\beta).$$

The equivalence classes are the \emph{$G$-orbits}.

The orbits are $G$-subspaces.

Any $G$-subspace is a union of orbits.

If there is only one orbit then $G$ is said to be \emph{transitive}.

\chapter{A Language for Permutation Group Theory}

If $H\leq G$ and $(G:H)=\{Hx\mid x\in G\}$ then $(G:H)$ is a $G$-space with action by right multiplication $Hx\to Hxg$. It is transitive.

\begin{thm}
If $\Omega$ is a transitive $G$-space then there is a subgroup $H$ such that $\Omega\cong (G:H)$.  
\end{thm}

\begin{proof}
Define $G_\alpha=\{g\in G\mid \alpha g=\alpha\}$ for $\alpha\in\Omega$ (\emph{stabilizer of }$\alpha$).  The map $\beta\to G_\alpha x$ if $\beta=\alpha x$ is an isomorphism $\Omega\to (G:G_\alpha)$.  
\end{proof}

If $\Omega_1, \Omega_2$ are $G$-spaces a map $f:\Omega_1\to \Omega_2$ is a $G$-morphism if 
$$(\omega g)f=(\omega f)g, \mbox{ for all } \omega\in \Omega_1, g\in G.$$

The class of $G$-spaces is a category. 

\begin{thm}   Let $\Omega_1, \Omega_2$ be $G$-spaces and suppose that $\Omega_1$ is transitive.  Let $\alpha_1\in\Omega_1$ and  $\alpha_2\in\Omega_2$.  Then there exists a morphism $f:\Omega_1\to \Omega_2$ such that $f:\alpha_1\to\alpha_2$ if and only if $G_{\alpha_1}\leq G_{\alpha_2}$.
\end{thm}

\begin{proof}
Suppose that $f:\Omega_1\to \Omega_2$ exists as required.  Let $g\in G_\alpha$.  Then 
$$\alpha_2g=(\alpha_1 f)g=(\alpha_1 g)f=\alpha_1 f=\alpha_2.$$
So $g\in G_{\alpha_2}$.  Thus $G_{\alpha_1}\leq G_{\alpha_2}$.

Conversely ...
\end{proof}

Here $f$ is unique.

\begin{note} If $\beta=\alpha x$ then $G_\beta=x^{-1}G_\alpha x$.
\end{note}

\begin{thm}
If $H, K$ are subgroups of $G$ then $(G:H)$ and $(G:K)$ are isomorphic if and only if there exists $x\in G$ such that $K=x^{-1}Hx$.
\end{thm}

\begin{proof}
Let $f:(G:H)\to (G:K)$ be an isomorphism.  If $f:H1\mapsto Kx$ then $H\leq G_{Kx}=x^{-1}Kx$ and same for $f^{-1}$.  Thus $K=x^{-1}Hx$.

Conversely look at $Ha\mapsto Kxa$ which is the required isomorphism.
\end{proof}

\begin{note}  If now $\Omega_1, \Omega_2$ are arbitrary $G$-spaces and if $f:\Omega_1\to \Omega_2$ is an isomorphism then $f$ maps orbits to orbits.
\end{note}

\begin{thm}
Let $\Omega$ be a transitive $G$-space, $H=G_\alpha$ for some $\alpha\in\Omega$.  Then 
$$\Aut(\Omega)=N_G(H)/H.$$
\end{thm}

\begin{proof}
Let $g\in N_G(H)$.  Define $f_g: \Omega\to \Omega$ by $\alpha x\to \alpha gx$.  If $\alpha x_1=\alpha x_2$ then $x_1x_2^{-1}\in H$ and so 
$$\alpha gx_1=\alpha gx_1x_2^{-1}x_2=\alpha gx_1x_2^{-1}g^{-1}gx_2
=\alpha gx_2.$$
Thus $f_g$ is well defined.  Check that $f_g$ is bijective and a $G$-morphism.
Thus we have a map $N_G(H)\to \Aut(\Omega)$.  Check that this is surjective and has kernel $H$.
\end{proof}

Define
\begin{align*}
\Orb(G,\Omega)&
=\{\Omega'\subseteq \Omega\mid\text{$\Omega'$ is a $G$-orbit}\},\\
\Omega^k&=\{(\omega_1, \ldots,\omega_k)\mid \omega_i\in \Omega\}
=\{f:k\to\Omega\},\\
\Omega^{(k)}&=\{(\omega_1, \ldots,\omega_k)\mid \omega_i\neq\omega_j, \mbox{ if }i\neq j\}
=\{f:k\to\Omega\mid f \mbox{ is injective}\},\\
\Omega^{\{k\}}&=\{\{\omega_1, \ldots,\omega_k\}\mid \omega_i\in \Omega, \omega_i\neq\omega_j, \mbox{ if }i\neq j\}
=\{\Gamma\subseteq \Omega\mid |\Gamma|=k\}.
\end{align*}

If $\Gamma\subseteq \Omega$ define 
\begin{align*}
G_{(\Gamma)}&=\{g\in G\mid \gamma g=\gamma, \forall \gamma\in\Gamma\}
=\bigcap_{\gamma\in\Gamma} G_\gamma,\\
G_{\{\Gamma\}}&=\{g\in G\mid \Gamma g=\Gamma\},\\ 
G^\Gamma&=G_{\{\Gamma\}}^\Gamma=\{g\in\Sym(\Gamma)\mid (\exists h\in G_{\{\Gamma\}})(h_{|\Gamma}=g)\}\cong G_{\{\Gamma\}}/G_{(\Gamma)}\qquad (G_{(\Gamma)}\triangleleft G_{\{\Gamma\}}). 
\end{align*}

The group $G$, or the $G$-space $\Omega$ is said to be
\begin{itemize}
\item \emph{$k$-transitive} if $G$ is transitive on $\Omega^{(k)}$,
\item \emph{highly transitive}:  if $G$ is $k$-transitive for all $k\in\N$,
\item \emph{$k$-homogeneous}:  if $G$ is transitive on $\Omega^{\{k\}}$,
\item \emph{highly homogeneous}:  if $G$ is $k$-homogeneous for all $k\in\N$.
\end{itemize}

\chapter{Separation Theorems I}
\setcounter{note}{0}
\begin{thm}[{{\hyperref[ref:lec3-4]{$\Pi$MN, 1976}}}\footnote{For those who do not know, Peter Neumann usually wrote his initials as $\Pi$MN rather than PMN.}]\label{3.1} Suppose that all $G$-orbits in $\Omega$ are infinite. If $\Gamma,\Delta$ are finite subsets of $\Omega$ then there exists $g\in G$ such that $\Gamma g\cap \Delta=\emptyset$.
\end{thm}
\begin{proof} Induction on $c=|\Gamma|$. For $c=0$ nothing to prove; $c=1$ follows from assumptions that all orbits are infinite.

Induction hypothesis: Assertion is true for sets $\Gamma$, $|\Gamma|=c-1$, and arbitrary finite subsets $\Delta'$. 

We may suppose that $\Gamma\not\subseteq \Delta$. Choose $\gamma_0\in \Gamma-\Delta$, let $\Gamma_0=\Gamma-\{\gamma_0\}$ and apply induction hypothesis to $\Gamma_0$.

Choose $g_0,g_1,\dots,g_d$ where $d=|\Delta|$ such that $\Gamma_0g_0\cap \Delta=\emptyset$, and
\begin{align*} \Gamma_0g_1&\cap (\Delta\cup\Delta g_0)=\emptyset,
\\ \Gamma_0g_2&\cap (\Delta\cup\Delta g_0\cup\Delta g_1)=\emptyset,
\\ &\vdots
\\ \Gamma_0g_d&\cap (\Delta\cup\Delta g_0\cup\cdots \cup\Delta g_{d-1})=\emptyset.\end{align*}

If $\gamma_0 g_i\not\in \Delta$ for some $i$, then $\Gamma g_i\cap \Delta=\emptyset$. So suppose that $\gamma_0 g_0,\gamma_0g_1,\dots,\gamma_0g_d\in \Delta$. Then since $|\Delta|=d$, there exist $p>q$ such that $\gamma_0g_p=\gamma_0g_q$. Let $g=g_pg_q^{-1}$; then $\gamma_0g=\gamma_0\not\in\Delta$ and $\Gamma_0g\cap \Delta=\emptyset$ because $\Gamma_0g_p\cap \Delta g_q=\emptyset$. So $\Gamma g\cap \Delta=\emptyset$.\end{proof}

\noindent
\textbf{Commentary on Theorem  \ref{3.1}:}

\begin{note} This is Lemma 2.3 in \hyperref[ref:lec3-4]{$\Pi$MN (1976)}. The proof is a small modification of the proof in \hyperref[ref:lec3-4]{Birch et al. (1976)}.
\end{note}

\begin{note}An immediate consequence\end{note}

\begin{cor}\label{3.2} Let $X$ be any group, $X_1,\dots,X_m$ cosets of subgroups $Y_1,\dots,Y_m$ such that $X=\bigcup X_i$. Then $[X:Y_i]$ is finite for some $i$.
\end{cor}
\begin{proof} Let $G=X$, $\Omega=\bigcup (X:Y_i)$ be the disjoint union of the coset spaces. Take $\Gamma=\{Y_1,\dots,Y_m\}$ and $\Delta=\{X_1,\dots,X_m\}$. If $g\in G$ then $g\in X_i$ for some $i$, so $Y_i g=X_i$. Thus $\Gamma g\cap \Delta\neq \emptyset$ for all $g\in G$, so we never get $\Gamma$ separate from $\Delta$. Hence at least one of the orbits is finite, so $[X:Y_i]<\infty$ for some $i$.
\end{proof}

\begin{thm}[{{\hyperref[ref:lec3-4]{BHN, 1954--55}}}]\label{3.3} If $X=\bigcup X_i$ (finite union) as in  \ref{3.2}, but the union is irredundant, then $[X:Y_i]$ is finite for all $i$ and $\sum 1/[X:Y_i]\geq 1$.
\end{thm}
\begin{proof} By  \ref{3.2} one of the $Y_i$ has finite index; so suppose that $[X:Y_m]$ is finite, and suppose that $Y_1,\dots,Y_p$ are different from $Y_m$, and $Y_{p+1}=\cdots =Y_m$. If $p=0$ then nothing to prove, since $Y_1=\cdots=Y_m$. So suppose that $p>0$ and for $1\leq j\leq p$, choose $x_j\in X_j-\bigcup_{i\neq j} X_i$. (Because of irredundancy of $\{X_i\}$.) The coset $Y_mx_j$ is covered by the finitely many cosets $Y_m x_j\cap X_i$ (finite or empty, discard the empty ones). If we reduce to an irredundant covering then $(Y_m x_j\cap X_j)$ remains. So $Y_m$ has an irredundant covering by finitely many cosets $(Y_m\cap X_i x_j^{-1})$, one of which is $(Y_m\cap X_jx_j^{-1})$, i.e., $Y_m\cap Y_j$. Induction tells us that $Y_m\cap Y_j$ has finite index in $Y_m$, hence in $X$. Thus $[X:Y_j]<\infty$.

Now let $Y=\bigcap Y_i$ (finite intersection of subgroups of finite index), so that $[X:Y]<\infty$. Counting cosets of $Y$, we get
\[ \sum_i [Y_i:Y]\geq [X:Y].\]
(Note that $[Y_i:Y]$ is the number of cosets of $Y$ in $X_i$.) Divide by $[X:Y]$ to get $\sum_i 1/[X:Y_i]\geq 1$.
\end{proof}

\begin{note} Conversely  \ref{3.1} is an immediate consequence of  \ref{3.2} (this is how it was originally proved).

For suppose that $\Gamma,\Delta$ are finite subsets of $\Omega$, and $\Gamma g\cap \Delta\neq \emptyset$ for all $g\in G$. For $\gamma\in\Gamma$, $\delta\in\Delta$, define
\[X_{\gamma,\delta}=\{g\in G\mid \gamma g=\delta\}.\]
Then $G=\bigcup X_{\gamma,\delta}$ (by our assumption). But also $X_{\gamma,\delta}=\emptyset$ (if $\gamma,\delta$ are in different orbits) as it is a coset of $G_\gamma$. By  \ref{3.2}, one of these stabilizers has finite index so there is a finite orbit.
\end{note}

\begin{note}
There is a quantitative version of  \ref{3.1}.\end{note}

\begin{thm}[{{\hyperref[ref:lec3-4]{Birch et al, 1976}}}] If every $G$-orbit has size greater than $|\Gamma|\cdot |\Delta|$, then there exists $g$ such that $\Gamma g\cap \Delta=\emptyset$.
\end{thm}

This bound is sharp: if $G=Z_c\times Z_d$ acting regularly (i.e., on the product set $Z_c\times Z_d$) and $\Gamma=Z_c$, $\Delta=Z_d$, then $\Gamma$ and $\Delta$ can \emph{not} be separated.

\chapter{Separation Theorems II}

\begin{thm}\label{4.1} Suppose that $k$ is an infinite cardinal number, $\Gamma,\Delta\subseteq\Omega$, with $\Gamma$ finite and $|\Delta|=k$, and that all orbits of $G$ have size greater than $k$. Then there exists $g\in G$ such that $\Gamma g\cap \Delta=\emptyset$.
\end{thm}
\begin{proof} Small modification of  \ref{3.1}. May assume $\Gamma\not\subseteq\Delta$. Proceed by induction on $|\Gamma|$.
\end{proof}

\begin{example}\label{4.2} Let $\Omega$ be any uncountable set. Let
\[ G=\FS(\Omega)=\{f\in\Sym(\Omega)\mid \deg f\text{ is finite}\}.\]
If $\Gamma=\Delta\subseteq\Omega$ and $\Gamma$ is infinite, then $\Gamma g\cap \Delta\neq \emptyset$ for all $g\in G$.

\smallskip
\noindent
\emph{Note.} $\FS(\Omega)$ is a union of countably many subgroups, each of index $n=|\Omega|$. In fact, if $\alpha_1,\alpha_2,\dots$ are distinct members of $\Omega$ then $G=\bigcup_{i=1}^\infty G_{\alpha_i}$ where $G=\FS(\Omega)$.
\end{example}

\begin{thm}[{{\hyperref[ref:lec3-4]{M.J. Tomkinson, 1987}}}]\label{4.3} Suppose that $|\Gamma|,|\Delta|\leq k$ (where $k$ is an infinite cardinal) and all orbits of $G$ in $\Omega$ have size greater than $2^k$. If $\Gamma\cap\Delta$ is finite, then there exists $g\in G$ such that $\Gamma g\cap \Delta=\emptyset$.
\end{thm}

The proof uses the following Ramsey-type theorem.

\begin{thm}[{{\hyperref[ref:lec3-4]{Erd\H{o}s--Hajnal--Rado, 1965}}}] \label{E+H+R}Let $k$ be an infinite cardinal number, and let $\Sigma$ be a set of cardinality greater than $2^k$. Suppose that the set $\Sigma^{\{2\}}=\bigcup_I \Pi_i$, where $|I|\leq k$. Then there exists $i\in I$ and there exist three distinct points $\alpha,\beta,\gamma\in \Sigma$ such that $\{\alpha,\beta\},\{\beta,\gamma\},\{\alpha,\gamma\}\in \Pi_i$. (Actually get a monochrome subset of cardinality $k+1$.)
\end{thm}

\begin{proof}[Proof of  \ref{4.3}] First reduce to the case where $|\Gamma\cap\Delta|=1$. Suppose we can do it for this case and use induction on $|\Gamma\cap\Delta|=d$. (Indeed, choose $\delta\in\Gamma\cap\Delta$ and let $\Gamma_0=\Gamma-\{\delta\}$. Then $|\Gamma_0\cap\Delta|=d-1$ and by inductive hypothesis there exists $g_0$ such that $\Gamma_0g_0\cap \Delta=\emptyset$. Then $|\Gamma g_0\cap\Delta|\leq 1$, and so by our supposed assumption there exists $g_1$ such that $\Gamma g_0g_1\cap\Delta=\emptyset$.)

So suppose now that $\Gamma\cap \Delta=\{\alpha\}$. Let $\Omega_0=\alpha G$ and (for convenience) choose a linear order $\leq$ on $\Omega_0$. Recall that $|\Omega_0|>2^k$. For $\xi\in\Omega_0$, choose $x_\xi\in G$ such that $\alpha x_\xi=\xi$. For $\gamma\in\Gamma$, $\delta\in\Delta$, define
\[ \Pi_{\gamma,\delta}=\left\{\{\xi,\eta\}\in\Omega_0^{\{2\}}\mid \xi<\eta\text{ and }\gamma x_\xi x_\eta^{-1}=\delta\right\}.\]

So suppose that $\Gamma g\cap \Delta\neq \emptyset$ for all $g\in G$. Then $\Omega_0^{\{2\}}=\bigcup_{\gamma,\delta} \Pi_{\gamma,\delta}$. By Theorem~ \ref{E+H+R}, there exist $\xi<\eta<\zeta$ such that, for some $\gamma,\delta$,
\[ \gamma x_\xi x_\eta^{-1}=\delta,\quad \gamma x_\eta x_\zeta^{-1}=\delta,\quad \gamma x_\xi x_\zeta^{-1}=\delta.\]
Then
\[ \gamma\underbrace{(x_\eta x_\zeta^{-1}\cdot x_\zeta x_\xi^{-1} \cdot x_\xi x_\eta^{-1})}_{=1}=\delta.\]
Now $\gamma\in\Gamma$, $\delta\in\Delta$ and $\Gamma\cap\Delta=\{\alpha\}$. Hence $\gamma=\delta=\alpha$. But $\alpha x_\xi x_\eta^{-1}=\xi x_\zeta^{-1}\neq \alpha$, a contradiction.
\end{proof}

\begin{example}  (Tomkinson, unpublished)  An example where $G$ has $k$ orbits, each of size $2^k$, and there are subsets $\Gamma, \Delta$ such that $|\Gamma|=|\Delta|=k$, $|\Gamma\cap\Delta|=1$ and $\Gamma g\cap \Delta\neq \emptyset$ for all $g\in G$.

\end{example}

\begin{thm}[Tomkinson] \label{Tomkinson} Suppose that $G$ is a group of cofinitary permutations of $\Omega$ (i.e., $\Fix(g)$ is finite for every $g\in G-\{1\}$). If $|\Gamma|=|\Delta|=k$ and all orbits have size greater than $2^{2^k}$, then there exists $g\in G$ such that $\Gamma g\cap \Delta=\emptyset$.
\end{thm}

\begin{prob} Does there exists a transitive group $G$ of degree $2^k$ such that there are subsets $\Gamma, \Delta$ of $\Omega$ with $|\Gamma|=|\Delta|=k$, $|\Gamma\cap\Delta|=1$ and $\Gamma g\cap \Delta\neq \emptyset$ for all $g\in G$?  Is it perhaps the case that if $G$ is transitive of degree $n\geq k^+$ and $|\Gamma|=|\Delta|=k$ with $|\Gamma\cap\Delta|=1$ then $(\exists g\in G)(\Gamma g\cap \Delta)=\emptyset$?  Can the bound in \ref{Tomkinson} be reduced?
\end{prob}

\noindent
{\bf Editors' note.}  The above lecture was given on 18.x.1988.  In the second lecture of the following week, on 25.x.1988, Peter presented brand new progress on the above problem by two members of the audience, Zo\'e Chatzidakis and  Peter Pappas.\footnote{This work appeared in Zo\'e Chatzidakis, Peter Pappas, M.\ J.\ Tomkinson, \lq Separation theorems for infinite permutation groups\rq. \emph{Bull.\ London Math.\ Soc.}, {22} (1990), 344--348.} 

\bigskip

\noindent
{\bf Theorem.}  \emph{Suppose $\Gamma, \Delta\subseteq \Omega$, $|\Gamma|,|\Delta|\leq k$ and let $G\leq\Sym(\Omega)$.  Suppose $S\subseteq G$ and $|S|>2^k$.  Then there exists $g\neq h\in S$ such that $\Gamma gh^{-1}\cap\Delta\subseteq \Fix(gh^{-1})$.}

\bigskip
\noindent 
{\em Proof.}  Assume there exists $S\subseteq G$ such that $|S|>2^k$ and for $\forall g\neq h\in S$, 
$$\Gamma gh^{-1}\cap\Delta\not\subseteq \Fix(gh^{-1}).$$
For such $g\neq h\in S$, $\exists \gamma\in \Gamma$, $\delta\in \Delta$ such that 
$\gamma gh^{-1}=\delta\nin\Fix(gh^{-1})$.  
Thus 
$$(*)\qquad \forall g\neq h\in S, \exists \gamma\in \Gamma, \delta\in \Delta, \gamma\nin \Fix(gh^{-1}) \mbox{ such that }\gamma g=\delta h, \mbox {in particular }\gamma\neq \delta.$$  Put a linear order $<$ on $S$ and $\forall (\gamma, \delta)\in \Gamma\times \Delta$ let 
\[\Pi_{\gamma, \delta}=\{\{g,h\}\in S\mid g<h, \gamma g=\delta h\}\quad\mbox{and}\quad I=\{(\gamma, \delta)\mid \gamma\neq\delta\}.\]
By (*), $S^{\{2\}}=\bigcup_I \Pi_{\gamma,\delta}$.  By Theorem~ \ref{E+H+R} $\exists g_1<g_2<g_3$ such that $\{g_1,g_2\}, \{g_2, g_3\}, \{g_1, g_3\}\in \Pi_{\gamma, \delta}$ for some $\gamma\neq \delta$.  Then $\gamma g_1=\gamma g_2$, $\gamma g_2=\gamma g_3$, $\gamma g_1=\gamma g_3$.  So $\delta g_2=\delta g_3$ and thus $\delta g_2g_3^{-1}=\delta=\gamma$ -- contradiction.  \hfill $\Box$

\begin{remark}
Note that $\Gamma gh^{-1}\cap\Delta=\Gamma\cap\Delta\cap \Fix(gh^{-1})$.
\end{remark}

\bigskip 
\noindent
{\bf Corollary.}
Suppose that $G$ is a group of cofinitary permutations of $\Omega$ (i.e., $\Fix(g)$ is finite for every $g\in G-\{1\}$). If $|\Gamma|=|\Delta|\leq k$ and all orbits have size greater than $2^k$, then there exists $g\in G$ such that $\Gamma g\cap \Delta=\emptyset$.

\bigskip

\noindent
\emph{Proof.}  \emph{Case 1.}  For some $g\in G$ the set $\Gamma g\cap \Delta$ is finite.  Then the result follows from Tomkinson's Theorem~\ref{Tomkinson}.

\emph{Case 2.}  For every $g\in G$ the set $\Gamma g\cap \Delta$ is infinite.  Apply previous theorem with $S=G$.  Then there exists $g\neq h\in G$ such that $\Gamma gh^{-1}\cap \Delta\subseteq \Fix(gh^{-1})$.  But the first set is infinite by assumption and the second set is finite -- contradiction.  \hfill $\Box$

\bigskip 

The Corollary is Tomkinson's Theorem.  The crucial case is where $\Gamma\cap\Delta=\{\alpha\}$.  The orbit $\alpha G$ has size $>2^k$.  Choose $S\subseteq G$ such that $\alpha g\neq \alpha h$ if $g\neq h\in S$ and $|S|>2^k$.  By remark there exists $g\neq h\in S$ such that 
\[\Gamma gh^{-1}\cap\Delta=\Gamma\cap\Delta\cap \Fix(gh^{-1})=\{\alpha\}\cap\Fix(gh^{-1})=\emptyset.\]

\begin{remark}
Tomkinson requires all orbits of length $>2^k$ and $|G|>2^{2^k}$.
\end{remark}

\chapter*{References for Lectures 3, 4}

\begin{itemize}
\item[1.] \label{ref:lec3-4}Peter M.\ Neumann, \lq The structure of finitary permutation groups\rq, \emph{Archiv der Math.}, 17 (1976), 3--17.
\item[2.] B.\ J.\  Birch, R.\ G.\ Burns, Sheila Oates Macdonald and Peter M.\ Neumann, \lq On the orbit sizes of permutation groups containing elements separating finite subsets\rq, \emph{Bull.\ Austral.\ Math.\ Soc.}, 14 (1976), 7--10.
\item[3.] B.\ H.\ Neumann, \lq Groups covered by permutable subsets\rq, \emph{J.\ London Math.\ Soc.}, 29 (1954), 236-248.
\item[4.] B.\ H.\ Neumann, \lq Groups covered by finitely many cosets\rq, \emph{Publ.\ Math.\ Debrecen}, 3 (1955), 227--242.
\item[5.] M.\ J.\ Tomkinson, \lq  Groups covered by abelian subgroups\rq, in \emph{Proc.\ Groups-St.\ Andrews 1985}, 332--334, L.M.S.\ Lecture Notes 121, CUP1987.
\item[6.]  P.\ Erd\"os, A.\ Hajnal and R.\ Rado, \lq Partition relations for cardinal numbers\rq, \emph{ Acta Sci.\ Math.\ Acad.\ Sci.\ Hung.}, 16 (1965), 93--196.  
\end{itemize}

[For 3, 4 see also \emph{Selected works of B.\ H.\ Neumann and Hanna Neumann} published by Charles Babbage Research Centre, Univ.\ of Manitoba, Winnipeg 1988.]

\chapter{The Symmetric Groups: Their Normal Subgroups}

Let $k$ be an infinite cardinal number. Let
\[ \BS(\Omega,k)=\{f\in\Sym(\Omega)\mid \deg f<k\},\]
and let
\[\FS(\Omega)=\BS(\Omega,\aleph_0)=\{f\mid \supp(f)\text{ is finite}\}.\]
This is the \emph{finitary symmetric group}. Let
\[ \Alt(\Omega)=\{f\in\FS(\Omega)\mid f\text{ is even}\}.\]

Note that:
\begin{enumerate}
\item $\deg f=\deg f^{-1}$.
\item $\supp(fg)\subseteq \supp(f)\cup\supp(g)$ (so $\deg fg\leq \deg f+\deg g$).
\end{enumerate}
Hence $\BS(\Omega,k)$ is a subgroup of $\Sym(\Omega)$.   Also, $\supp(g^{-1}fg)=(\supp(f))g$ (so $\deg g^{-1}fg=\deg f$). Hence $\BS(\Omega,k)\normal \Sym(\Omega)$.

\begin{thm}\label{5.1} Let $\Omega$ be an infinite set, $n=|\Omega|$, $S=\Sym(\Omega)$. Then the sequence
\[ \{1\}<\Alt(\Omega)<\FS(\Omega)<\BS(\Omega,\aleph_1)<\cdots<\BS(\Omega,n)<S\]
contains \emph{all} the normal subgroups of $S$: indeed, all subnormal subgroups. In particular, $\BS(\Omega,k^+)/\BS(\Omega,k)$ is simple.
\end{thm}

Define a type sequence of degree $n$ to be a sequence of the form $\aleph_0^{a_0}1^{a_1}2^{a_2}\dots$, where $a_i$ are cardinal numbers such that $a_0\aleph_0+a_1+2a_2+\cdots=n$. If $f\in S$ then the cycle type of $f$ is defined to be $\aleph_0^{a_0}1^{a_1}2^{a_2}\dots$, where $f$ has $a_0$ infinite cycles and $a_r$ finite $r$-cycles (for $r=1,2,\dots)$ in its cycle decomposition.

\begin{lem} Permutations $f_1$, $f_2$ are conjugate in $S$ if and only if they have the same cycle type. If $f_1,f_2$ lie in $\BS(\Omega,k)$ the same is true there. The same is true in $\Alt(\Omega)$.
\end{lem}

\begin{lem}\label{5.3} If $f\in S$ then there exist $t_1,t_2\in S$ such that $f=t_1t_2$ and $t_1^2=t_2^2=1$. Furthermore, $t_1,t_2$ can be chosen so that $\supp(t_i)\subseteq \supp(f)$.
\end{lem}
\begin{proof} Define $t_1,t_2$ on each cycle of $f$: if $f:\alpha_i\mapsto \alpha_{i+1}$ ($i\in \Z$ if the cycle is infinite, $i\in\Z/r\Z$ if it is a $r$-cycle), define
\begin{align*} t_1&:\alpha_i\mapsto \alpha_{-i},
\\ t_2&:\alpha_i\mapsto \alpha_{-i+1}.
\end{align*}
Then $t_1^2=1=t_2^2$ and $t_1t_2=f$. Also if $\alpha$ is fixed by $f$ then it is also fixed by $t_1,t_2$. So $\supp(t_i)\subseteq \supp(f)$.
\end{proof}

\begin{lem}\label{5.4} Suppose that $f\in S$ and $\aleph_0\leq \deg f<k$. If $t\in S$ and $t^2=1$ and $\deg t\leq \deg f$ then $t=f_1f_2f_3f_4$, where each $f_i$ is conjugate to one of $1$, $f$, $f^{-1}$.
\end{lem}
\begin{proof} Let $m$ be the number of non-trivial cycles of $f$. Let $t$ be of cycle type $1^{a_1}2^{a_2}$ (so $2a_2\leq \deg f$). We consider three cases.

\medskip

\noindent \textbf{Case 1:} There exists $a$ such that $a_2=2a$ and $m\geq a_2$. Then we can choose elements $\alpha_i,\beta_i$ ($i\in I$) such that $|I|=a$ and they all lie in different non-trivial cycles of $f$. Let $t_1=\prod_{i\in I}(\alpha_i,\beta_i)$ and $t_2=[t_1,f]=t_1^{-1}f^{-1}t_1f$ (a product of two conjugates of $f$). Then
\[ t_2=\prod_{i\in I} (\alpha_i,\beta_i)\prod_{i\in I} (\alpha_if,\beta_if).\]
Since all visible transpositions are pairwise disjoint, this is the cycle decomposition of $t_2$. Thus the cycle type of $t_2$ is $1^x2^{a_2}$ for some $x$. If $a_2<n$ then $x=n=a_1$. So $t_1,t_2$ are conjugate in $\BS(\Omega,k)$. If $a_2=n$, then $0\leq x\leq n$ but $t$ can be written as a product of two conjugates of $t_2$.

\medskip

\noindent \textbf{Case 2:} $a_2=2a$ and $m<a_2$. Then $m<\deg f$, so $\deg f=\aleph_0$ and $m$ is infinite. Therefore $f$ has at least one infinite cycle:
\[ f:\alpha_i\mapsto \alpha_{i+1}\qquad (i\in\Z).\]
Take $t_1=\prod_{j\in\N}(\alpha_{4j},\alpha_{4j+2})$ and $t_2=[t_1,f]$. Then
\[ t_2=\prod (\alpha_{4j},\alpha_{4j+2})\prod (\alpha_{4j+1},\alpha_{4j+3}).\]
So $t_2$ has cycle type $1^n2^{\aleph_0}$. Again $t$ is conjugate to $t_2$, or is a product of two conjugates.

\medskip

\noindent \textbf{Case 3:} $a_2$ is finite and odd. Using the methods of Cases 1 and 2 as appropriate we find $t_2$ of type $1^n2^{\aleph_0}$
\[ t_2=\prod_{j=0}^\infty (\alpha_{2j},\alpha_{2j+1}).\]
If $t_3=\prod_{j=a_2}^\infty (\alpha_{2j},\alpha_{2j+1})$ then $t_3$ is conjugate to $t_2$ and $t_1t_2$ has cycle type $1^n2^{a_2}$ so it is conjugate to $t$.
\end{proof}

\begin{proof}[Proof of Theorem  \ref{5.1}] Let $f\in S$ and $m=\deg f$. By Lemma  \ref{5.4}, all involutions of degree $\leq m$ lie in
\[ \gen f^{\BS(\Omega,m^+)},\]
and so by Lemma  \ref{5.3} if $g\in\BS(\Omega,m^+)$, then $g=f_1\dots f_8$ where each $f_i$ is conjugate to $1$, $f$ or $f^{-1}$. Hence the theorem. (Normal closures of elements of finite support work similarly -- $\Alt(\Omega)$ is a union of a chain of simple groups so is simple.)
\end{proof}

\chapter*{References for Lecture 5}

Theorem 5.1 (the normal subgroups of $\Sym(\Omega)$) goes back to

\medskip

\noindent
J.\ Schreier and S.\ Ulam, \lq Ueber die Permutationsgruppe der nat\"urlichen Zahlenfolge\rq, \emph{Studia Math.} 4 (1933), 134-141.

\medskip\noindent
and
\medskip

\noindent
R.\ Baer, \lq Die Kompositionreihe der Gruppe aller einein deutigigen Abbildungen einer unendlicher Menge auf sich\rq,  \emph{Studia Math.} 5 (1934), 15--17.

 \medskip\noindent
See also
\medskip

\noindent
W.\ R.\ Scott, \emph{Group theory} (Prentice-Hall, N.J., 1964), Ch.\ 11 and references cited there (p.\ 318).

\chapter{The Symmetric Groups: A Miscellany of Useful Facts}

\begin{fact}\label{6.1} If $\Omega$ is infinite and $n=|\Omega|$, then $|\Sym(\Omega)|=2^n$.
\end{fact}
\begin{proof} Certainly $|\Sym(\Omega)|\leq n^n\leq (2^n)^n=2^{n^2}=2^n$ by the axiom of choice.

Express $\Omega$ as a disjoint union $\bigcup_I \{\alpha_i,\beta_i\}$ for $\alpha_i\neq\beta_i$, where $|I|=n$. If $J\subseteq I$ define
\[ t_J=\prod_{i\in J} (\alpha_i,\beta_i).\]
Clearly $t_J\neq t_{J'}$ if $J\neq J'$; therefore $2^n\leq |\Sym(\Omega)|$.

Hence $2^n=|\Sym(\Omega)|$.
\end{proof}

\noindent
\emph{Idle question:} How far can this be proved in ZF? Is Fact  \ref{6.1} equivalent to AC?\footnote{Without AC it is possible to show that if $|\Omega|\geq 3$ then $|\Omega|<|\Sym(\Omega)|$, see John W.\ Dawson, Jr. and Paul E.\ Howard, \lq Factorials of infinite cardinals\rq, \emph{Fund.\ Math.} 93 (1976), 186--195 and for a short recent proof see Matt Getzen, \lq On the cardinality of infinite symmetric groups\rq,  \emph{Involve} 8 (2015), 749--751.  $\Pi$MN's idle question is mentioned in the paper by Dawson and Howard.  The question is answered negatively in David Pincus, \lq A note on the cardinal factorial\rq, \emph{Fund.\ Math.} 98 (1978), 21--24. }

\begin{fact} The symmetric group $S=\Sym(\Omega)$ carries a natural topology, which is Hausdorff.

Take as a base for the open sets containing $1$ the subgroups $S_{(\Phi)}$ (the pointwise stabilizers of $\Phi\subseteq S$) as $\Phi$ ranges over the finite subsets of $\Omega$.

Therefore the set
\[ \{S_{(\Phi)} f\mid \Phi\subseteq\Omega,\;\Phi\text{ finite}, f\in S\}\]
is a base for the topology. [Note: $S_{(\Phi)}f\cap S_{(\Phi')}g=S_{(\Phi\cap\Phi')}h$.]   Check that this defines a Hausdorff topology on $S$.

Each set $S_{(\Phi)}f$ is open and closed, since
\[ S_{(\Phi)}f=S-\bigcup \{S_{(\Phi)}g\mid g\not\in S_{(\Phi)}f\}.\]
So $S$ is totally disconnected, i.e., zero-dimensional.
\end{fact}

There are variations on this theme: let $\ms F\subseteq\mc P(\Omega)$ such that
\begin{enumerate}
\item $(\forall \Phi_1,\Phi_2\in\ms F)(\exists \Phi_3\in\ms F)(\Phi_1\cup\Phi_2\subseteq\Phi_3)$,
\item $\bigcup_{\ms F} \Phi=\Omega$.
\end{enumerate}
We can take the sets $\{S_{(\Phi)}f\mid \Phi\in\ms F,f\in S\}$ as a base for a topology on $S$.

E.g., $\Omega$ may have a topology and $\ms F$ the family of closed sets, or $\ms F$ the family of compact sets.

If $\Omega$ is countable, then $S$ (with the topology defined by stabilizers of finite sets) is a complete metric space. Choose an enumeration $\Omega=\{\omega_1,\omega_2,\omega_3,\dots\}$ ($\omega_i\neq \omega_j$ for $i\neq j$). Define $d(f,g)=2^{-n}$, where $n=\min\{i\mid \omega_i f\neq \omega_i g\text{ or }\omega_if^{-1}\neq \omega_ig^{-1}\}$.
Check that $d$ is a metric on $S$, that $d$ defines the topology mentioned above, that $(S,d)$ is a complete metric space.

\newpage

\chapter{Symmetric Groups: Some Free Subgroups}
\setcounter{note}{0}

\begin{thm}\label{7.1}Let $\Omega$ be an infinite set, $n=|\Omega|$. There is a free subgroup of rank $2^n$ in $\Sym(\Omega)$.
\end{thm}
\begin{proof} Let $\Omega=\bdisun_I \Omega_i$, where $|\Omega_i|=\aleph_0$ and $|I|=n$.

There is a free subgroup $F_0$ of rank $\aleph_0$ in $\Sym(\Omega_i)$. Hence $F_0^I\leq \prod_I \Sym(\Omega_i)\leq \Sym(\Omega)$. Therefore it is sufficient to prove that $F_0^I$ has a free subgroup of rank $2^n$.

\begin{lem}[Combinatorial lemma]\label{7.2} Let $I,A$ be sets such that $|A|=\aleph_0\leq n=|I|$. There exists $B\subseteq A^I$ such that
\begin{enumerate}
\item $|B|=2^n$,
\item if $f_1,\dots,f_r$ is any finite set of distinct members of $B$ then there exists $i\in I$ such that $f_1(i),\dots,f_r(i)$ are all different.
\end{enumerate}
\end{lem}
\begin{proof} Identify $I$ with the free semigroup (i.e., the set of words) on an alphabet $I_0$ of size $n$; identify $A$ with the free semigroup on a countably infinite set $A_0$. Let
\[ B=\{f:I\to A\mid f\text{ is a semigroup homomorphism induced by a map $I_0\to A_0$}\}.\]
Then $|B|=|A_0^{I_0}|=\aleph_0^n=2^n$.
Given that $f_1,\dots,f_r$ are distinct members of $B$, if $1\leq p<q\leq r$ there exists $x_{pq}\in I_0$ such that $f_p(x_{pq})\neq f_q(x_{pq})$.

Let $w=\prod_{p,q} x_{pq}=x_{12}\ldots x_{1r}x_{23}\ldots x_{r-1}x_r$.  Then $w\in I$ and $f_p(w)\neq f_q(w)$ because these two words have the same length and differ at least in the $x_{pq}$ letter. Thus $f_1(w),\dots,f_r(w)$ are distinct.
\end{proof}

Now complete the proof of  \ref{7.1}. Let $F_0=\gen A$, where $|A|=\aleph_0$. Inside $F_0^I$ let $B$ be a subset of $A^I$ as given by  \ref{7.2}. Claim that $B$ is a free generating set for the group $F$ that it generates.

Let $w=f_1^{\ep_1}f_2^{\ep_2}\ldots f_\ell^{\ep_\ell}$ be a reduced word in members of $B$: so $f_1,\dots,f_\ell\in B$, $\ep_i=\pm 1$, and if $f_i=f_{i+1}$ then $\ep_i=\ep_{i+1}$. Need to prove that $w$ is not the identity on $F$. Choose $i\in I$ such that $f_{j_1}(i)=f_{j_2}(i)\Leftrightarrow f_{j_1}=f_{j_2}$. Viewing $w$ as a function, yields
\[ w(i)=f_1(i)^{\ep_1}\ldots f_\ell(i)^{\ep_\ell}\]
is reduced as written in the free generators $A$ of $F_0$. Therefore, since $F_0$ is free, $w(i)\neq 1$. So $w\neq 1$. Hence $B$ is a free basis for $F$.
\end{proof}

\noindent \textbf{Commentary}:

\begin{note} This free group $F$ has $n$ orbits, each of length $\aleph_0$. This particular free group is highly non-transitive.
\end{note}

\begin{cor} There is a free subgroup of $\Sym(\Omega)$ of rank $2^n$, which is \emph{highly-transitive}.
\end{cor}
\begin{proof} Choose a surjection $B\mathop{\to}\limits^\sigma \bigcup_{k\in\N} \Omega^{(k)}\times \Omega^{(k)}$. (
Possible since $|B|=2^n$ while $|\bigcup \Omega^{(k)}\times \Omega^{(k)}|=n$.) For each $f\in B$ choose a finitary permutation $f^*$ mapping $\boldsymbol\alpha f$ to $\boldsymbol\beta$ where $(\boldsymbol\alpha, \boldsymbol\beta)\in\Omega^{(k)}\times \Omega^{(k)}$ is the image $\sigma f=(\alpha,\beta)$. Then
\[ ff^*:\alpha\mapsto\beta.\]
Let $B^*=\{ff^*\mid f\in B\}$, $F^*=\gen{B^*}$. Then if $k\in\N$ and $\boldsymbol\alpha, \boldsymbol\beta\in\Omega^{(k)}$ then $(\boldsymbol\alpha,\boldsymbol\beta)=\sigma(f)$ for some $f\in B$ (because $\sigma$ is surjective).  And then $\boldsymbol\alpha ff^*=\boldsymbol\beta$. Hence $F^*$ is highly transitive.

Let $w(\phi_1,\dots,\phi_n)=\phi_1^{\ep_1}\ldots \phi_n^{\ep_n}$ be a non-trivial reduced word. Then $w(f_1f_1^*,\dots,f_nf_n^*)=w(f_1,\dots,f_n)\cdot g^*$ where $g^*\in\FS(\Omega)$. Since $f_1,\dots,f_n\in B$ which is a free basis for $F$, we have $w(f_1,\dots,f_n)\neq 1$ in $F$, and since $F$ is torsion free, $w(f_1,\dots,f_n)\not\in \FS(\Omega)$. Hence $w(f_1f_1^*,\dots,f_nf_n^*)=w(f_1,\dots,f_n)\cdot g^*\neq 1$. Thus $B^*$ is a free basis of $F^*$. (Trick due to Tits.)
\end{proof}

\begin{variant} In fact, the same construction allows us to make $F$ $m$-transitive for all $m<n$.
\end{variant}

\begin{note} Free groups are known to be residually finite. It follows that $F_0$, hence $F$ in Theorem  \ref{7.1}, can be chosen to have all orbits finite. If $n=\aleph_0$, then we can also choose $F$ to be cofinitary.
\end{note}

\chapter{Symmetric Groups: Some Lemmas}

A subset $\Sigma$ of $\Omega$ such that $|\Sigma|=|\Omega-\Sigma|=|\Omega|=n$ will be called a \emph{moiety} of $\Omega$. Also, if $\Gamma\subseteq\Delta$ we identify $S_{(\Omega-\Gamma)}$ with $\Sym(\Gamma)$, where $S=\Sym(\Omega)$.

\begin{lem}\label{8.1} Suppose that $\Gamma_1,\Gamma_2\subseteq\Omega$ and $|\Gamma_1\cap\Gamma_2|=|\Gamma_1\cup\Gamma_2|\;(=|\Gamma_1|=|\Gamma_2|)$. Then $\gen{\Sym(\Gamma_1),\Sym(\Gamma_2)}=\Sym(\Gamma_1\cup\Gamma_2)$.
\end{lem}
\begin{proof} We may assume that $\Gamma_1,\Gamma_2$ infinite. Clearly
\[ \gen{\Sym(\Gamma_1), \Sym(\Gamma_2)}\leq \Sym(\Gamma_1\cup\Gamma_2).\]
We now prove the reverse inclusion. Let $g\in\Sym(\Gamma_1\cup\Gamma_2)$. Choose a moiety $\Phi$ of $\Gamma_1\cap\Gamma_2$. Then
\[ |\Phi g\cap \Gamma_1|=|\Gamma_1|\text{ or }|\Phi g\cap \Gamma_2|=|\Gamma_2|.\]

Suppose without loss of generality that the first equality holds.  Then there is a moiety $\Psi$ of $\Gamma_1\cap\Gamma_2$ such that $\Psi g$ is a moiety of $\Gamma_1$.  There is a permutation $h_1\in\Sym(\Gamma_1)$ such that $h_1\ha\Psi=g\ha\Psi$. Then $\Psi\subseteq\Fix(gh_1^{-1})$. But $|\Gamma_1-\Gamma_2|\leq |\Psi|$ and so there exists $h_2\in\Sym(\Gamma_1)$ such that $\Gamma_1-\Gamma_2\subseteq\Psi h_2$. Then $\Gamma_1-\Gamma_2\subseteq\Fix(h_2^{-1}gh_1^{-1}h_2)$. Thus $h_3=h_2^{-1}gh_1^{-1}h_2\in\Sym(\Gamma_2)$, and $g=h_2h_3h_2^{-1}h_1$ is an element of $\gen{\Sym(\Gamma_1),\Sym(\Gamma_2)}$, and the result follows.
\end{proof}

\begin{lem}\label{8.2} Suppose that $G\leq S=\Sym(\Omega)$ and there exists a moiety $\Sigma$ such that $G^\Sigma=\Sym(\Sigma)$ (here $G^\Sigma=G_{\{\Sigma\}}^\Sigma$). Then $(\exists x\in S)(S=\gen{G,x})$.
\end{lem}
\begin{proof} Let $\Sigma_0=\Omega-\Sigma$, and let $\Sigma_1\dot\cup \Sigma_2\dot\cup \Sigma_3\dot\cup \Sigma_4=\Sigma$ be a partition of $\Sigma$ as a union of four moieties. Choose bijections $x_0:\Sigma_0\to\Sigma_1$, $x_1:\Sigma_1\to\Sigma_2$, $x_3:\Sigma_3\to\Sigma_4$.

Define $x\in S$ by
\begin{align*} x\ha \Sigma_0&=x_0,
\\x\ha \Sigma_1&=x_1,
\\x\ha \Sigma_2&=x_1^{-1}x_0^{-1},
\\x\ha \Sigma_3&=x_3,
\\x\ha \Sigma_4&=x_3^{-1}.
\end{align*}
So $x$ permutes $\Sigma_0,\Sigma_1,\Sigma_2$ cyclically and is a product of $3$-cycles there   and $x$ transposes $\Sigma_3,\Sigma_4$ (and is a product of $2$-cycles). [So cycle type of $x$ is $2^n3^n$.]

Consider $\gen{G_{\{\Sigma\}},x^3}$. Since $x^3$ fixes $\Sigma_0,\Sigma_1,\Sigma_2$ pointwise and transposes $\Sigma_3,\Sigma_4$ it lies in $\Sym(\Sigma)$ and has cycle type $1^n2^n$ there. Since $G_{\{\Sigma\}}$ induces $\Sym(\Sigma)$ all conjugates of $x^3$ in $\Sym(\Sigma)$ lie in $\gen{G_{\{\Sigma\}},x^3}$. Hence $\Sym(\Sigma)\leq \gen{G_{\{\Sigma\}},x^3}\leq \gen{G_{\{\Sigma\}},x}$.

But $\gen{G_{\{\Sigma\}},x}$ also contains $x^{-1}\Sym(\Sigma)x=\Sym(\Sigma x)=\Sym(\Omega-\Sigma_1)$. Thus \[\gen{G_{\{\Sigma\}},x}\geq \gen{\Sym(\Sigma),\Sym(\Omega-\Sigma_1)}=\Sym(\Omega)\] by Lemma~ \ref{8.1}.
\end{proof}

\begin{lem}\label{8.3} Let $S=\Sym(\Omega)$, $|\Omega|=n\geq \aleph_0$ and let $(X_i)_{i\in I}$ be a chain of proper subgroups of $S$, such that $S=\bigcup X_i$. Then $|I|>n$.
\end{lem}
\begin{proof} Suppose first that there exist $i\in I$ and $\Sigma$ a moiety of $\Omega$ such that $X_i^\Sigma=\Sym(\Sigma)$. By Lemma  \ref{8.2} there exists $x$ such that $\gen{X_i,x}=S$. But $x\in X_j$ for some $j$ and either $X_i\leq X_j$ or $X_j\leq X_i$.  If $X_i\leq X_j$ then $X_j=S$ and if $X_j\leq X_i$ then $X_i=S$. This contradicts the fact that $X_i,X_j$ are proper subgroups. Thus for all $i$ and for all moieties $\Sigma$, $X_i^\Sigma\neq \Sym(\Sigma)$. Suppose $|I|\leq n$. Partition $\Omega$ as $\bdisun_{i\in I} \Sigma_i$ where $\Sigma_i$ are moieties. Choose $z_i\in \Sym(\Sigma_i)-X_i^{\Sigma_i}$. Now let $z=\prod_I z_i$. Then $z\not\in X_i$ because $z\ha \Sigma_i\not\in X_i^{\Sigma_i}$. Hence $z\not\in \bigcup X_i$, contradicting $|I|>n$.
\end{proof}

\noindent
\textbf{Commentary.} Compare with K\"onig's lemma: If $a_i<b_i$ then $\sigma a_i<\prod b_i$.

\begin{cor} Suppose that $S=\gen{G\cup A}$ where $|A|\leq n$. Then there exists a finite subset $A_0$ of $A$ such that $S=\gen{G\cup A_0}$.
\end{cor}
\begin{proof} Choose a subset $A_0$ of $A$ of least cardinality such that $S=\gen{G,A_0}$. If $A_0$ is infinite, then list it as
\[ A_0=\{a_\xi\mid \xi<\lambda\}\]
for some initial cardinal $\lambda$. Let $X_\eta=\gen{G\cup \{a_\xi\mid \xi<\eta\}}$ for $\eta<\lambda$. Then $S=\bigcup X_\eta$ and $|\lambda|\leq n$, and so by  \ref{8.3}, $X_\eta=S$ for some $\eta$. But $X_\eta=\gen{G\cup A_1}$ where $|A_1|=|\eta|<|\lambda|$. This contradicts choice of $A_0$. Therefore $A_0$ is finite.
\end{proof}

\chapter{Symmetric Groups: Some Problems}

\begin{idleq}Suppose that $S=\gen{G\cup A}$ where $A$ is finite. Does it follow that $S=\gen{G\cup\{x\}}$ for some $x\in S$? Recall, $S$ is $\Sym(\Omega)$ and $n=|\Omega|$.\footnote{This idle question is stated as Question 3.2 in the paper H.\ D.\ Macpherson and Peter M.\ Neumann, \lq Subgroups of infinite symmetric groups\rq, 
\emph{J.\ London Math.\ Soc. (2)}  42 (1990),  64--84 and it is answered in Fred Galvin, \lq Generating countable sets of permutations\rq,
\emph{J.\ London Math.\ Soc.\ (2)} 51 (1995),  230--242.} 
\end{idleq}

\begin{prob} Let $\ell=\ell(n)$ be the least cardinal such that $S$ is expressible as the union of a chain of $\ell$ proper subgroups. What is $\ell$?\footnote{For recent work on this question see e.g.\ Heike Mildenberger and Saharon Shelah, \lq
The minimal cofinality of an ultrapower of $\omega$ and the cofinality of the symmetric group can be larger than $\mathfrak{b}^+$\rq, 
\emph{J.\ Symbolic Logic} 76 (2011), 1322--1340.} 
\end{prob}

\noindent
\emph{Note:} If $G$ is any uncountable group, say of cardinality $m$, then $G$ can be expressed as the union of a chain of length $\cof(m)$ of subgroups of cardinality less than $m$.

\begin{proof} Let $\mu$ be the initial ordinal such that $|\mu|=m$. Write $G=\{g_\xi\mid \xi<\mu\}$. Let $G_\eta=\{g_\xi\mid \xi<\eta\}$. Then $|G_\eta|=|\eta|<m$ (as long as $\eta$ is infinite) and $\bigcup G_\eta=G$. If $k=\cof(m)$, then $\mu=\bigcup_{\alpha\in\kappa} \eta_\alpha$, where $\kappa$ is the initial ordinal with $|\kappa|=k$. Hence $G=\bigcup_{\alpha\in \kappa} G_{\eta_\alpha}$.
\end{proof}

So in Problem 2 we know that $\ell(n)\leq \cof(2^n)$. On the other hand we know that $n<\ell(n)$ (via Lemma  \ref{8.3}).

\smallskip

\noindent
\emph{Problem} 2*.  What is the least cardinal $\ell^*=\ell^*(n)$ such that $S$ is the union of $\ell^*$ proper subgroups?

\smallskip

We know:
\begin{enumerate}
\item If $n=\aleph_0$ then $\ell^*(n)>\aleph_0$ [Macpherson].
\item If $n>\aleph_0$ and $n$ is regular (i.e., $n=\cof(n)$), or if $n=n^{\aleph_0}$, then $\ell^*(n)\leq n$.
\end{enumerate}

\noindent
\emph{Half-hearted conjecture:} $\ell^*(n)=n$ if $\cof(n)>\aleph_0$.

\smallskip

\noindent
\emph{Problem} 2**.   What is the least cardinal $\ell^{**}$ such that $S$ is a union of $\ell^{**}$ proper cosets?

\smallskip

Clearly $\ell^{**}\leq n$, because if $G=S_\alpha$ (the stabilizer of $\alpha$) then $[S:G]=n$. Presumably $\ell^{**}\not< n$.

\begin{prob} Describe the maximal proper subgroups of $S$.
\end{prob}

\noindent
\emph{Problem} 3*. Is it true that for all $G<S$ there exists $M$ maximal proper in $S$ such that $G\leq M$?\footnote{A negative answer needing extra set theoretic assumptions is given in  James E.~Baumgartner, Saharon Shelah and Simon Thomas,
\lq Maximal subgroups of infinite symmetric groups\rq, \emph{Notre Dame J.\ Formal Logic} 34 (1993), 1--11, and  Simon Thomas, \lq Aspects of infinite symmetric groups\rq, in \emph{Infinite groups and group rings (Tuscaloosa, AL, 1992)}, 139--145, Ser.\ Algebra, 1, World Sci.\ Publ., River Edge, NJ, 1993. }

\smallskip

Note that by Lemma  \ref{8.2}, if $G^\Sigma=\Sym(\Sigma)$ for some moiety $\Sigma$, then $G$ is contained in a maximal proper subgroup. For, if $S=\gen{G,x}$ and $M$ is maximal amongst subgroups containing $G$ and avoiding $x$ (exists by Maximal Property) then $M$ is a maximal subgroup of $S$.


\chapter{Symmetric Groups: Supplements to Normal Subgroups}

Recall that if $S=\Sym(\Omega)$, where $|\Omega|=n\geq \aleph_0$, the normal subgroups are
\[ \{1\}<\Alt(\Omega)<\FS(\Omega)<\BS(\Omega,\aleph_1)<\cdots <\BS(\Omega,k)<\BS(\Omega,n)<S.\]
All factors are simple.

Fix $k$, $\aleph_0\leq k\leq n$, let $B=\BS(\Omega,k)$. We seek $G\leq S$ such that $S=B\cdot G$.

\begin{thm}\label{10.1}  Assume that if $m$ is a cardinal such that $k\leq m\leq n$ then the cofinality of $m$ is at least equal to $k$.\footnote{In the lecture and in Neumann and Macpherson (1990) this theorem is stated without this extra assumption.  It was pointed out by Stephen Bigelow in \lq Supplements of bounded permutation groups\rq, 
\emph{J.\ Symbolic Logic}, 63 (1998), 89--102 that extra set theoretic assumptions were needed in this theorem.  The assumption added here is just what makes the argument given in the lecture work, but Bigelow's paper contains a detailed discussion of what kind of set theoretic conditions are needed for the conclusion in this theorem to hold.  Recall that the cofinality of a cardinal number $m$ is the least cardinal number $l$ such that if $\mu$ is the initial ordinal of cardinality $m$ then every sequence of cardinality strictly smaller than $l$ is bounded.}   
If $G\leq S$ then $S=B\cdot G$ if and only if there exists $\Delta\subseteq\Omega$ such that $|\Delta|<k$ and $G^{(\Omega-\Delta)}=\Sym(\Omega-\Delta)$.
\end{thm}
\begin{proof} Suppose that $\Delta$ exists as in the assertion of the theorem. Let $f\in S$. Since $|\Delta|<k$ there exists $g_1\in B$ such that
\[ g_1 \ha \Delta=f\ha \Delta.\]
Then $fg_1^{-1}$ fixes $\Delta$ pointwise. By assumption there exists $g_2\in G_{\{\Delta\}}$ such that $g_2\ha {(\Omega-\Delta)}=fg_1^{-1}\ha (\Omega-\Delta)$. Then $fg_1^{-1}g_2^{-1}$ fixes $\Omega-\Delta$ pointwise. So if $g_3=fg_1^{-1}g_2^{-1}$ then $g_3\in B$ and $f=g_3g_2g_1\in BGB=GB$. Hence $S=BG$.

Suppose, conversely, $S=BG$.

\medskip

\noindent \textbf{Step 1.} There exists a moiety $\Sigma_0$ of $\Omega$ such that $G^{\Sigma_0}=\Sym(\Sigma_0)$.

\begin{proof}[Proof of Step 1] Let $(\Sigma_i)_{i\in I}$ be a family of $n$ pairwise disjoint moieties of $\Omega$. Suppose by way of contradiction that $G^{\Sigma_i}\neq \Sym(\Sigma_i)$ for all $i\in I$. Choose $z_i\in\Sym(\Sigma_i)-G^{\Sigma_i}$ and let $z=\prod z_i$. Since $S=BG$, there exist $x\in B$, $y\in G$ such that $z=xy$. Now $\deg(x)<k$ and so $\Sigma_i\subseteq\Fix(x)$ for all but at most $\deg(x)$ members of $i$. Choose $i_0\in I$ such that $\Sigma_{i_0}\subseteq\Fix(x)$. Then $z^{\Sigma_{i_0}}=y^{\Sigma_{i_0}}$ and yet $z_{i_0}=y^{\Sigma_{i_0}}$, contradicting our choice of $z_{i_0}$. Thus there exists $i\in I$ such that $G^{\Sigma_i}=\Sym(\Sigma_i)$.
\end{proof}

\noindent \textbf{Step 2.} There exists a moiety $\Sigma_1$ of $\Omega$ and a subgroup $H\leq G_{\{\Sigma_1\}}$ such that $H^{\Sigma_1}=\Sym({\Sigma_1})$ and $H^{\Omega-\Sigma_1}\leq \BS(\Omega-\Sigma_1,k)$.

\begin{proof}[Proof of Step 2] Let $\Sigma_0$ be as in Step 1. Let $t\in \Sym(\Sigma_0)\leq S$, $t$ of cycle type $2^n$ on $\Sigma_0$. Since $S=BG$, there exist $x\in B$, $y\in G$ such that $t=xy$. Let $\Gamma_0=\supp(x)\cup\supp(x)\cdot t$. Then $|\Gamma_0|<k$, and $\Gamma_0t=\Gamma_0$. Define $\Sigma_1=\Sigma_0-\Gamma_0$. Then $\Sigma_1$ is a moiety and $G^{\Sigma_1}=\Sym(\Sigma_1)$. Also $t$ and  $x$ fix $\Sigma_1$ setwise, so $y\in G_{\{\Sigma_1\}}$.   Furthermore $x$ fixes
$\Sigma_1$ pointwise   and so $y\ha {\Sigma_1}=t\ha {\Sigma_1}$. Whence $y\ha \Sigma_1$ is of cycle type $2^n$, so the conjugates of $y$ under $\Sym(\Sigma_1)$ generate $\Sym(\Sigma_1)$. Thus if $H=\gen{y}^{G_{\{\Sigma_1\}}}$, then $H^{\Sigma_1}=\Sym(\Sigma_1)$; but since $H^{(\Omega-\Sigma_1)}$ is generated by conjugates in $G_{\{\Omega-\Sigma_1\}}$ of $y\ha {(\Omega-\Sigma_1)}$, $H^{(\Omega-\Sigma_1)}\leq \BS(\Omega-\Sigma_1,k)$.
\end{proof}

\noindent \textbf{Step 3.} There exists a moiety $\Sigma_1$ of $\Omega$ and a subgroup $K\leq G_{\{\Sigma_1\}}$ such that $K^{\Sigma_1}=\Sym(\Sigma_1)$, and $|(\Omega-\Sigma_1)\cap \supp(K)|<k$.
\begin{proof}[Proof of Step 3] Consider subsets $\Gamma$ of $\Omega-\Sigma_1$ ($\Sigma_1$ as given by Step 2) such that there exists $K\leq G_{\{\Sigma_1\}}$ with $K^{\Sigma_1}=\Sym(\Sigma_1)$, and $\supp(K)\subseteq\Sigma_1\cup\Gamma$ \emph{and} $K^\Gamma\leq \BS(\Gamma,k)$. Choose such $\Gamma$ of least cardinality $m$. Suppose by way of contradiction that $m\geq k$. Take $\mu$ to be the initial ordinal for the cardinal $m$, and express $\Gamma$ as a union $\bigcup \Gamma_\xi$ (for $\xi<\mu$) where $|\Gamma_\xi|<m$, $\Gamma_{\xi_1}\subseteq \Gamma_{\xi_2}$ if $\xi_1\leq \xi_2$.

Define $K_\xi=\{f\in K\mid \supp(f)\subseteq \Sigma_1\cup\Gamma_\xi\}$. Then $(K_\xi)_{\xi<\mu}$ is a chain of subgroups of $K$. Also, $K=\bigcup K_\xi$.\footnote{Here we use the assumption that the cofinality of $m$ is at least equal to $k$.  That assumption implies that any subset of $\Gamma$ of cardinality less than $k$ is contained in one of the sets $\Gamma_\xi$.  Thus, if $g\in K$ then there exists $\xi$ such that $(\Omega-\Sigma_1)\cap \supp(K)\subseteq \Gamma_\xi$ and $g\in K_\xi$.} Therefore $K^{\Sigma_1}=\bigcup_{\xi<\mu} K_\xi^{\Sigma_1}$. By Lemma  \ref{8.3} there exists $\eta<\mu$ such that
\[ K_\eta^{\Sigma_1}=\Sym(\Sigma_1),\]
and thus $\Gamma_\eta$ contradicts the minimality of $m$. Hence $m<k$.
\end{proof}

\begin{definition} We say that a moiety $\Sigma$ is \emph{good} if it has the property determined in Step 3. (Thus $\Sigma_1$ is good.)
\end{definition}

\noindent \textbf{Step 4.} If $\Sigma'$ is \emph{any} moiety then there exists $\Sigma$ such that $\Sigma$ is good, $\Sigma\subseteq \Sigma'$, and $|\Sigma'-\Sigma|<k$.
\begin{proof}[Proof of Step 4] Choose $f\in S$ such that $\Sigma_1f=\Sigma'$ (within the symmetric group we can move any moiety to any other moiety). Then $f=xy$ for suitable $x\in B$, $y\in G$. If $\Sigma=\Sigma'-\supp(x)$, then $\Sigma=\Sigma''y$ for $\Sigma''\subseteq \Sigma$ and $|\Sigma_1-\Sigma''|<k$. Then $y^{-1}K_{\{\Sigma'\}}y$ is a suitable `$K$' for $\Sigma$.
\end{proof}

\noindent\emph{Finish the proof of the theorem.} Write $\Omega=\Sigma_1'\cup\Sigma_2'$ where $|\Sigma_1'\cap \Sigma_2'|=n$.

Choose $\Sigma_1''\subseteq \Sigma_1'$, $\Sigma_2''\subseteq \Sigma_2'$ such that $|\Sigma_1'-\Sigma_1''|<k$, $|\Sigma_2'-\Sigma_2''|<k$, and there are groups $K_1,K_2$ as in Step 3. Choose
\begin{align*} \Sigma_1&=\Sigma_1'-\left(\supp(K_2)\cap (\Omega-\Sigma_2)\right)
\\ \Sigma_2&=\Sigma_2'-\left(\supp(K_1)\cap (\Omega-\Sigma_1)\right).
\end{align*}
There are subgroups $L_1,L_2$ such that
\begin{itemize}
\item $L_1$ fixes $\Sigma_2-\Sigma_1$ pointwise and $L_1^{\Sigma_1}=\Sym(\Sigma_1)$,
\item $L_2$ fixes $\Sigma_1-\Sigma_2$ pointwise and $L_2^{\Sigma_2}=\Sym(\Sigma_2)$.
\end{itemize}

\end{proof}

Consider the proposition
\[ \mathrm{SP}:\; S=BG\Leftrightarrow (\exists \Delta\subseteq \Omega)(|\Delta|<k\text{ and }S_{(\Delta)}\leq G).\]
Semmes (1983) proved GCH $\Rightarrow$ SP.

\begin{thm}\label{10.2}\footnote{This is proved in Bigelow's paper referred to above.  The proof given in the lecture was based on Theorem 10.1 without the cofinality condition and is thus not valid.} $\mathrm{SP}_{k,n}\Leftrightarrow (m<k\Rightarrow 2^m<2^n)$.
\end{thm}

Therefore SP is universally true if and only if exponentiation is strictly monotone. (Note, GCH implies $x<y\Rightarrow 2^x<2^y$.)




\chapter{Symmetric Groups: Subgroups of Small Index I}

\begin{thm}[{{\hyperref[ref:lec10-13]{John D.\ Dixon, $\Pi$MN, Simon Thomas, 1986; Semmes, 1982}}}]\footnote{This theorem led to a large number of results characterizing \lq subgroups of small index\rq\ in automorphism groups of a wide variety of structures.  One such result is discussed in Lecture 16.  For more information the reader can start by consulting Wilfrid Hodges, Ian Hodkinson, Daniel Lascar and Saharon Shelah, \lq The small index property for $\omega$-stable $\omega$-categorical structures and for the random graph\rq, 
\emph{J.\ London Math.\ Soc.\ (2)} 48 (1993), 204--218.} 
\label{11.1} Let $\Omega$ be a countably infinite set, $S=\Sym(\Omega)$, $G\leq S$. Then $|S:G|<2^{\aleph_0}$ if and only if there exists a finite set $\Delta_0$ such that $S_{(\Delta_0)}\leq G\leq S_{\{\Delta_0\}}$. (It follows that there are only countably many subgroups of index $<2^{\aleph_0}$ and apart from $S$ itself they all have index $\aleph_0$.)
\end{thm}
\noindent\emph{Proof.} Certainly if $S_{(\Delta_0)}\leq G$ then 
\[|S:G|\leq |S:S_{(\Delta_0)}|=|\{f:\Delta_0\to \Omega\mid f\text{inj.\ map}\}|\leq \aleph_0^{|\Delta_0|}=\aleph_0.\]

Now suppose that $|S:G|<2^{\aleph_0}$.

\begin{lem}\label{11.2} There exists a moiety $\Sigma$ such that $\Sym(\Sigma)\leq G$.
\end{lem}\begin{proof} Let $(\Sigma_i)_{i\in I}$ be a family of $\aleph_0$ pairwise disjoint moieties of $\Omega$. Let $S_i=\Sym(\Sigma_i)$, and let $P=\prod_I S_i\leq S$. Let $H=P\cap G$. Then $|P:H|\leq |S:G|<2^{\aleph_0}$. Let $H_i=H^{\Sigma_i}\leq S_i$. But $H\leq \prod H_i$ and so $|P:H|\geq |P:\prod H_i|=\prod_i |S_i:H_i|$. Therefore $H_i=S_i$ for all except finitely many $i$. Then $G^{\Sigma_i}=\Sym(\Sigma_i)$ for infinitely many $i$. Choose $\Sigma$ to be any moiety such that $G_{\{\Sigma\}}^\Sigma=\Sym(\Sigma)$. Let $K=G_{(\Omega-\Sigma)}=G\cap \Sym(\Sigma)$. Then $K$ is normal in $\Sym(\Sigma)$ and so $K^\Sigma\normal G^\Sigma=\Sym(\Sigma)$. But $|\Sym(\Sigma):K|<2^{\aleph_0}$. Hence $K^\Sigma=\Sym(\Sigma)$. Thus $K=\Sym(\Sigma)$, so $\Sym(\Sigma)\leq G$.
\end{proof}

\begin{lem}[{{\hyperref[ref:lec10-13]{Sierpinski, 1928}}}] \label{11.3} Let $\Gamma$ be any countably infinite set. There exists a family $(\Gamma_i)_{i\in I}$ of moieties such that
\begin{enumerate}
\item $\Gamma_i\cap\Gamma_j$ is finite for $i\neq j$,
\item $|I|=2^{\aleph_0}$.
\end{enumerate}
[``Almost disjoint family of sets.'']
\end{lem}\begin{proof} The lemma is invariant under bijections, so take $\Gamma=\Q$, $I=\R$, and take $\Gamma_i$ to be the set of points in a strictly monotone increasing sequence converging to $i$.
\end{proof}

\begin{proof}[Proof of Theorem  \ref{11.1} continued] 
We now continue the proof of Theorem  \ref{11.1}. Let $\Sigma$ be as in Lemma  \ref{11.2}. Let $(\Gamma_i)_{i\in I}$ be a pairwise disjoint family of moieties of $\Sigma$, as given in  \ref{11.3}. Choose $x_i\in S$ such that $x_i$ interchanges $\Omega-\Sigma$ with $\Gamma_i$ and fixes $\Sigma-\Gamma_i$ pointwise.

Since there are $2^{\aleph_0}$ permutations $x_i$ and $<2^{\aleph_0}$ cosets of $G$ in $S$, there exist $i,j$, $i\neq j$, such that $g=x_ix_j^{-1}\in G$. Let $\Sigma'=\Sigma g$.

\medskip

\noindent\textbf{Claim:} $\Sigma\cap\Sigma'$ is a moiety and $\Omega-(\Sigma'\cup\Sigma)$ is finite.

\noindent
We have
\[ \Omega-(\Sigma\cup \Sigma')=(\Omega-\Sigma)\cap(\Omega-\Sigma')
=(\Omega-\Sigma)\cap (\Omega-\Sigma)g
=(\Gamma_i\cap\Gamma_j')x_j^{-1},
\]
which is finite. Also,
\[ \Sigma\cap\Sigma'=\Sigma\cap \Sigma x_ix_j^{-1}
\supseteq \Sigma \cap (\Omega-\Gamma_i)x_j^{-1}
\supseteq \Gamma_j\quad \text{(a moiety)}.
\]

By Lemma \ref{8.1}, $\gen{\Sym(\Sigma),\Sym(\Sigma')}=\Sym(\Sigma\cup\Sigma')$. But \[\gen{\Sym(\Sigma),\Sym(\Sigma')}=\gen{\Sym(\Sigma),g^{-1}\Sym(\Sigma)g}\leq G.\] Thus there exists $\Delta\subseteq\Omega$, $\Delta$ finite, $S_{(\Delta)}\leq G$.

Now choose $\Delta_0$ a smallest finite set such that $S_{(\Delta_0)}\leq G$.
If $g\in G$ and $\Delta_0'=\Delta_0g$, then
\[ G\geq \gen{\Sym(\Omega-\Delta_0),\Sym(\Omega-\Delta_0')}=\Sym(\Omega-(\Delta_0\cap\Delta_0')).\]
By minimality $\Delta_0'=\Delta_0$. Hence $G\leq S_{\{\Delta_0\}}$.

\end{proof}

\chapter{Symmetric Groups: Subgroups of Small Index II}

\begin{question} What happens if $|\Omega|=n>\aleph_0$?
\end{question}

\begin{thm}[{{\hyperref[ref:lec10-13]{John Dixon, $\Pi$MN and Simon Thomas, 1986}}}]\label{12.1} Suppose that $G\leq S=\Sym(\Omega)$ and that $|S:G|\leq n$. Then there exists $\Delta\subseteq \Omega$ such that $|\Delta|<n$ and $S_{(\Delta)}\leq G$. (Converse fails with $n=\aleph_\omega$ and $\Delta$ countable, then $|S:S_{(\Delta)}|>n$.)
\end{thm}

The first step in the proof is:

\begin{lem}\label{12.2} There exists a moiety $\Sigma$ such that $\Sym(\Sigma)\leq G$. [Requires only $|S:G|<2^n$.]
\end{lem}

\begin{lem}[{{\hyperref[ref:lec10-13]{Sierpinski, 1928}}}] If $|\Gamma|=n$ then there exists a family $(\Gamma_i)_{i\in I}$ of moieties of $\Gamma$ such that $|\Gamma_i\cap\Gamma_j|<n$ when $i\neq j$ and $|I|>n$.
\end{lem}
\begin{proof} Let $\nu$ be the initial ordinal of cardinality $n$. Let's call a set $X$ of functions $\nu\to\nu$ \emph{good} if $f,g\in X$,
\[ f\neq g\Rightarrow |\{\xi\in \nu\mid f(\xi)=g(\xi)\}|<n.\]
The union of a chain of good sets of functions is itself a good set. By maximality principles, there exists a maximal good set $X_0$.

Suppose by way of contradiction that $|X_0|=m\leq n$. Index $X_0$ as $\{f_\eta\mid \eta<\mu\}$ where $\mu$ is the initial ordinal of cardinality $m$. Define
\[ \Phi_\xi=\{f_\eta(\xi)\mid \eta<\xi\}\]
for $\xi<\nu$.

Then $|\Phi_\xi|\leq |\xi|<n$. Now define $g:\nu\to\nu$ by choosing $g(\xi)\in \nu-\Phi_\xi$. Then $X_0\cup\{g\}$ is good. If $f\in X_0$, then $f=f_\eta$ for some $\eta$ and $g(\xi)\neq f(\xi)$ as soon as $\xi>\eta$.

This contradicts maximality. Hence $|X_0|>n$.

Now choose a bijection $\phi:\nu\times \nu\to \Gamma$, and let $\Gamma_i=\Psi_i\phi$ where $\Psi_i$ is the graph of $i\in X_0$, $\{(\xi,i(\xi)\mid \xi\in\nu\}$. Then $|\Gamma_i\cap\Gamma_j|<n$ when $i\neq j$.
\end{proof}

The proof of  \ref{12.1} now follows exactly as in the countable case.

\begin{question}
\begin{enumerate}
\item Can the (very strong) hypothesis $|S:G|\leq n$ be weakened to $|S:G|<2^n$?
\item Can the conclusion be strengthened to give $\Delta\subseteq \Omega$ such that $S_{(\Delta)}\leq G$ and $|\Delta|<\cof(n)$? (Reason: If $|\Delta|\geq \cof(n)$ then $|S:S_{(\Delta)}|>n$.)
\end{enumerate}
\end{question}

\smallskip

If $\ms F\subseteq \mc P(\Omega)$, define $\ms F$ to be a \emph{filter} if $\ms F\neq \emptyset$ and
\begin{enumerate}
\item $\Gamma_1,\Gamma_2\in\ms F\Rightarrow \Gamma_1\cap\Gamma_2\in\ms F$,
\item $\Gamma_1\in\ms F$ and $\Gamma_1\subseteq \Gamma_2\Rightarrow \Gamma_2\in\ms F$.
\end{enumerate}
If $\ms F$ is a filter, let
\[ S_{(\ms F)}=\{f\in S\mid \Fix(f)\in\ms F\},\]
\[ S_{\{\ms F\}}=\{f\in S\mid \Gamma\in\ms F\Rightarrow \Gamma f,\Gamma f^{-1}\in\ms F\}=\Aut(\ms F).\]

\begin{lem} These are groups and $S_{(\ms F)}\normal S_{\{\ms F\}}$.
\end{lem}
\begin{proof} Exercise.\end{proof}\newpage

\begin{thm}\label{12.5} If $G\leq S$ and $|S:G|<2^n$ then there is a filter $\ms G$ such that
\begin{enumerate}
\item $\ms G$ contains moieties of $\Omega$, 
\item $S_{(\ms G)}\leq G\leq S_{\{\ms G\}}$.
\end{enumerate}
\end{thm}
\begin{proof} As in  \ref{12.2}.\end{proof}

\chapter{Symmetric Groups: Subgroups of Small Index III: Commentary}



\setcounter{note}{0}
\begin{note}
Essentially the only candidates for subgroups of small index are the groups $S_{\{\ms F\}}$ for suitable filters $\ms F$.   The following is a version of Theorem \ref{12.5}.
\end{note}

\begin{thm}\label{13.1} If $|S:G|<2^n$ then there is a filter $\ms G\subseteq \mc P(\Omega)$ such that $S_{(\ms G)}\leq G\leq S_{\{\ms G\}}$ and there is a moiety of $\Omega$ in $\ms G$.
\end{thm}
\begin{proof} Define
\[ \ms G=\{\Gamma\subseteq \Omega\mid (\exists \Delta\subseteq \Gamma)(|\Omega-\Delta|=n\text{ and }S_{(\Delta)}\leq G)\}.\]
Certainly there is a moiety in $\ms G$, so $\ms G\neq\emptyset$. Clearly if $\Gamma_2\supseteq \Gamma_1\in\ms G$ then $\Gamma_2\in\ms G$. So all we need to prove is that $\Gamma_1,\Gamma_2\in\ms G$ implies $\Gamma_1\cap\Gamma_2\in\ms G$.

Now suppose $\Delta_1\subseteq \Gamma_1$, $\Delta_2\subseteq \Gamma_2$, and $|\Omega-\Delta_1|=n$ and $|\Omega-\Delta_2|=n$, and $S_{(\Delta_1)}\leq G$ and $S_{(\Delta_2)}\leq G$.
Let $\Phi_i=\Omega-\Delta_i$ for $i=1,2$. If $|\Phi_1\cap\Phi_2|=n$ then \[\gen{\Sym(\Phi_1),\Sym(\Phi_2)}=\Sym(\Phi_1\cup\Phi_2)=S_{(\Omega-(\Phi_1\cup\Phi_2))}=S_{(\Delta_1\cap\Delta_2)}.\] So $\Delta_1\cap\Delta_2$ `witnesses' $\Gamma_1\cap\Gamma_2\in\ms G$.

Suppose then that $|\Phi_1\cap\Phi_2|<n$. Choose a family $(\Sigma_i)_{i\in I}$ of pairwise disjoint moieties such that $\Sigma_i\cap\Phi_1$ and $\Sigma_i\cap \Phi_2$ are moieties for each $i\in I$.

By Lemma  \ref{12.2}, $\Sym(\Sigma_i)\leq G$ for almost all $i$. Therefore there is  a moiety $\Phi\subseteq \Phi_1\cup\Phi_2$ such that $\Phi\cap \Phi_1$, $\Phi\cap\Phi_2$ are moieties, and $\Sym(\Phi)\leq G$.

But now $\Sym(\Phi_1\cup\Phi_2)=\gen{\Sym(\Phi), \Sym(\Phi_1),\Sym(\Phi_2)}\leq G$. Hence if $\Delta=\Delta_1\cap\Delta_2$ then $S_{(\Delta)}\leq G$, so $\Delta$ is a witness to the assertion that $\Gamma_1\cap\Gamma_2\in\ms G$.
\end{proof}

\begin{note} Compare with Theorem  \ref{12.1}: if $|S:G|\leq n$ then $\exists \Delta(|\Delta|<n$ and $S_{(\Delta)}\leq G)$. We would \emph{like} to be able to prove that if $|S:G|<2^n$ then $\exists\Delta$ ($|\Delta|<n$ and $S_{(\Delta)}\leq G$).

\hyperref[ref:lec10-13]{Shelah and Thomas (1989)} have proved: If there exists a measurable cardinal then there is a counterexample. In fact, there is a filter $\ms F$ such that $|S:S_{\{\ms F\}}|<2^n$ and $S_{(\Delta)}\not\leq S_{\{\ms F\}}$ for all $\Delta$ such that $|\Delta|< n$.
\end{note}

\begin{note} We'd also like to prove that if $|S:G|\leq n$ then there exists $\Delta$ such that $|\Delta|<\cof(n)$ and $S_{(\Delta)}\leq G$.

Note that if $k=\cof(n)$ then
\[ n^m\begin{cases} >n&m\geq k,
\\ \leq n&m<k.\end{cases}\]
(Essentially K\"onig's theorem.) Hence $|S:S_{(\Delta)}|\leq n$ if and only if $|\Delta|<k$.  Whence the question.

Shelah and Thomas have also proved that if there exists a measurable cardinal it implies that there is a counterexample in a suitable extension of ZFC+$\exists$MC. So it is not provable in ZFC.
\end{note}

\begin{note} By contrast, in some (wild) versions of set theory, \emph{stronger} theorems are provable. For example, if $2^{\aleph_0}=2^{\aleph_1}=\aleph_2$, $n=\aleph_1$, then $|S:G|<2^{\aleph_1}$ implies $\exists\Delta_0$ ($\Delta_0$ is finite and $S_{(\Delta_0)}\leq G\leq S_{\{\Delta_0\}}$).
[See Dixon, $\Pi$MN and Simon Thomas, 1986.]
\end{note}
\begin{note}
Suppose now that $|\Omega|=n=\aleph_0$. Then Theorem  \ref{11.1} can be expressed in the form $|S:G|<2^{\aleph_0}$ if and only if $G$ is an open subgroup. David Evans has proved that if $G$ and $T$ are closed subgroup of $S$ and $G\leq T$ then $|T:G|<2^{\aleph_0}$ if and only if $G$ is open in $T$ [unpublished]\footnote{See David M.\ Evans, \lq A note on automorphism groups of countably infinite structures\rq, \emph{Arch. Math. (Basel)}, 49 (1987), 479--483.}.
\end{note}

\chapter*{References for Lectures 10, 11, 12, 13}
\label{ref:lec10-13}
\begin{itemize}
\item[1.] Stephen W.\ Semmes, \lq Infinite symmetric groups, maximal subgroups, and filters\rq, Unpublished, but see preliminary report in \emph{Abstracts Amer.\ Math.\ Soc.}, 3 (Jan.\ 1982), 38.
\item[2.] Stephen W.\ Semmes, \lq Endomorphisms of infinite symmetric groups\rq, Unpublished, but see preliminary report in  \emph{Abstracts Amer.\ Math.\ Soc.}, 2 (August 1981), 426.
\item[3.] John D.\ Dixon,  Peter.\ M.\ Neumann and Simon Thomas, \lq Subgroups of small index in infinite symmetric groups\rq, \emph{Bull.\ London Math.\ Soc.}, 18 (1986), 580--586.
\item[4.] Saharon Shelah and Simon Thomas, \lq Subgroups of small index in infinite symmetric groups II\rq, \emph{J.~Symbolic Logic} 54 (1989), 95--99.
\item[5.] W.\ Sierpi\'nski, \lq  Sur une d\'ecomposition des ensembles\rq, \emph{Monatshefte f\"ur Math.\ u.\ Phys.}, 35 (1928), 239--242.
\end{itemize}

\chapter{Automorphisms of the Rational Line I}

\begin{terminology}
A linearly ordered set is \emph{dense} if $\alpha<\beta\Rightarrow \exists \gamma(\alpha<\gamma<\beta)$. It is said to be \emph{open} if it has neither a maximum nor a minimum.
\end{terminology}

\begin{thm}[Cantor] \label{14.1} A countable dense open linearly ordered set is order isomorphic to $(\Q,<)$.
\end{thm}
\begin{proof} Let $(\Lambda,<)$ be a countable dense open linear order. Say
\begin{align*} \Lambda&=\{\lambda_0,\lambda_1,\lambda_2,\dots\},
\\\Q&=\{q_0,q_1,q_2,\dots\}.
\end{align*}
Define $\phi:\Lambda\to\Q$ by specifying 
\begin{enumerate}
\item[(1)] $\lambda_0\phi=q_0$, and
\item[(2)] if $\phi$ is defined, order-preserving and injective on $\{\lambda_0,\lambda_1,\dots,\lambda_r\}$ then: let \[\{\lambda_0,\lambda_1,\dots,\lambda_r\}=\{\mu_0,\mu_1,\dots,\mu_r\}\] where $\mu_0<\mu_1<\cdots<\mu_r$, and if $\lambda_{r+1}$ lies in $(\mu_i,\mu_{i+1})$ (or possibly $(-\infty,\mu_0)$ or $(\mu_r,\infty)$ ) then choose $s$ to be the least natural number such that $q_s\in(\mu_i\phi,\mu_{i+1}\phi)$ (or the other intervals), and define $\lambda_{r+1}\phi=q_s$.
\end{enumerate}
Ultimately $\phi$ is defined, order-preserving, injective, from $\Lambda$ to $\Q$.

The problem is to prove that $\phi$ is surjective. Suppose not: Let $m$ be the least natural number such that $q_m\not\in \mathrm{Im}(\phi)$. Choose $r$ so large that each of $q_0,q_1,\dots,q_{m-1}$ lies in $\{\lambda_0,\lambda_1,\dots,\lambda_r\}\phi$. Let $\{\mu_0,\dots,\mu_r\}=\{\lambda_0,\dots,\lambda_r\}$ be as before. Let $\nu_i=\mu_i\phi$, for $0\leq i\leq r$. Then $\nu_0<\nu_1<\cdots <\nu_r$ and $\nu_i\in\Q$ for each $i$ and $q_m\in (\nu_i,\nu_{i+1})$ for some $i$ (or $q_m\in(-\infty,\nu_0)$ or $q_m\in(\nu_r,\infty)$). Let $t$ be minimal such that $\lambda_t\in (\mu_i,\mu_{i+1})$ (or $\lambda_t\in(-\infty,\mu_0)$ or $\lambda_t\in(\mu_r,\infty)$). Then when $\lambda_t\phi$ comes to be defined it \emph{must} be $q_m$, by the construction.
\end{proof}

\begin{cor}\label{14.2} If $\Lambda$ is any countable linearly ordered set, then $\Lambda$ is order-isomorphic to a subset of $\Q$.
\end{cor}
\begin{thm} \label{14.3} If $k\in\N$ and $\alpha_1<\alpha_2<\cdots<\alpha_k$, $\beta_1<\beta_2<\cdots<\beta_k$, ($\alpha_i,\beta_i\in\Q$) then there exists $g\in\Aut(\Q,\leq)$ such that $\alpha_i g=\beta_i$ for all $1\leq i\leq k$, i.e., the automorphism group $\Aut(\Q,\leq)$ is \emph{order-$k$-transitive} for all $k\in\N$.
\end{thm}
\begin{proof} Define
\[ \omega g=\begin{cases} \omega+(\beta_1-\alpha_1)&\omega\leq \alpha_1,
\\ \beta_i+\frac{\omega-\alpha_i}{\alpha_{i+1}-\alpha_i}(\beta_{i+1}-\beta_i)&\alpha_i\leq \omega\leq \alpha_{i+1},
\\\omega+\beta_k-\alpha_k&\alpha_k\leq \omega.\end{cases}\]
This is the required element.
\end{proof}

\begin{thm}\label{14.4} Let $A=\Aut(\Q,\leq)$, and define
\begin{align*}
B&=\{f\in A\mid\supp(f)\text{ is bounded}\},
\\ R&=\{f\in A\mid\supp(f)\text{ is bounded below}\},
\\ L&=\{f\in A\mid\supp(f)\text{ is bounded above}\}.
\end{align*}
Then:
\begin{enumerate}
\item $B$, $L$ and $R$ are normal subgroups of $A$.
\item $B=L\cap R$ and $A=LR$.
\item $B$, $R/B$ and $L/B$ (hence also $A/L$, $A/R$) are simple.
\end{enumerate}

Thus 
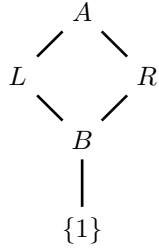
\begin{center}
\begin{tikzpicture}[node distance=1.2cm,line width=1pt]
\node(A) at (0,0)     {$A$};
\node(L)      [below left of=A]  {$L$};
\node(R)      [below right of =A]       {$R$};
\node(B)      [below right of =L]      {$B$};
\node(I)      [below of =B]      {$\{1\}$};

\draw (A) -- (L); 
\draw (A) -- (R); 
\draw (L) -- (B); 
\draw (R) -- (B); 
\draw (B) -- (I); 
\end{tikzpicture}
\end{center}
is a diagram of the lattice of all normal (and all subnormal) subgroups of $A$. (Recall: `subnormal subgroup' means transitive closure of normal subgroup).
\end{thm}

\begin{thm}\label{14.5} Let $\Sigma\subseteq \Q$. Then $\Sigma$ has $<2^{\aleph_0}$ transforms under $A$ if and only if $\Sigma$ is a union of finitely many intervals with rational endpoints (possibly $\pm \infty$ as endpoints, possibly open/closed/half-open, etc).
\end{thm}
\begin{proof} If $\Sigma$ is such a union then $\{\Sigma g\mid g\in A\}$ is countable. So suppose $|\{\Sigma g\mid g\in A\}|<2^{\aleph_0}$. Let $H=A_{\{\Sigma\}}$, so that $|A:H|<2^{\aleph_0}$.

The proof takes three steps:
\begin{itemize}
\item Step 1: All $H$-orbits are convex.
\item Step 2: There are only finitely many $H$-orbits.
\item Step 3: The end-points of the $H$-orbits are rational.
\end{itemize}
\noindent \textbf{Step 1.} Let $\Gamma$ be a non-convex $H$-orbit. There exist $\alpha_0,\alpha_1\in\Gamma$ with $\alpha_0<\alpha_1$, and $\beta_0\not\in\Gamma$ with $\alpha_0<\beta_0<\alpha_1$. Then there exists $h\in H$ such that $\alpha_0h=\alpha_1$. Let $\alpha_i=\alpha_0h^i$. Then
\[ \cdots \alpha_{-1}<\alpha_0<\alpha_1<\alpha_2<\cdots.\]
If $\beta_i=\beta_0h^i$ then $\alpha_i<\beta_i<\alpha_{i+1}$ for each $i$. For any subset $I$ of $\Z$ there exists $x_I\in A$ such that
\begin{align*} \alpha_ix_I&=\begin{cases} \beta_i&i\in I,\\\alpha_i&i\not\in I,\end{cases}
\\ \beta_ix_I&=\beta_i \qquad i\not\in I.
\end{align*}
Now if $I\neq J$ then $x_Ix_J^{-1}\not\in H$ (because $\alpha_ix_Ix_J^{-1}=\beta_i$ for some $i$). Hence, $|A:H|=2^{\aleph_0}$, which is false.

\medskip

\noindent \textbf{Step 2.} Similar (exercise).

\medskip

\noindent \textbf{Step 3.} Easy.
\end{proof}

\chapter{Automorphisms of the Rational Line II}

Notation and terminology:
\begin{itemize}
\item $A=\Aut(\Q,\leq)$.
\item $B=\{f\in A\mid \supp(f)$ is bounded (on both sides)$\}$.
\item If $\alpha,\beta\in \Q$ then $A(\alpha,\beta)=\{f\in A\mid \supp(f)\subseteq (\alpha,\beta)\}$, so $A(\alpha,\beta)=\Aut((\alpha,\beta),\leq)\cong A$ and $B=\bigcup A(\alpha,\beta)$.
\item If $f\in A$ then $f$ extends to $f^*:\R\to\R$.
\item The translation $\xi\mapsto \xi+1$ will be called $a$.
\item If $u\in \R$, $f\in A$ and $\xi f-\xi\mapsto u$ as $|f|\to\infty$ then $f$ will be called a rough translation through $u$.
\item A \emph{$\Z$-sequence} is a sequence $(\xi_n)_{n\in\Z}$ such that $\xi_n<\xi_{n+1}$ for all $n$ and $\xi_n\to\pm\infty$ as $n\to\pm\infty$.
\item A subset $\Gamma$ will be said to of $\Z$-type if it may be indexed as a $\Z$-sequence.
\end{itemize}

\begin{lem}\label{15.1} $A$ acts transitively on the set of $\Z$-sequences.
\end{lem}
\begin{proof}
Let $(\xi_n)$, $(\eta_n)$ be $\Z$-sequences. Let $h_n:[\xi_n,\xi_{n+1})\to[\eta_n,\eta_{n+1})$ be an order-isomorphism and define $h:\Q\to\Q$ by $\xi h=\xi h_n$ if $\xi\in [\xi_n,\xi_{n+1})$. Then $h\in A$ as desired.
\end{proof}

\begin{lem}\label{15.2} The following conditions on $f\in A$ are equivalent:
\begin{enumerate}
\item $f^*$ is fixed-point-free on $\R$;
\item one cycle of $f$ is a $\Z$-sequence;
\item all cycles are $\Z$-sequences;
\item $f$ is conjugate in $A$ to $a$ or $a^{-1}$.
\end{enumerate}
\end{lem}
\begin{proof} 
(iv)$\Rightarrow$(iii)$\Rightarrow$(ii)$\Rightarrow$(i) are trivial. 

(i)$\Rightarrow$(iii) because if $\gen f$ had an orbit that was bounded above (or below) then its supremum (or infimum) in $\R$ would be a fixed point of $f^*$.

For (iii)$\Rightarrow$(iv), suppose that $\gen f$ has an orbit $\Gamma$ which is a $\Z$-sequence. We can suppose that $\Gamma=\{\xi_n\mid n\in\Z\}$ and replacing $f$ by $f^{-1}$ if necessary that $\xi_n f=\xi_{n+1}$ ($\xi_n<\xi_{n+1}$ and $\xi_n\to\pm\infty$ as $n\to\pm\infty$). Let $g_0:[0,1)\to[\xi_0,\xi_1)$ be an order-isomorphism and define $g:\Q\to\Q$ by the rule $\xi g=(\xi-n)g_0 f^n$ if $\xi\in [n,n+1)$. Check that $g^{-1}ag=f$. If $\xi\in[\xi_n,\xi_{n+1})$ then $\xi=(\eta-n)g_0f^n$ for some $\eta$ with $\eta\in[n,n+1)$, so 
\[ \xi g^{-1}ag=\eta ag=(\eta+1)g=(\eta+1-(n+1))g_0f^{n+1}=\xi f.\]
\end{proof}

\begin{lem} \label{15.3} If $f\in A$ there exist conjugates $f_1,f_2$ of $a$ such that $f=f_1f_2^{-1}$.
\end{lem}\begin{proof}
Choose a $\Z$-sequence $(\xi_n)$ such that $\xi_{n+1}>\max\{n+1,\xi_n,\xi_n f,\xi_n f^{-1}\}$ (actually a little redundant, but no matter), if $n\geq 0$. If $n\leq 0$ choose $\xi_{n-1}<\min \{n-1,\xi_n,\xi_n f,\xi_n f^{-1}\}$. Can do this recursively. Then
\[ \xi_{n-1}<\xi_n f^{\pm 1}<\xi_{n+1}.\]
Choose $f_2\in A$ such that $f_2:\xi_n\to\xi_{n+2}$ (Lemma  \ref{15.1}). Define $f_1=ff_2$. Then $f_2$ is conjugate to $a$ by Lemma  \ref{15.2}. Let $\xi\in \R$ and $\xi_n\leq \xi<\xi_{n+1}$. Then $\xi f_1^*=\xi f^*f_2^*$ but $\xi_{n+1}\leq \xi f^*<\xi_{n+2}$, and so $\xi_{n+3}\leq \xi f^*f_2^*<\xi_{n+4}$. Hence $\xi f^*\neq \xi$. So $f$ is conjugate to $a$ in $A$.
\end{proof}

\begin{lem}
\begin{enumerate}
\item[(1)] For all $u\in\R$ there exist rough translations through $u$ in $A$.
\item[(2)] If $f_1$ is a rough translation through $u_1$ and $f_2$ is a rough translation through $u_2$, then $f_1f_2^{-1}$ is a rough translation through $u_1-u_2$.
\end{enumerate}
\end{lem}
\begin{proof}
We prove (1), leaving (2) as an exercise. Let $(\xi_n)$ be a $\Z$-sequence
\begin{center}$\dots,-2\frac12,-2,-1,0,1,2,2\frac12,3,3\frac13,3\frac23,4,4\frac14,\dots$\end{center}
i.e., $0,\pm m\pm \frac rm$ for $0\leq r\leq m-1$ and $m=1,2,\dots$.

Let $(\eta_n)$ be a monotone increasing sequence of rationals converging to $u$, and similarly $(\eta_n')$ monotone decreasing to $u$. Define
\[ b_n=\begin{cases}\xi_n+\eta_n&n>0,
\\ \xi_n+\eta_n'&n<0
\\ \xi_n & n=0.\end{cases}\qquad \text{(Need $\eta_1'-\eta_1<1$.)}\]
Then $(b_n)$ is a $\Z$-sequence.  If $f\in A$, $f:\xi_n\to b_n$ (Lemma  \ref{15.1}), then $f$ is a rough translation through $u$.
\end{proof}

\begin{lem}\label{15.5} Suppose that $u>0$ and $f$ is a rough translation through $u$. For any $g\in A$ there exist conjugates $f_1,\dots,f_6$ of $f$ such that $g=f_1^{-1}f_2f_3f_4^{-1}f_5^{-1}f_6$.
\end{lem}
\begin{proof} Quite easy.\end{proof}

As an aside, we have the following.

\begin{cor}\label{15.6} (The heart of J. Truss's argument in the next lecture.) Let $I$ be any set and let $N$ be a normal subgroup of the Cartesian power $A^I$. If $|A^I:N|<2^{\aleph_0}$ then $N=A^I$.
\end{cor}
\begin{proof} For each real number $u$, let $(a_u)\in A^I$ be a sequence each of whose components is a rough translation through $u$. Since $|\R|=2^{\aleph_0}>|A^I:N|$, so there exist $u_1,u_2\in \R$, $u_1>u_2$ and $(a_{u_1})\equiv (a_{u_2})\bmod N$. Then if $(f)=(a_{u_1})(a_{u_2})^{-1}\in N$, then all components of $f$ are rough translations through $u=u_1-u_2$. If $g\in A^I$ then performing the necessary conjugations componentwise, we obtain
\[ (g)=(f_1)^{-1}(f_2)(f_3)(f_4)^{-1}(f_5)^{-1}(f_6).\]
Hence $(g)\in N$. So $N=A^I$.
\end{proof}

\chapter{Truss's Theorem}
\begin{thm}\label{16.1}
If $G\leq A$ and $|A:G|<2^{\aleph_0}$ then there is a finite set $\Delta_0\subseteq \Q$ such that $G=A_{(\Delta_0)}$.  
\end{thm}

\begin{proof}
{\bf Step 1:}  We can assume that $G$ is transitive on $\Q$.

\smallskip

So assume  \ref{16.1} is true for transitive groups. By  \ref{14.5}, any orbit of $G$ is an interval with rational endpoints, and $G$ has only finitely many orbits.  Let $\Delta_0$ be the set of endpoints.  Then $G\subseteq A_{(\Delta_0)}$ and if $\Delta_0=\{\alpha_1,\ldots,\alpha_n\}$ then $G$ is transitive on the intervals $(\alpha_i, \alpha_{i+1})$ (also write $(\alpha_0, \alpha_1)=(-\infty, a_1)$ and $(\alpha_n,\alpha_{n+1})=(\alpha_n, \infty)$). For each of these intervals $(\alpha_i, \alpha_{i+1})$ the group $G^{(\alpha_i, \alpha_{i+1})}$ is transitive and has index $<2^{\aleph_0}$ in $\Aut((\alpha_i,\alpha_{i+1}), \leq)$ so $G^{(\alpha_i, \alpha_{i+1})}=\Aut((\alpha_i, \alpha_{i+1}), \leq)$.  Then $G\cap A(\alpha_i, \alpha_{i+1})\trianglelefteq A(\alpha_i, \alpha_{i+1})$ has index $<2^{\aleph_0}$ and so $G\cap A(\alpha_i,\alpha_{i+1})= A(\alpha_i, \alpha_{i+1})$. Thus $A_{(\Delta_0)}=A(\alpha_0, \alpha_{1})\times\cdots\times A(\alpha_n, \alpha_{n+1})$ and $G=A_{(\Delta_0)}$. 

\smallskip 

Now assume that $G$ is transitive on $\Q$.

\smallskip

{\bf Step 2:}  $B\leq G$.

\medskip

For any sequence $\xi_1<\xi_2<\cdots$ let $A_n=A(\xi_n,\xi_{n+1})$ and let 
$P=\prod A_n\leq A$. Since $|P:P\cap G|<2^{\aleph_0}$ and $P\cap G\leq \prod H_n$, where $H_n$ is the projection of $P\cap G$ into $A_n$ it follows that $P\cap G$ projects onto $A_n$ for almost all $n$.  If $G^{(\xi_n, \xi_{n+1})}=A(\xi_n, \xi_{n+1})$ then $G\cap A_n\leq A_n$ and so $A_n\leq G$.  Then there exists $n_0$ such that $A(n,n+1)\leq G$ if $n\geq n_0$ and there exists $n_1$ such that $ A(n+\frac{1}{2}, n+ \frac{3}{2})\leq G$ if $n\geq n_1$.  So if $m=\max(n_0,n_1)$ then 
\[\left\langle\left. A(n,n+1), A\left(n+\tfrac{1}{2},n+\tfrac{3}{2}\right)\;\right|\; n\geq m\right\rangle \leq G.\]
But this group fixes all rational numbers $\leq m$ and is transitive on $(m,+\infty)$.  But if $q\in\Q$, then since $G$ is transitive there exists $g\in G$ such that $mg=q$.  Therefore $G$ contains a subgroup fixing $(-\infty, q)$ pointwise and transitive on $(q,+\infty)$.  Hence $G$ is order 2-transitive (i.e.\ transitive on the set of ordered pairs $(\alpha, \beta)$ with $\alpha<\beta$).

But $A(n_0, n_0+1)\leq G$ and so $A(\alpha, \beta)\leq G$ for any $\alpha, \beta$ with $\alpha<\beta$.  Thus $B=\bigcup  A(\alpha, \beta)\leq G$.

\smallskip 

{\bf Step 3:}  If $\Gamma\subseteq \Q$ and $\Gamma$ is of $\Z$-type then $A_{(\Gamma)}\leq G$.

\medskip

Suppose $\Gamma=\{\xi_n\mid n\in\Z\}$.  Express $\Z$ as $\{n(i,j)\mid i\in\Z, j\in\N\}$, such that $n(i,j)+1<n(i+1,j)$ for all $i,j$.  Let $P_j=\prod A(\xi_{n(i,j)}, \xi_{n(i,j)+1})$. Then $A_{(\Gamma)} =P=\prod P_j=\prod  A(\xi_{n(i,j)}, \xi_{n(i,j)+1})$.  The projection of $P\cap G$ to $P_j$ must be surjective for almost all $j$.  Choose $j$ such that $P\cap G$ projects onto $P_j$. Then $P_j\cap G\trianglelefteq P_j$ and by Corollary \ref{15.6}, $P_j\leq G$.  Choose $h\in A$ such that (by Lemma \ref{15.1}) 
$$h=\begin{cases}\xi_{n(i,j)}\mapsto 4i,\\\xi_{n(i,j+1)}\mapsto 4i+3\end{cases}$$
Then $h^{-1}Gh\geq \prod A(4i, 4i+3)$.  Define  
$$U=\{u\in \R\mid h^{-1}Gh\text{ contains a rough translation through }u\}.$$
Then $U$ is a subgroup of $\R$ and $|\R:U|<2^{\aleph_0}$ so $U$ is dense. Choose $u_0\in U$ so that $1\frac{7}{8}<u_0<2\frac{1}{8}$. Take $g\in h^{-1}Gh$ a rough translation through $u_0$.  Then 
\[g^{-1}A(4i, 4i+3)g= A(4i+u_0+\eta_{4i}, 4i+3+u_0+\eta_{4i+3})
=A(4i+2+\eta'_{4i}, 4i+5+\eta'_{4i+3}),\]
where $|\eta'_{4i}|, |\eta'_{4i+3}|<\frac{1}{4}$ if $|i|>i_0$.  So
$$\prod_{i>|i_0|} A(4i+3, 4i+4)\leq h^{-1}Gh.$$
But 
$$\prod_{i\leq |i_0|} A(4i+3, 4i+4)\leq B\leq h^{-1}Gh.$$
Thus
$$\prod_{i\in\Z} A(4i+3, 4i+4)\leq h^{-1}Gh.$$
Hence $A_{(\Gamma_0)}=h^{-1}Gh$.  where $\Gamma_0$ is the set of integers $4i, 4i+3$.  Therefore 
$$A_{(\Gamma)}=\prod A(\xi_{n(i,j)}, \xi_{n(i,j)+1})\times A(\xi_{n(i,j)+1}, \xi_{n(i+1,j)})\leq G.$$

\smallskip 

{\bf Step 4:}  But $\langle A_{(\Gamma)}\mid \Gamma\text{ is a $\Z$-sequence}\rangle=A$. So $A=G$.

\medskip
\end{proof}

\chapter*{References for Lectures 14, 15, 16}

\begin{itemize}
\item[1.] A.\ M.\ W.\ Glass, \emph{Ordered permutation groups}  (L.M.S\ Lecture Note Series, 55), CUP, Cambridge 1981.
\item[2.] J.\ K.\ Truss, \lq Infinite permutation groups. II. Subgroups of small index\rq, \emph{J.\ Algebra}, 120 (1989), 494--515.
\item[3.]  Peter.\ M.\ Neumann, \lq Subgroups of small index in the order-automorphism group of the rationals\rq, Unpublished typescript (8 pp.), Nov.\ 1986.
\end{itemize}

\subsection*{Notation and terminlogy}

\noindent
$A=\Aut(\Q, \leq)$,

\noindent
$B=\{f\in A\mid \supp(f)\mbox{ is bounded above and below}\}$.

\noindent
If $\alpha<\beta$ then $A(\alpha, \beta)=\{f\in A\mid \supp(f)\subseteq (\alpha, \beta)\}$:  thus $A(\alpha, \beta)\cong A$ and $B=\bigcup_{\alpha<\beta} A(\alpha, \beta)$.

\noindent
If $f\in A$ then $f^*$ will be its natural extension $\R\to \R$.

\noindent
We use $a$ to denote the translation $a: \xi\mapsto \xi+1,\ \xi\in\Q$.  

\noindent
If $u\in \R, f\in A$ and $\xi f-\xi\longrightarrow u$ as $|\xi|\longrightarrow \infty$ then $f$ will be called a \emph{rough translation through} $u$.

\noindent
A $\Z$-\emph{sequence} in $\Q$ is a sequence $(\xi_n)_{n\in\Z}$ such that $\xi_n<\xi_{n+1}$ for all $n$ and $\xi_n\longrightarrow \pm\infty$ as $n\longrightarrow\pm\infty$.  A subset $\Gamma$ of $\Q$ will be said to be of $\Z$-type if it may be enumerated as a $\Z$-sequence; thus $\Gamma$ is of $\Z$-\emph{type} if and only if it is discrete, closed and umbounded above and below as  subset of $\R$.

\newpage
\blankpage

\part{Hilary Term 1989}

\blankpage

\chapter{Imprimitivity}

\begin{defn}
A $G$-invariant equivalence relation is known as a ($G$-) congruence on $\Omega$. 
\end{defn}

Here $G$-invariant means (1)  $\rho\subseteq \Omega\times \Omega$ and $\rho g=\rho$ for all $g\in G$, or 
$$\mbox{(2)}\quad  \omega_1\equiv\omega_2 \mod \rho
\Leftrightarrow \omega_1 g\equiv\omega_2 g \mod \rho.$$

\begin{observation}  [Universal algebra]
If $f:\Omega_1\to \Omega_2$ is a $G$-morphism and $\rho$ is defined by
\[\alpha\equiv \beta \mod\rho \Leftrightarrow \alpha f=\beta f\]
then 
\begin{enumerate}
\item $\rho$ is a $G$-congruence,
\item $f$ factorizes through $\Omega_1/\rho$, that is to say, we have a diagram
\begin{center}
\begin{tikzpicture}[scale=.8]
\node(A) at (0,0)     {$\Omega_1$};

\node(L)  at (4,0)  {$\Omega_2$};

[node distance=1cm,line width=1pt]

\node(R)   at (2.1,-2)      {$\Omega_1/R$};

\draw [->] (A) -- (L); 
\draw [->] (A) -- (R); 
\draw [->] (R) -- (L); 

\node at (2.1,0.3) {$f$};
\node at (0,-1) {nat.~surj.};
\node at (4.2,-1) {injection};
\end{tikzpicture}
\end{center}
and $\Omega_1/\rho\cong \im(f)$.
\end{enumerate}
\end{observation}

Now suppose that $G$ is transitive on $\Omega$ and $\rho$ is a $G$-congruence.
Then $G$ is transitive on $\Omega/\rho$.  

\begin{defn}
A non-empty subset $\Gamma$ of $\Omega$ is called a \emph{block} (or a \emph{block of imprimitivity}) if 
$$\big(\forall g\in G\big) \big(\Gamma g=\Gamma \wedge \Gamma g\cap\Gamma=\emptyset\big).$$
\end{defn}

\begin{observation}
If $\rho$ is a congruence then each $\rho$-class is a block.   Conversely, if $\Gamma$ is a block then there is a unique (given transitivity) congruence $\rho$ of which $\Gamma$ is an equivalence class.
\end{observation}

\begin{proof}
Suppose that $\Gamma$ is a block. Check that $\{\Gamma g\mid g\in G\}$ is a partition of $\Omega$.  The associated equivalence relation is a congruence.
\end{proof}

\begin{defn}
$G$ is said to be \emph{primitive} if $G$ is transitive and the only $G$-congruences are the universal and the trivial relations.  Equivalently, the only blocks are $\Omega$ and singleton sets.  
\end{defn}

Note that a $G$-space $\Omega$ is \emph{congruences simple} if and only if either $|\Omega|=2$ and $G$ acts trivially or $G$ is primitive on $\Omega$.

\begin{thm}
Suppose that $\Omega$ is transitive and $\alpha\in\Omega$.  There is a lattice isomorphism from the lattice ${\cal R}=\{\rho\mid \rho\mbox{ a congruence on }\Omega\}$
to the lattice ${\cal S}=\{H\mid G_\alpha\leq H\leq G\}$.
\end{thm}

\begin{proof}
Given $\rho\in{\cal R}$ define $H(\rho)=\{g\in G\mid \alpha g\equiv \alpha\mod\rho\}$.  Certainly $H(\rho)\leq G$ and $G_\alpha\leq H(\rho)$.  Notice that $\rho(\alpha)$ is the $H(\rho)$-orbit $\alpha H(\rho)$.  Therefore the map ${\cal R}\to {\cal S}$ given by $\rho\mapsto H(\rho)$ is injective and order-preserving.  Now suppose that $H\in {\cal S}$.  let $\Gamma(H)=\alpha H$.  Prove that  $\Gamma(H)$ is a block and thus corresponds to a congruence $\rho(H)$.  The map ${\cal S}\to {\cal R}; H\mapsto \rho(H)$ is inverse to the map the other way.
\end{proof}

\begin{cor}
The transitive $G$-space $\Omega$ is primitive if and only if $G_\alpha$ is a maximal proper subgroup.  A coset space $(G:H)$ is primitive if and only if $H$ is maximal in $G$.
\end{cor}

Now let $\rho$ be an equivalence relation on $\Omega$, such that all $\rho$-classes have the same cardinality (a homogeneous equivalence relation).

\begin{thm}
Let $W=\Aut(\rho)$, let $\Gamma=\rho(\alpha)\ (\alpha\in\Omega)$ and let $\Delta=\Omega/\rho$.  Then

(1)  $W$ is transitive on $\Omega$;

(2)  $W$ has a normal subgroup $U$ isomorphic to $\Sym(\Gamma)^\Delta$;

(3)  $U$ has a complement $V$ isomorphic to $\Sym(\Delta)$;

(4)  elements of $V$ act by conjugation on $U$ in such a way that permutes the factors.
\end{thm}

\begin{observation}
This is a prototype for a wreath product (in this case $\Sym(\Gamma)\Wrr\Sym(\Delta)$).
\end{observation}

\begin{proof} 
(1)  The map $\alpha\to\beta$ extends to a bijection $\rho(\alpha)\to \rho(\beta)$ and this extends to a permutation on all $\Omega$.

(2)  Define
$$U=\{g\in \Sym(\Omega)\mid \Gamma'g=\Gamma'\mbox{ for all }\rho\mbox{-classes }\Gamma'\}=\ker(W^\Delta).$$
Then $U\trianglelefteq G$ and $U=\prod\Sym(\Gamma')\cong \Sym(\Gamma)^\Delta$.

(3)  For each $\delta\in \Delta$ choose a bijection of $\Omega$ with $\Gamma\times\Delta$.  Take 
$$V=\{(1\times g)\in \Sym (\Gamma\times\Delta)\mid g\in\Sym(\Delta)\}
\leq \Aut(\rho)=W.$$
The rest is routine.
\end{proof}

\chapter{Wreath Products}

Let $A$ be an abstract group and $B$ a group acting on a set $\Delta$.

Define $$K=A^\Delta=\{f:\Delta\to A\}=\{(a_\delta)_{\delta\in \Delta}\mid a_\delta\in A\}.$$
Define an action $f\mapsto f^b$ of $B$ on $K$ such that $f^b(\delta)=f(\delta b^{-1})$.

Let $W$ be the \emph{semidirect product} of $K$ by $B$.  That is 
$W=\{fb\mid f\in K, b\in B\}$  $(=K\times B$ as a set).  Multiplication is given by
$$(fb)(gc)=(fg^{b^{-1}})(bc).$$
In other words $(fb)(gc)=f(bgb^{-1})bc$.

If we identify $K$ with $\{fb\mid b=1\}$ and $B$ with $\{fb\mid f=1\}$ then $K\trianglelefteq W$, $B\leq W$, $K\cap B=\{1\}$, $W=KB$.

This group $W$ is known as `the wreath product' of $A$ by $B$.  Notation $A\Wrr B$ or $A\Wrr_\Delta B$.  [Notation like $A\wrr B$ or $A\wrr_\Delta B$, or even $A\wr B$ or $A\wr_\Delta B$, occurs also in the literature.]

If $\Gamma$ is an $A$-space then $W$ has a natural action on $\Gamma\times\Delta$.  Namely $(\gamma, \delta)fb=(\gamma f(\delta), \delta b)$.  

Think of it pictorially:
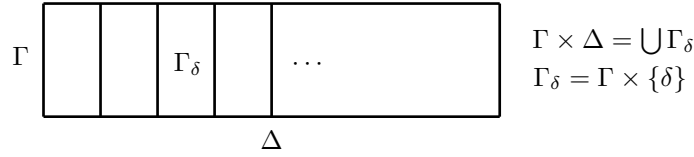
\begin{center}
\begin{tikzpicture}[node distance=5cm,line width=1pt]
\node(A) at (0,0)   {} ;

\node(L)  at (6,0) {} ;

[node distance=1cm,line width=1pt]

\draw [-] (0,0) -- (6,0); 
\draw [-] (0,0) -- (0,-1.5); 
\draw [-] (0,-1.5) -- (6,-1.5); 
\draw [-] (6,-1.5) -- (6,0); 

\draw [-] (0.75,0) -- (0.75,-1.5); 
\draw [-] (1.5,0) -- (1.5,-1.5); 
\draw [-] (2.25,0) -- (2.25,-1.5); 
\draw [-] (3.0,0) -- (3.0,-1.5); 

\node at (3.5, -0.8) {$\cdots$};
\node at (3.0, -1.8) {$\Delta$};
\node at (1.9, -0.8) {$\Gamma_\delta$};
\node at (-0.3, -0.7) {$\Gamma$};
\node at (7.5, -0.5) {$\Gamma\times\Delta=\bigcup \Gamma_\delta$};
\node at (7.45, -1.0) {$\Gamma_\delta=\Gamma\times\{\delta\}$};

\end{tikzpicture}
\end{center}

Need to check that this does satisfy the condition describing an action.

\begin{thm}
(1)  If $A$ is faithful on $\Gamma$ and $B$ is faithful on $\Delta$ then $A\Wrr B$ is faithful on $\Gamma\times \Delta$.

(2)  If $A$ is transitive on $\Gamma$ and $B$ is transitive on $\Delta$ then  $A\Wrr B$ is transitive on $\Gamma\times \Delta$.

(3)  If $\rho$ is defined on $\Gamma\times \Delta$ by 
$$(\gamma_1, \delta_1)\equiv (\gamma_2, \delta_2) \mod \rho\Leftrightarrow \delta_1=\delta_2,$$
then $\rho$ is a $W$-congruence, $W^{(\Gamma\times\Delta)/\rho}=B^\Delta$ and if $\Gamma'$ is a $\rho$-class then $\Gamma'$ is identifiable with $\Gamma$ and $W_{\{\Gamma\}}=A$.  
\end{thm}

\begin{thm}
In this form the wreath product operation is associative in the sense that if $A$ acts on $\Gamma$, and $B$ acts on $\Delta$ and $C$ on $\Phi$ then 
$(A\Wrr_\Delta B)\Wrr_\Phi C$ is equivalent as a permutation group to $A\Wrr_\Delta (B\Wrr_\Phi C)$.
\end{thm}

\begin{proof}
Exercise.
\end{proof}

Suppose now that $G$ acts transitively on $\Omega$ and $\rho$ is a $G$-congruence on $\Omega$.  Set $\Gamma=\rho(\alpha)$ and $\Delta=\Omega/\rho$.

Define $A=G_{\{\Gamma\}}^\Gamma$ and $B=G^\Delta$.

\begin{thm}\label{18.3}
There exists $\phi:\Omega\to \Gamma\times \Delta$ and $\psi:G\to A\Wrr_\Delta B$ such that 

(1)  $\phi$ is bijective and $\psi$ is injective;

(2)  $(\omega^g)\phi=(\omega g)^{g\psi}$.
\end{thm}

\begin{proof}
Choose a transversal $T=\{t_\delta\mid \delta\in \Delta\}$ for $G_{\{\Gamma\}}$ in $G$ chosen so that if $\delta_0\in \Gamma\in\Omega/\rho$ then $\delta_0t_\delta=\delta$.  Define $\phi:\Omega\to \Gamma\times \Delta; \omega\mapsto (\omega t_\delta^{-1}, \delta)$.  Define $\psi$ so that $g\mapsto fb$ where $f(\delta)=t_\delta g t_{\delta g}^{-1}\in A$ and $b=g\ha \Delta$.  Check that all is well.
\end{proof}

\chapter{Variations on Wreath Products}

\noindent 
{\bf Variation I.} [$\Pi$MN, unpublished]

Ingredients:
\begin{itemize}
\item $A$ is a group
\item $\Gamma$ is an $A$-space
\item $B$ is a (top) group
\item $\Delta,\Phi$ are $B$-spaces
\item $\pi:\Delta\to\Phi$ is a $B$-morphism (need not be surjective, though may be if we like)
\item $W=A\Wrr_\Phi B$, and $K=A^\Phi$, $B$ acts on $\Phi$ so $W=K\ltimes B$
\item $W$ acts on $\Gamma\times \Delta$ by $(\gamma,\delta)(fb)=(\gamma f(\delta\pi),\delta b)$.
\end{itemize}

\noindent 
\emph{Notes:}
\begin{enumerate}
\item If $\Phi=\Delta$ and $\pi$ is the identity then we get the ordinary action of $A\Wrr_\Delta B$ on $\Gamma\times\Delta$.
\item If $\Phi$ is a singleton $\{\phi\}$ then $W=A\times B$ acting componentwise on $\Gamma\times \Delta$.
\item Iterations of this construction provide the generalized wreath products over finite posets studied by \hyperref[ref:lec17-20]{Rosemary Bailey et al.\ (1983)}.
\end{enumerate}

\noindent 
{\bf Variation II.} Restricted wreath products.

Replace the Cartesian product by the restricted Cartesian product or an `interdirect product'\footnote{A group sitting between the restricted Cartesian product and the full Cartesian product.}.

Ingredients: $A$, $\Gamma$, $B$, $\Delta$ as before. Let $m$ be an infinite cardinal number.

Let $K=A^{(\Delta,m)}=\{f:\Delta\to A\mid |\supp(f)|<m\}$ where $\supp(f)=\{\delta\mid f(\delta)\neq 1\}$. Again, $B$ acts on $K$ as a group of automorphisms permuting the factors and we write $A\wrr_{(\Delta,m)} B$ to be $K\ltimes B$ (semi-direct product).

\begin{thm}\label{19.1} Let $m$ be an infinite cardinal number. Suppose that $|\Gamma|<m$, and that $B\leq \BS(\Delta,m)$. Then $A\wrr_{(\Delta,m)}B\leq \BS(\Gamma\times\Delta,m)$. Furthermore, if $G$ is a faithful transitive group of $m$-bounded permutations on $\Omega$, and $\rho$ is a proper $G$-congruence, and if $\Gamma=\rho(\alpha)$, $\Delta=\Omega/\rho$, $A=G_{\{\Gamma\}}^\Gamma$, $B=G^\Delta$, then $|\Gamma|<m$, $B\leq \BS(\Delta,m)$, and when we identify $\Omega$ with $\Gamma\times \Delta$ as in Theorem 18.3, we have $G\leq A\wrr_{(\Delta,m)} B$.
\end{thm}

\noindent 
{\bf Variation III.}  \hyperref[ref:lec17-20]{P.~Hall's Wreath Towers (1962)}

Ingredients:
\begin{itemize}
\item $(I,\leq)$ linearly ordered (index) set.
\item $(H_i)_{i\in I}$ a family of groups.
\item $(\Gamma_i)_{i\in I}$ a family of sets, $\Gamma_i$ a faithful $H_i$-space.
\item $\ast=\ast_i\in \Gamma_i$ (a \lq base point\rq).
\item $\Pi=\{(\gamma_i)_{i\in I}\mid \gamma_i=\ast\text{ a.e. }\}=\{(\gamma_i)_{i\in I}\mid \{j\mid \gamma_j\neq \ast\text{ is finite}\}\}$ 
\end{itemize}
$H_j$ acts on $\Pi$ as follows: $(\gamma_i)h=(\gamma_i')$ where
\begin{enumerate}
\item $\gamma_i'=\gamma_i$ for all $i$ if $\exists i_0>j$, $\gamma_{i_0}\neq \ast$,
\item if $\forall i_0>j$ then $\gamma_{i_0}=\ast$, then $\gamma_i'=\begin{cases} \gamma_i,&\text{if $i\neq j$,}\\ \gamma_jh,&\text{if $i=j$.}\end{cases}$
\end{enumerate}
Define $\wrr (H_i)_{i\in I}=\gen{H_j\mid j\in I}\leq \Sym(\Pi)$.

$\Pi$ contains precisely
\[ (\ast,\dots,\ast,\gamma_i,\dots,\gamma_j,\dots,\gamma_k,\ast,\dots).\]
$\Pi=\Pi_{L,j}\times \Pi_{R,j}=\Pi_1\times \Pi_2$.
where $\Pi_{L,j}=(\gamma_i)_{i\leq j}$, $\Pi_{R,j}=(\gamma_i)_{i>j}$. So
\[ \Pi=\bigcup_{\sigma\in \Pi_2} \Pi_1\times \sigma.\]

\begin{example} Let $I=\N$. Then $\wrr (H_i)_{i\in \N}$ is described as follows.
\begin{align*} W_1&=H_1\qquad \text{(acting on $\Gamma_1$)}
\\ W_2&=H_1\wrr H_2\qquad \text{(acting on $\Gamma_1\times \Gamma_2$)}
\\ W_3&=(H_1\wrr H_2)\wrr H_3\qquad \text{(acting on $(\Gamma_1\times \Gamma_2)\times \Gamma_3$)}
\\ &\vdots
\\ W_{n+1}&=(H_1\wrr H_2\wrr\cdots)\wrr H_{n+1}=W_n\wrr H_{n+1}
\\ &\vdots
\end{align*}

There is an embedding $W_n\to W_{n+1}$, $W_n\to K$ (as first factor).

Then $\displaystyle\wrr (H_i)_{i\in \N}=\lim_\rightarrow W_n$ acting on $\lim_\rightarrow (\Gamma_1\times\cdots\times \Gamma_n)$.
\end{example}

\begin{thm} If $I=I_1\dot\cup I_2$ is a (Dedekind) section of $I$, then $\mathrm{wr}(H_i)_{i\in I}=[\mathrm{wr}(H_i)_{i\in I_1}]\wrr [\mathrm{wr}(H_i)_{i\in I_2}]$.
\end{thm}

\chapter{More Variations on the Wreath Product Theme}

\noindent 
{\bf Variation IV.} Generalized wreath product, [cf., \hyperref[ref:lec17-20]{Charles Holland (1969)}, \hyperref[ref:lec17-20]{R. Bailey et al.\ (1983)}].

Ingredients:
\begin{itemize}
\item $(I,\leq)$ a partially ordered set
\item $(H_i,\Gamma_i)_{i\in I}$ permutation groups (transitive on each $\Gamma_i$)
\item $\ast=\ast_i\in\Gamma_i$ base point.
\end{itemize}

If $i\in I$ define $U_i=\{j\in I\mid i<j\}$. Let $\Pi=\{(\gamma_i)_{i\in I}\mid \supp((\gamma_i)_{i\in I})\text{ satisfies MAX (i.e., ACC)}\}$. Here, $\supp (\gamma_i)=\{j\mid \gamma_j\neq \ast\}$.  Set
$\Pi_j=\{(\gamma_i)_{i\in U_j}\mid \supp(\gamma_i)\text{ satisfies MAX}\}$ and
$\pi_j:\Pi\to\Pi_j$ the projection.

Define
\[ W=\{(f_j)_{j\in I}\mid f_j:\Pi_j\to H_j\text{ and $\{j\mid f_j\neq 1\}$ satisfies MAX}\}.\]
Then $W\subseteq\Pi_I (H_j^{\Pi_j})$.

Let $(f_j)_{j\in I}\in W$, let $(\gamma_i)_{i\in I}\in \Pi$ and define
\[ (\gamma_i)(f_j)=(\gamma_i')\in \Pi,\]
where
\[\gamma_i'=\underbrace{\gamma_i\underbrace{f_i\underbrace{\left((\gamma_k)\pi_i\right)}_{\in\Pi_i}}_{\in H_i}}_{\in \Gamma_i}.\]

Assertion: $(f_i)$ acts as a permutation of $\Pi$ and this identifies $W$ with a subgroup of $\Sym(\Pi)$, [cf., \hyperref[ref:lec17-20]{Charles Holland, 1969}]

Idea: Say we have $H_1$ on $\Gamma_1$, $H_2$ on $\Gamma_2$, $H_3$ on $\Gamma_3$.

$H_1\Wrr_{\Gamma_1\times \Gamma_2}(H_2\Wrr_{\Gamma_3} H_3)$. If $\omega$ belongs to this wreath product, then $\omega=f_1(f_2f_3)$, where
\[ f_1:\{\phi\}\to H_3,\quad f_2:\Gamma_3\to H_2,\quad f_1:\Gamma_2\times \Gamma_3\to H_1.\]

\newpage
\emph{Notes:}
\begin{enumerate}
\item In the finite case (\hyperref[ref:lec17-20]{Bailey et al., 1983}) this wreath product is precisely what is obtained by iteration of Variation I.
\item Charles Holland shows that his generalized wreath product is universal for embedding of (faithful) transitive groups with given lattices of congruences. But if we start from $G$ on $\Omega$ with congruences $(\rho_i)_{i\in I}$ satisfying a certain completeness condition, one gets `components' $H_i$ on $\Gamma_i$ of $G$ and embeddings
\begin{align*} \phi&:\Omega\to \Pi,
\\ \psi&:G\to W,
\end{align*}
($W$ on $\Pi$ is the generalized wreath product) but $\phi$ need not be surjective.
\end{enumerate}

\begin{prob} Investigate the structure of $A\Wrr_\Delta B$ as abstract group:
\begin{itemize}
\item When is it directly indecomposable?
\item When is the base group characteristic?
\item What is its automorphism group?
\end{itemize}
\end{prob}

\noindent
\emph{Problem 4$^\#$} Do the same for $\mathrm{wr} (H_i)_{i\in I}$ (Hall's restricted wreath product with fully ordered $I$) and Charles Holland's $\mathrm{Wr} (H_i)_{i\in I}$.

\begin{prob} Let $(\Omega,(\rho_i)_{i\in I})$ be a structure which is a set with a family of equivalence relations. Let
\begin{align*} A&=\Aut(\Omega,(\rho_i)_{i\in I})
\\ &=\{f\in\Sym(\Omega)\mid \rho_if=\rho_i,\;\forall i\}.
\end{align*}
What conditions ensure that $A$ is transitive? What is the structure of $A$?
\end{prob}
Work of \hyperref[ref:lec17-20]{Bailey et al.\ (1983)}, \hyperref[ref:lec17-20]{Behrendt (1986) and (1988)}:

If $|I|=1$ then $A$ is transitive if and only if $\rho$ is `homogeneous' (all classes have same size) and in general
\[ A=\prod_{k\in K} \Sym(\Gamma_k)\Wrr \Sym (A_k).\]
($K$ is the set of cardinalities occurring as sizes of equivalence classes, and $A_k$ is the equivalence classes of size $k$.)

\chapter*{References for Lectures 17--20 on Wreath Products}

\begin{itemize}
\item[1.] \label{ref:lec17-20}R.A.\ Bailey, Cheryl E.\ Praeger, C.A.\ Rowley and T.P.\ Speed, \lq Generalized wreath products of permutation groups\rq, \emph{Proc.\ London Math.\ Soc.\ (3)}, 47 (1983), 69--82.
\item[2.] Gerhard Behrendt, \lq Systems of equivalence relations on countable sets\rq, \emph{Arch.\ Math.}, 47 (1986), 238--242.
\item[3.] Gerhard Behrendt, \lq Equivalence structures and their automorphisms\rq, \emph{European.\ J.\ Combin.}, 9 (1988), 199--206.
\item[4.] P.\ Hall, \lq Wreath powers and characteristically simple groups\rq, \emph{Proc.\ Camb.\ Philos.\ Soc.} 58 (1962), 170--184 = \emph{Collected works of Philip Hall}, Clarendon Press, Oxford, 1988, pp. 609--625.
\item[5.] W.\ Charles Holland, \lq  The characterization of generalized wreath products\rq, \emph{J.\ Algebra}, 13 (1969), 152--172.
\item[6.] Marc Krasner and L\'eo Kaloujine, \lq  Produit complet des groups de permutation et probl\`eme d'extension de groupes, I\rq, \emph{Acta Sci.\ Math.\ (Szeged)}, 13 (1950), 208--230.
\item[7.] Marc Krasner and L\'eo Kaloujine, \lq  Produit complet des groups de permutation et probl\`eme d'extension de groupes, II\rq, \emph{Acta Sci.\ Math.\ (Szeged)}, 14 (1951), 39--66.
\item[8.] Marc Krasner and L\'eo Kaloujine, \lq  Produit complet des groups de permutation et probl\`eme d'extension de groupes, III\rq, \emph{Acta Sci.\ Math.\ (Szeged)}, 14 (1951), 69--82.
\item[9.] Peter M.\ Neumann, \lq On the structure of standard wreath products of groups\rq, \emph{Math.\ Z.}, 84 (1964), 343--373.
\end{itemize}

\chapter{Primitive Permutation Groups}

Recall that $G$ is primitive on $\Omega$ if
\begin{enumerate}
\item $G$ is transitive
\item there are no non-trivial proper $G$-congruences.
\end{enumerate}

(Note that (i) is needed since the example $G=\gen 1$, $\Omega=\{a,b\}$.)

\begin{thm}\label{21.1} The following are equivalent:
\begin{enumerate}
\item $G$ is primitive;
\item $G$ is transitive and there are no non-trivial proper blocks (if $\Gamma$ is a block then for all $g\in G$, $\Gamma=\Gamma g$ or $\Gamma\cap \Gamma g=\emptyset$);
\item $G$ is transitive and a stabilizer $G_\alpha$ is a maximal proper subgroup (that is, $\Omega$ is $G$-isomorphic to a coset space $(G:H)$ where $H$ is maximal).
\end{enumerate}
\end{thm}

\begin{terminology} Let us say $\Gamma\subseteq\Omega$ \emph{separates pairs} if whenever $\alpha,\beta \in \Omega$, $\alpha\neq \beta$, there exists $g\in G$ such that $\Gamma g$ contains one but not both of $\alpha,\beta$.
\end{terminology}

\begin{thm}\label{21.2} Suppose that $G$ is transitive on $\Omega$. Then $G$ is primitive if and only if every non-empty proper subset $\Gamma$ of $\Omega$ separates pairs.
\end{thm}
\begin{proof} Suppose first that $G$ is imprimitive. Let $\Gamma$ be a non-trivial proper block. Then $\Gamma$ does \emph{not} separate pairs since by definition
\[ \forall g\in G\; (\alpha,\beta\in \Gamma g\;\vee\; \alpha,\beta\not\in \Gamma g).\]
Suppose, conversely, that $G$ is primitive, and let $\Gamma$ be a proper, non-empty subset of $\Omega$. Define $\rho_0$ to be the equivalence relation on $\Omega$ whose classes are $\Gamma$ and $\Omega-\Gamma$. Define
\[ \rho=\bigcap_{g\in G} \rho_0^g.\]
Then $\rho$ is a $G$-congruence and proper (i.e., is not the universal relation). Therefore $\rho$ is trivial, i.e., is equality. Now if $\alpha\neq \beta$ then $\alpha\not\equiv \beta\bmod \rho$ and so there exists $g\in G$ such that $\alpha\not\equiv \beta\bmod \rho_0^g$. Thus one of $\alpha,\beta$ lies in $\Gamma g$, the other does not. Hence $\Gamma$ separates pairs.
\end{proof}

\begin{thm}\label{21.3} If $G$ is primitive on $\Omega$ and $N$ is a non-trivial normal subgroup then $N$ is transitive.
\end{thm}
\begin{proof} Given $N\normal G$, define $\rho$ to be equivalence modulo $N$, that is $\alpha\equiv\beta\bmod \rho\Leftrightarrow \exists x\in N\;(\alpha x=\beta)$. Let $g\in G$ and $\alpha\equiv \beta\bmod \rho$. Then $ax=\beta$ ($x\in N$) and
\[ \beta g=\alpha xg=\alpha gg^{-1}xg,\]
so $\beta g=\alpha g\bmod \rho$. So $\rho$ is a congruence. Since $G$ is primitive and $N$ is non-trivial, $\rho$ is universal, hence $N$ is transitive.
\end{proof}

\begin{remark} If $\Omega$ is finite, then this leads to a classification of primitive groups on $\Omega$ by minimal normal subgroups (O'Nan-Scott Theorem).
\end{remark}

\begin{thm}\label{21.4} Suppose that $G$ is faithful and primitive on $\Omega$. Suppose that $G$ has a non-trivial abelian normal subgroup $A$. Then
\begin{enumerate}
\item $A$ is regular on $\Omega$ (i.e., $A$ acts on $\Omega$ like on itself \`a la Cayley);
\item if $H=G_\alpha$ then $G=HA$, $H\cap A=\{1\}$;
\item $A=C_G(A)$ (i.e., $C_H(A)=\{1\}$);
\item if $N$ is any non-trivial normal subgroup then $A\leq N$;
\item $A$ is the additive group of an $F$-vector space for some prime field $F$ and as an $FH$-module (with $H$ acting by conjugation) $A$ is faithful and irreducible;
\item $G\leq \AGL(A)=\{x\mapsto x\alpha+t\mid \alpha\in\Aut(A), t\in A\}\leq \Sym(A)$, the affine general linear group.
\end{enumerate}
\end{thm}
\begin{proof}
\begin{enumerate}
\item $A$ is transitive by  \ref{21.3}. $A_\alpha=H\cap A\normal A$ (because $A$ is abelian) so $A_\alpha$ fixes every point of $\Omega$, hence $A_\alpha=\{1\}$. Thus $A$ acts regularly on $\Omega$.
\item Since $A$ is transitive, $HA=G$. And we just proved that $A_\alpha=H\cap A=\{1\}$.
\item Let $C=C_G(A)$. Then certainly $A\leq C$. As $A$ is transitive, $C=C_\alpha A$ (same argument as before). But $C_\beta=a^{-1}C_\alpha a$ if $a:\alpha\mapsto \beta$ ($a\in A$). Thus $C_\alpha$ fixes all points of $\Omega$, so $C_\alpha=\{1\}$, and therefore $C=A$.
\item Let $A_0\normal G$, $A_0<A$, $A_0\neq \{1\}$. Then $A_0$ is transitive (contains $1$ element from each coset of $\Stab$). So does $A$, hence $A_0=A$. Thus $A$ is a minimal normal subgroup. Let $N\normal G$, $N\neq \{1\}$. If $N\cap A=\{1\}$ then $N\leq C_G(A)=A$ (for if $x\in N$, $y\in A$, then $[x,y]=x^{-1}y^{-1}xy\in N\cap A=\{1\}$) and hence $N=\{1\}$. Thus, if $N\neq\{1\}$ then $N\cap A\neq\{1\}$.  So by minimality of $A$ (as a normal subgroup) $A\leq N$.
\item $A$ is characteristically simple (by minimality of $A$). So if $A$ has torsion then $x^p=1$ for some prime $p$. So $A$ is an elementary $p$-group since the subgroup consisting of $p$-torsion elements is a characteristic subgroup. Otherwise $A$ is torsion-free divisible.
\item Exercise.
\end{enumerate}
\end{proof}

\chapter{Strong Primitivity}

Here $G$ is transitive on $\Omega$,  $g\in G$ and  $\alpha\in \Omega$, .

\begin{definition}[Wielandt] A subset $\Gamma$ of $\Omega$ is called an \emph{$\alpha$-semiblock} if
\begin{enumerate}
\item $\alpha\in\Gamma$
\item $\forall g\in G$, $\alpha g\in \Gamma\Rightarrow \Gamma g\subseteq\Gamma$.
\end{enumerate}
\end{definition}

\begin{comm} Could have defined an $\alpha$-block to be a subset $\Gamma$ such that
\begin{enumerate}
\item $\alpha\in\Gamma$
\item $\forall g\in G$, $\alpha g\in \Gamma\Rightarrow \Gamma g=\Gamma$.
\end{enumerate}
\end{comm}

\noindent
Easy: $\Gamma$ is an $\alpha$-block $\Leftrightarrow$ $\Gamma$ is a block and $\alpha\in\Gamma$.

\begin{ex} $\Omega=\Q$, $G=\Aut(\Q,\leq)$. If $\Gamma=(-\infty,\alpha]$ then $\Gamma$ is an $\alpha$-semiblock.
\end{ex}

Let
\begin{align*} \ms L_1&=\{\Gamma\subseteq\Omega\mid \Gamma\text{ is an $\alpha$-semiblock}\},
\\ \ms L_2&=\{R\subseteq\Omega^2\mid R\text{ is $G$-invariant pre-order}\},
\\ \ms L_3&=\{S\subseteq G\mid S\text{ is a subsemigroup and $G_\alpha\subseteq S$}\}.
\end{align*}
(A pre-order is a reflexive and transitive relation.)

\begin{thm} There are natural isomorphism between $\ms L_1$, $\ms L_2$, $\ms L_3$ (as lattices).
\end{thm}
\begin{proof} Define $\phi_{12}:\ms L_1\to \ms L_2$, $\Gamma\mapsto R_\Gamma$, where $\omega_1 R_\Gamma \omega_2\Leftrightarrow \exists g\in G$ such that $\alpha g=\omega_2\wedge \omega_1\in \Gamma g$. Certainly $R_\Gamma$ is reflexive. Suppose that $\omega_1R_\Gamma\omega_2$ (witnessed by $g$), $\omega_2R_\Gamma\omega_3$ (witnessed by $h$). Then
$$ \alpha g=\omega_2,\quad \omega_1\in \Gamma g,\qquad\mbox{and}\qquad
\alpha h=\omega_3,\quad \omega_2\in \Gamma h.
$$
Then $\alpha gh^{-1}\in \Gamma$ so $\Gamma gh^{-1}\subseteq\Gamma$, so $\Gamma g\subseteq\Gamma h$. Therefore $\omega_1\in \Gamma h$, so $h$ witnesses $\omega_1 R_\Gamma \omega_3$.

Define $\phi_{21}:\ms L_2\to \ms L_1$, $R\mapsto \Gamma_R=\{\omega\in\Omega\mid \omega R\alpha\}$.  Prove that $\Gamma_R$ is an $\alpha$-block, that $\phi_{12}$ and $\phi_{21}$ are inverse maps, that they are complete lattice isomorphisms.

Define $\phi_{13}:\ms L_1\to \ms L_3$, $\Gamma\mapsto S_\Gamma=\{x\in G\mid \alpha x\in \Gamma\}$. Condition (i) gives that $G_\alpha\subseteq S_\Gamma$; if $g\in S_\Gamma$, $h\in S_\Gamma$, then $\alpha g\in \Gamma$, so $\Gamma g\subseteq\Gamma$, whence $\alpha gh\in \Gamma gh\subseteq\Gamma h\subseteq\Gamma$, so $gh\in S_\Gamma$.

Define $\phi_{31}:\ms L_3\to\ms L_1$, $S\mapsto \Gamma_S=\{\alpha x\mid x\in S\}$. Then $\alpha\in\Gamma_S$. Suppose $\alpha g\in \Gamma_S$. There exists $h\in S$ such that $\alpha g=\alpha h$. Then 
\[\Gamma_S g=\{\alpha xg\mid x\in S\}=\{\alpha \underbrace{x}_{\in S}\underbrace{gh^{-1}}_{\in G_\alpha\subseteq S}\underbrace{h}_{\in S}\mid x\in S\}\subseteq\{\alpha y\mid y\in S\}=\Gamma_S.\]
\end{proof}

\begin{definition}[Wielandt] The group $G$ is said to be \emph{strongly primitive} if
\begin{enumerate}
\item $G$ is transitive;
\item there are no non-trivial proper $G$-invariant pre-order relations.
\end{enumerate}
\end{definition}
(Strong primitivity implies primitivity.)

\begin{cor} Suppose that $G$ is transitive. The following are equivalent:
\begin{enumerate}
\item $G$ is strongly primitive;
\item the only $\alpha$-semiblocks are $\{\alpha\}$ and $\Omega$;
\item $G_\alpha$ is a maximal proper subsemigroup of $G$.
\end{enumerate}
\end{cor}

\begin{terminology} Let's say that $\Gamma$ \emph{separates ordered pairs} if, whenever $\omega_1,\omega_2\in\Omega$, $\omega_1\neq \omega_2$, there exists $g\in G$ such that $\omega_1\in\Gamma g$, $\omega_2\not\in\Gamma g$.
\end{terminology}

\begin{thm} Suppose that $G$  is transitive on $\Omega$. Then $G$ is strongly primitive if and only if every non-empty proper subset separates ordered pairs.
\end{thm}
\begin{proof} Suppose $G$ not strongly primitive. A non-trivial proper $\alpha$-semiblock $\Gamma$ does not separate the pair $\alpha, \beta$ if $\beta\in\Gamma$, $\beta\neq\alpha$. Suppose that $G$ is strongly primitive, $\Gamma$ a non-empty proper subset. Suppose $\alpha\in\Gamma$. Define $\Gamma'=\bigcap \{\Gamma g\mid \alpha\in \Gamma g\}$. Then $\Gamma'$ is an $\alpha$-semiblock, so $\Gamma'\in\{\alpha\}$.  Hence $\Gamma$ separates ordered pairs.
\end{proof}

\begin{observation} If $G$ is primitive but not strongly primitive, then there is a $G$-invariant connected partial order relation on $\Omega$.  (A point $\omega_1$ is \lq connected to\rq\ a point $\omega_2$ if there is a sequence $\gamma_0, \ldots, \gamma_n$ such that $\gamma_0=\omega_1$, $\gamma_n=\omega_2$ and $\gamma_i\leq \gamma_{i+1}$ or $\gamma_{i+1}\leq \gamma_{i}$ for all $i=0, \ldots, n-1$.)   
\end{observation}
\begin{proof} There is a proper $G$-invariant binary relation $\leq$ on $\Omega$ which is reflexive and transitive. Define
\[ \omega_1\rho\omega_2\Leftrightarrow \omega_1\leq \omega_2\text{ and }\omega_2\leq \omega_1.\]
Then $\rho$ is a $G$-invariant equivalence relation and is proper. Therefore $\rho$ is `equality'. Hence $\leq$ is antisymmetric, hence a partial order. The relation `connected to' is a non-trivial $G$-congruence, hence is universal, so $(\Omega,\leq)$ is connected.
\end{proof}

Say that $\Gamma$ is \emph{idealistic} if $\Gamma$ separates pairs but does not separate ordered pairs.

\begin{thm} If $\Gamma$ is idealistic and
\[ \omega_1\leq \omega_2\Leftrightarrow \forall g\in G\;(\omega_2 g\in \Gamma\Rightarrow \omega_1g\in\Gamma)\]
then $\leq$ is a $G$-invariant connected partial order on $\Omega$ and $\Gamma$ is an ideal in $(\Omega, \leq)$, i.e., closed downward ($\omega\in\Gamma$ and $\sigma\leq \omega$ implies $\sigma\in\Gamma$).
\end{thm}

\begin{remark} If $G$ is periodic then whenever $\Omega$ is a primitive $G$-space it is strongly primitive.
\end{remark}

\chapter{Suborbits and Orbital Graphs}

\begin{terminology} An orbit of $G_\alpha$ is known as a \emph{suborbit} of $G$ in $\Omega$. (More generally, we'll speak of an orbit of $G_{\alpha_1,\dots,\alpha_d}$ as a \emph{suborbit of depth $d$} (at least if $\alpha_1,\dots,\alpha_d$ are distinct).)

Then $\{\alpha\}$ is referred to as the \emph{trivial} suborbit.
\end{terminology}

\begin{lem} Suppose that $G$ is transitive on $\Omega$. There is a one-one correspondence between the suborbits of $G$ in $\Omega$ and the $G$-orbits in $\Omega^2$. The trivial suborbit corresponds to the diagonal.
\end{lem}
\begin{proof} Given $\Delta\subseteq\Omega^2$, define $\Delta(\alpha)=\{\beta\mid (\alpha,\beta)\in \Delta\}$. The correspondence $\Delta\mapsto \Delta(\alpha)$ from $G$-orbits in $\Omega^2$ to subsets of $\Omega$ is the required one-one correspondence.
\begin{enumerate}
\item $\Delta(\alpha)$ is a $G_\alpha$-orbit ($\neq\emptyset$ because $G$ is transitive).
\item $\Delta\mapsto \Delta(\alpha)$ is one-one.
\item For $\beta G_\alpha$ look at the $G$-orbit of $(\alpha,\beta)$. This proves surjectivity.
\item $\{(\omega,\omega)\mid \omega\in\Omega\}\mapsto \{\alpha\}$.
\end{enumerate}
\end{proof}

\begin{terminology} The $G$-orbits on $\Omega^2$ are known as \emph{orbitals} of $G$. If $\Delta$ is one of these we define the \emph{orbital graph}:
\begin{enumerate}
\item Vertex set is $\Omega$;
\item Edge set is $\Delta$.
\end{enumerate}
This is a directed graph whose automorphism group is transitive (given that $G$ is transitive on $\Omega$) on vertices and on directed edges.
\end{terminology}

If $\Delta$ is an orbital then so also is $\Delta^*=\{(\omega_1,\omega_2)\mid (\omega_2,\omega_1)\in \Delta\}$. This \emph{pairing of orbitals} gives a \emph{pairing of suborbits}.

\begin{lem} Let $\Gamma_1,\Gamma_2$ be suborbits of $G$ in $\Omega$, and let $\gamma_1\in\Gamma_1$, $\gamma_2\in\Gamma_2$. Then $\Gamma_2=\Gamma_1^*$ if and only if there exists $g\in G$ such that $g:\gamma_2\mapsto \alpha\mapsto \gamma_1$.
\end{lem}
\begin{proof} Let $\Delta_1,\Delta_2$ be the orbitals corresponding to $\Gamma_1$ and $\Gamma_2$. Then $(\alpha,\gamma_1)\in \Delta_1$ and $(\alpha,\gamma_2)\in\Delta_2$. So $\Delta_2=\Delta_1^*$ if and only if $(\gamma_2,\alpha)\equiv (\alpha,\gamma_1)\bmod G$.
\end{proof}

\noindent
\emph{Note:} Wielandt (following W.A.~Manning) uses this as the \emph{definition} of pairing: ``$(\dots,\gamma_2,\alpha,\gamma_1,\dots)$ implies $\gamma_2$ is the reflection of $\gamma_1$ by $\alpha$.''

\smallskip

\noindent
\emph{Notes:}\begin{enumerate}
\item The trivial orbit is self-paired.
\item The orbital graph for a self-paired suborbit may be viewed as an \emph{undirected} graph.
\end{enumerate}

\begin{lem} Suppose that $G$ is transitive on $\Omega$. There is a one-one correspondence between suborbits of $G$ and double cosets of $G_\alpha$ in $G$. If $\,\Gamma\leftrightarrow G_\alpha xG_\alpha$ then $\Gamma^*\leftrightarrow G_\alpha x^{-1}G_\alpha$ ($=(G_\alpha x G_\alpha)^{-1}$.)
\end{lem}
\begin{proof} In the isomorphism $\Omega\to (G:G_\alpha)$ (i.e., the right-coset space $\{G_\alpha x\mid x\in G\}$), $G_\alpha$ acts by right multiplication, and so if $\alpha x=\beta$ then the $G_\alpha$-orbit of $\beta$ is the set $\{G_\alpha xh\mid h\in G_\alpha\}$ of right cosets of $G_\alpha$. The union of these is the double coset $(G_\alpha xG_\alpha)$.
\end{proof}

\begin{terminology} The `length' (sizes) of the suborbits are known as the \emph{subdegrees} of $G$.
\end{terminology}

\begin{lem} Suppose that $G$ is transitive on $\Omega$. Let $\Gamma$ be a suborbit, $\beta\in\Gamma$. Then
\[ |\Gamma|=|G_\alpha:G_{\alpha,\beta}|=|G_\alpha:G_\alpha\cap g^{-1}G_\alpha g|,\]
where $\alpha g=\beta$. Thus
\[ |\Gamma^*|=|G_\beta:G_{\alpha,\beta}|=|g^{-1}G_\alpha g:G_\alpha\cap g^{-1}G_\alpha g|.\]
\end{lem}

\begin{lem} If $\Omega$ is finite then $|\Gamma|=|\Gamma^*|$ for any suborbit $\Gamma$.
\end{lem}
\begin{proof} $|\Omega|\cdot|\Gamma|=|\Delta|=|\Delta^*|=|\Omega|\cdot|\Gamma^*|$.
\end{proof}

\begin{example} Let $m_1,m_2$ be non-zero cardinal numbers. If $n\geq \max\{m_1,m_2,\aleph_0\}$ then there is a transitive permutation group $G$ of degree $n$, which has a suborbit $\Gamma$ such that $|\Gamma|=m_1$ and $|\Gamma^*|=m_2$.
\end{example}

An idea of the construction: Choose groups $A_1,A_2$ of orders $m_1$, $m_2$ respectively. Let $A=A_1\times A_2$, let $K=A^\Z$ (direct power, doubly infinite sequences, all but finitely many are 1), and let
\begin{align*} G&=K\cdot \gen b\text{ (the semi-direct product)}
\\ &=A\wrr \gen b.
\end{align*}
Let $H=\{(a_i)_{i\in \Z}\mid (a_i)\in K,a_i\in A_2\text{ if }i<0,\;a_i\in A_1\text{ if } i\geq 0\}$. So $H\leq K\leq G$. Let $\Omega=(G:H)$. Let
\[ X=\{(a_i)_{i\in \Z}\in H\mid a_0=1\}.\]
Then $H=X\times A_1$, $b^{-1}Hb=X\times A_2$. Since $H$ is the stabilizer of a single point, we get
\[ |H:H\cap b^{-1}Hb|=m_1,\quad |b^{-1}Hb:H\cap b^{-1}Hb|=m_2.\]

\chapter{Suborbits in Primitive Groups}

\begin{thm}[{{\hyperref[ref:lec23-28]{D.\ G.\ Higman}}}]\label{24.1} The group $G$ is primitive on $\Omega$ if and only if all the non-trivial orbital graphs are connected (in the undirected sense).
\end{thm}\begin{proof} Suppose first that $\Delta$ is a non-trivial orbital and $(\Omega,\Delta)$ is not connected. The connected components are blocks. So $G$ is imprimitive.

Now suppose that $G$ is imprimitive and let $\Gamma$ be a non-trivial proper block (i.e., $\Gamma g=\Gamma$ or $\Gamma g\cap\Gamma=\emptyset$ for all $g\in G$). Choose $\alpha,\beta\in\Gamma$, $\alpha\neq\beta$. let $\Delta=(\alpha,\beta)G$. Suppose that $\gamma\in\Gamma$ and $(\gamma,\delta)\in\Delta$ or $(\delta,\gamma)\in\Delta$. If $(\gamma,\delta)\in\Delta$ then there exists $g\in G$ such that $(\alpha,\beta)g=(\gamma,\delta)$. Then $\Gamma g\cap\Gamma\neq \emptyset$, so $\Gamma g=\Gamma$, and so $\delta\in\Gamma$. Likewise, if $(\delta,\gamma)\in\Delta$ then $\delta\in\Gamma$. Thus points of $\Gamma$ are not joined to points of $\Omega-\Gamma$ in $(\Omega,\Delta)$, so $(\Omega,\Delta)$ is not connected.
\end{proof}

\begin{observation}\label{24.2}
Let $G$ be a primitive permutation group with a non-trivial suborbit of length $1$. Then $\Omega$ is finite of prime cardinality, and $G$ is the cyclic group acting regularly. 

$C_p$ acting on $C_p$ by multiplication has all suborbits of length $1$. This is primitive because $\{1\}$ is maximal.
\end{observation}
\begin{proof}[Proof of observation] Let $\{\beta\}$ be a non-trivial $G_\alpha$-orbit. Then $G_\alpha\leq G_\beta$. By primitivity, $G_\alpha$ is maximal so $G_\alpha=G_\beta$. Now choose $g$ mapping $\alpha$ to $\beta$. Then $g^{-1}G_ag=G_\beta=G_\alpha$. Thus $g\in N_G(G_\alpha)$, and $G_\alpha<N_G(G_\alpha)$. By maximality, $N_G(G_\alpha)=G$, that is $G_\alpha\normal G$; so $G_\alpha=\{1\}$. But $\{1\}$ is maximal in $G$ if and only if $G$ is cyclic of prime order.
\end{proof}

\noindent\textbf{Exercise.} Give a graph-theoretic proof.

\begin{observation}\label{24.3} Suppose that $G$ is primitive and has some non-trivial finite suborbits, but only finitely many. Then $\Omega$ is finite.
\end{observation}
\begin{proof}[Proof of observation] Let $\Gamma_0=\{\alpha\}$, $\Gamma_1,\dots,\Gamma_s$ be the finite $G_\alpha$-orbits, let $\Delta_i$ be the orbital corresponding to $\Gamma_i$, and let $R=\bigcup \Delta_i\subseteq \Omega\times \Omega$. (We view it as a relation and will show that it is an equivalence relation.)

Certainly $R$ is a $G$-invariant relation on $\Omega$. We prove that it is a congruence. 

Since $\Delta_0\subseteq R$, $R$ is reflexive. 

Also $R$ is transitive: suppose (wlog) $(\alpha,\beta)\in R$, $(\beta,\gamma)\in R$; then there are only finitely many possibilities for $\beta g$ ($g\in G_\alpha$) and given $\alpha,\beta$ there are only finitely many possibilities for $\gamma$.
(Look at paths of length $2$ in the graph.) Hence the $G_\alpha$-orbit containing $\gamma$ is finite; thus $(\alpha,\gamma)\in R$. 

To prove symmetry, let $(\alpha,\beta)\in R$. Then by transitivity $R(\beta)\subseteq R(\alpha)$; but $R(\beta),R(\alpha)$ are finite sets of the same size. So $R(\beta)=R(\alpha)$. Since $\alpha\in R(\alpha)$, also $\alpha\in R(\beta)$, that is $(\beta,\alpha)\in R$. Thus $R$ is an equivalence relation, hence the universal relation, and $\Omega$ is finite.
\end{proof}

\begin{thm} Suppose that $G$ is primitive on $\Omega$ with a non-trivial suborbit $\Gamma$ such that $|\Gamma^*|\leq|\Gamma|=m$. Then $n=|\Omega|\leq \max\{m,\aleph_0\}$. Furthermore, if $m$ is finite then all suborbits are finite.
\end{thm}
\begin{proof} Let $\Delta$ be the orbital corresponding to $\Gamma$, and let $\Delta^{**}=\Delta\cup\Delta^*$. Then $(\Omega,\Delta^{**})$ is an undirected graph, regular (all vertices have same valency $v$). Here $v=m+m^*$ if $\Gamma^*\neq \Gamma$ and $v=m$ if $\Gamma=\Gamma^*$. Define
\[ d(\omega_1,\omega_2)=\min\{\ell\mid \text{ there is a $\Delta^{**}$-path of length $\ell$ between $\omega_1$ and $\omega_2$}\}.\]
By D.\ G.\ Higman's result earlier (24.1) this is a well-defined metric. Define $\Sigma^d(\alpha)=\{\omega\in\Omega\mid d(\alpha,\omega)=d\}$. Since $(\Omega,\Delta^{**})$ is connected (\hyperref[ref:lec23-28]{D.\ G.\ Higman}),
\[ \Omega=\bigcup_{d=0}^\infty \Sigma^d(\alpha).\]
But also $|\Sigma^d(\alpha)|\leq v(v-1)^d$ (proof by induction). Thus if $v$ is infinite then $|\Omega|=v$. If $v$ is finite then $|\Omega|\leq \aleph_0$ and furthermore, since $\Omega$ is a union of finite $G_\alpha$-invariant sets, all suborbits are finite.
\end{proof}

\begin{prob}
For which pairs $m,m^*$ of cardinal numbers $\geq 2$ do there exist primitive groups with paired suborbits $\Gamma,\Gamma^*$ of lengths $m,m^*$?\footnote{David Evans proved that if $m$ is an arbitrary infinite cardinal number then there exist examples with $m^*$ finite, see David M.\ Evans, \lq 
Suborbits in infinite primitive permutation groups\rq,
\emph{Bull.\ London Math.\ Soc.} 33 (2001), 583--590. }
\end{prob}

Partial answers: (1) If $m,m^*\geq \aleph_0$ then such $G$ does exist.  (2)  Yes if $m=\aleph_0$, $m^*=2$ [Hrushovsky]. (3)  If $m,m^*=2$ then $G$ is a finite dihedral group of twice prime order.

\begin{prob} Classify primitive groups with a subdegree equal to $2$.
\end{prob}

\chapter{Transitive Groups with All Suborbits Finite I}

\begin{thm}Suppose that $G$ is primitive and has a pair of non-trivial paired suborbits $\Gamma,\Gamma^*$ (perhaps $\Gamma=\Gamma^*$), both of which are finite. Then $|\Omega|\leq \aleph_0$ and all the suborbits are finite.
\end{thm}

Standard notation: The subdegrees $1=n_0\leq n_1\leq n_2\leq \cdots$ (i.e., lengths of suborbits). Furthermore, if $|\Gamma|=m$, $|\Gamma^*|=m^*$ then $n_k\leq (m+m^*-1)^k$.

\begin{thm} Suppose that $G$ is primitive and that a shortest non-trivial suborbit is self-paired and of length $m$ ($m=n_1$). Then $n_k\leq m(m-1)^k$.
\end{thm}
\begin{prob} How close can one get to this?\footnote{Simon Smith tackled this problem in Simon M.\ Smith,  \lq Subdegree growth rates of infinite primitive permutation groups\rq,
\emph{J.\ Lond.\ Math.\ Soc.\ (2)},  82 (2010), 526--548.}
\end{prob}

Observe that we cannot have $n_k=m(m-1)^k$ for all $k$; for if $G$ were a group with $n_k=m(m-1)^{k-1}$ then the graph $(\Omega,\Gamma)$ where $\Gamma$ is the appropriate orbital would be the \emph{tree} of valency $m$, but the automorphism group of this tree is imprimitive (indeed if $\omega_1\equiv \omega_2\bmod \rho\Leftrightarrow d(\omega_1,\omega_2)\equiv 0\bmod 2$, then this equivalence relation is invariant under the automorphism). Let $G_0$ be the stabilizer of the two blocks; then $G_0$ is primitive on each of the blocks and $G_0$ has subdegrees $1,n_1,n_2,\dots$, where 
\[n_k=[\underbrace{m(m-1)}_{\text{new $m$ from $G_0$}}][(m-1)^2]^{k-1}.\]

\begin{thm} Suppose that $G$ is strongly primitive and has a non-trivial finite suborbit. Then $|\Omega|\leq \aleph_0$, all suborbits are finite and $n_k\leq n_1^k$.
\end{thm}
\begin{proof} Let $\Gamma$ be a non-trivial finite suborbit and let $\Delta$ be the corresponding orbital graph. Let
\[ \Sigma^d=\{\beta\in\Omega\mid\text{there is a directed path from $\alpha$ to $\beta$ of length $d$ and no shorter path exists}\}.\]
Then $\Sigma^d$ is $G_\alpha$-invariant, hence a union of suborbits. Also, $|\Sigma^d|\leq m^d$ where $m=|\Gamma|$. Define
\begin{align*} \Sigma&=\bigcup_d \Sigma^d
\\ &=\{\omega\in\Omega\mid\text{there is a directed path from $\alpha$ to $\omega$}\}.
\end{align*}
Choose $\beta\in\Gamma=\Delta(\alpha)$. By strong primitivity if $\Sigma\neq\Omega$ there exists $g\in G$ such that $\alpha\in\Sigma g$, $\beta\not\in \Sigma g$. But $\Sigma g$ has the property that if $\omega\in\Sigma g$ then $\Delta(\omega)\subseteq \Sigma g$. But this is a contradiction and therefore $\Sigma=\Omega$. (Idea is all suborbits finite and do argument over again.)
\end{proof}

From now on we assume:
\[ (\ast)\begin{cases}(1)\text{ $G$ is transitive};
\\ (2)\text{ that there exists $m\in \N$ such that $|\Gamma|\leq m$ for all suborbits $\Gamma$.}\end{cases}\]

\begin{thm} Suppose that $G$ satisfies $(\ast)$. Then $|\Gamma|=|\Gamma^*|$ for all suborbits $\Gamma$.
\end{thm}
\begin{proof} (Cheryl Praeger, 1988 unpublished.)  Let $\Delta$ be a non-trivial orbital with $|\Delta(\alpha)|\leq d$ and let $d^*=|\Delta^*(\alpha)|$ and wlog $d^*\leq d$. Want to prove $d^*=d$. Let $s\geq 1$ and consider 
\[\Sigma^s(\alpha)=\{\omega\in\Omega\mid\text{there is a directed $\Delta$ path of length $s$ from $\alpha$ to $\omega$}\}.\]
Then $\Sigma^s(\alpha)=\bigcup \Gamma_i$ where $(\Gamma_i)_{i\in I}$ is a family of suborbits. Let $n_i=|\Gamma_i|$. Count sequences $(\omega_0,\omega_1,\dots,\omega_s)$ where $\omega_0=\alpha$, $(\omega_,\omega_{i+1})\in\Delta$ for all $i$. Call such a sequence a \emph{path}.

Let $k_i$ be the number of paths of length $s$ with $\omega_0=\alpha$, $\omega_s$ some fixed point in $\Gamma_i$. Then $d^s=\sum_{i\in I} k_i n_i$. Now do the same for $\Delta^*$. Then the points that can be reached by a $\Delta^*$-path from $\beta$ lie in $\bigcup \Gamma_i^*(\beta)$. So, because the $k_i$ are the same,
\[ (d^*)^s=\sum k_in_i^*.\]
Thus $d^s\leq (\sum k_i)m$, $(d^*)^s\geq \sum k_i$, $d^s\leq m(d^*)^s$. Hence $(d/d^*)^s\leq m$. This is true for all $s$, so $d/d^*\leq 1$, hence $d=d^*$.
\end{proof}

\chapter{Transitive Groups with All Suborbits Finite II: The Bergman--Lenstra Theorem}

\begin{terminology} Let's say that a transitive group $G$ is \emph{almost regular} on $\Omega$ if $G_\alpha$ is finite [Recall: \emph{regular} means transitive and $G_\alpha=\{1\}$]. Let's say that $G$ is \emph{finite-by-almost-regular} if there is a congruence $\rho$ on $\Omega$ such that $|\rho|=|\rho(\alpha)|$ is finite and $G^{\Omega/\rho}$ is almost regular.
\end{terminology}

\begin{thm}[{{\hyperref[ref:lec23-28]{Bergman and Lenstra, 1989}}}]\label{26.1}\footnote{At the time of this lecture $\Pi$MN was unaware that this result had previously been proved by G\"unther Schlichting  in G.~Schlichting, \lq Operationen mit periodischen Stabilisatoren\rq,  \emph{Arch. Math. (Basel)}, 34 (1980),  97--99.   Schlichting's proof uses concepts from the theory of locally compact groups, e.g., the Haar-measure.}  Suppose that $G$ is transitive on $\Omega$. Then $G$ is finite-by-almost-regular if and only if there exists an integer $m$ such that $|\Gamma|\leq m$ for all suborbits $\Gamma$.
\end{thm}
\begin{proof} Suppose that $G$ is finite-by-almost-regular. Let $\rho$ be a congruence with $|\rho|$ finite such that $G^{\Omega/\rho}$ is almost regular. Let $H=\{g\in G\mid \alpha g\equiv \alpha\bmod \rho\}$. Then $H$ is the stabilizer of $\rho(\alpha)$ in $\Omega/\rho$ so $H^{\Omega/\rho}$ is finite of order $m_1$, say. Let $m_2=|\rho|$ and let $m=m_1m_2$. Let $\Gamma$ be an $G_\alpha$-orbit. Then $\rho(\Gamma)=\bigcup_{\omega\in\Gamma}\rho(\omega)$ and $\rho(\Gamma)$ is an $H$-orbit, hence $|\rho(\Gamma)/\rho|\leq m_1$. Therefore $|\Gamma|\leq |\rho(\Gamma)|\leq m_1m_2=m$.

Suppose now that $|\Gamma|\leq m$ for all $G_\alpha$-suborbits $\Gamma$. For each non-empty finite subset $\Phi$ of $\Omega$, define
\[ m(\Phi)=\max\{|\Gamma|\mid \Gamma\text{ a $G_{(\Phi)}$-orbit}\}.\]
Thus $m(\Phi)\leq m$. Now define
\[ m_0=\min\{m(\Phi)\mid \Phi\text{ finite non-empty subset of $\Omega$}\}.\]

\noindent
[Note: $m_0=1$ $\Leftrightarrow$ $G$ is almost regular. For if $m_0=1$ then $\exists\Phi=\{\alpha_0,\dots,\alpha_d\}$ such that $G_{(\Phi)}=\{1\}$. But
\[ G_{(\Phi)}=G_{\alpha_0\alpha_1}\cap G_{\alpha_0\alpha_2}\cap\cdots\cap G_{\alpha_0\alpha_d}.\]
But $|G_{\alpha_0}:G_{\alpha_0\alpha_i}|\leq m$. So $|G_{\alpha_0}:G_{(\Phi)}|\leq m^d$, i.e., $|G_{\alpha_0}|\leq m^d$.]

Back to the proof. Now define
\begin{align*}N=\{x\in G\mid\ &\text{there exists a non-empty finite set $\Phi$ such that $m(\Phi)=m_0$,}\\&\text{$\Gamma x=\Gamma$ for all $G_{(\Phi)}$-orbits $\Gamma$ of length $m_0$}\}.
\end{align*}
We will prove $N\normal G$:

Obviously $1\in N$. Suppose $x,y\in N$. Choose $\Phi,\Psi$ finite sets which witness the fact that $x\in N$, $y\in N$ respectively.

Now let $\Theta=\Phi\cup\Psi$. So $\Theta$ is finite, non-empty and $m(\Theta)\leq m(\Phi)$, $m(\Theta)\leq m(\Psi)$. Since $m_0$ was chosen to be minimal, $m(\Theta)=m_0$. Let $\Gamma$ be a $G_{(\Theta)}$-orbit of length $m_0$. Then $\Gamma\subseteq \Gamma'$ where $\Gamma'$ is a $G_{(\Phi)}$-orbit. But $m_0=|\Gamma|\leq |\Gamma'|\leq m_0$. Thus $\Gamma=\Gamma'$, that is $\Gamma$ is a $G_{(\Phi)}$-orbit. Hence $\Gamma x=\Gamma$.  Likewise $\Gamma y=\Gamma$. Then $\Gamma xy^{-1}=\Gamma$, and so $\Theta$ witnesses the fact that $xy^{-1}\in N$. Finally, if $g\in G$, $x\in N$ and $\Phi$ witnesses $x\in N$ then $\Phi g$ witnesses the fact that $g^{-1}xg\in N$. So $N\normal G$.

Let $\rho$ be equivalence modulo $N$, that is $\omega_1\equiv \omega_2\bmod \rho\Leftrightarrow \exists x\in N$ such that $\omega_1x=\omega_2$. Then $\rho$ is a $G$-congruence on $\Omega$.

Let $H=G_\alpha N$. Then $H$ is the stabilizer of $\rho(\alpha)$ in $G$, and so $H^{\Omega/\rho}$ is the stabilizer of $\rho(\alpha)$ in $\Omega/\rho$. Choose a finite set $\Phi_0$ so that $\alpha\in \Phi_0$ and $m(\Phi_0)=m_0$. Then $G_{(\Phi_0)}\leq G_\alpha\cap N$.

But if $\Phi_0=\{\alpha,\alpha_1,\dots,\alpha_d\}$ then
\[ G_{(\Phi_0)}=G_{\alpha\alpha_1}\cap G_{\alpha\alpha_2}\cap \cdots \cap G_{\alpha\alpha_d},\]
so $|G_\alpha:G_{(\Phi_0)}|\leq m^d$. Therefore $|H:N|=|G_\alpha N:N|=|G_\alpha:G_\alpha\cap N|\leq m^d$. Hence $|H^{\Omega/\rho}|\leq m^d$. Thus $G^{\Omega/\rho}$ is almost regular.

Now let $\Omega_1=\alpha N=\rho(\alpha)$. Suppose by way of contradiction that $|\Omega_1|>m$. Then there exist $x_0,x_1,\dots,x_m\in N$ such that $\alpha x_0,\alpha x_1,\dots,\alpha x_n$ are all different. Let $\Phi_0,\Phi_1,\dots,\Phi_m$ be non-empty finite sets such that $\Phi_i$ witnesses $x_i\in N$. Let
\[ \Phi=\{\alpha x_0,\alpha x_1,\dots,\alpha x_m\}\cup\Phi_0\cup\Phi_1\cup\cdots \cup \Phi_m.\]
Then $\Phi$ is finite and non-empty. Let $\Gamma$ be a $G_{(\Phi)}$-orbit of length $m_0$. Then $\Gamma$ is a $G_{(\Phi_i)}$-orbit for each $i$. So $\Gamma x_i=\Gamma$ for all $i$. Let $\gamma\in \Gamma$. For each $i$ there exists $y_i\in G_{(\Phi)}$ such that $\gamma x_i=\gamma y_i$. Thus $x_0y_0^{-1},x_1y_1^{-1},\dots,x_my_m^{-1}\in G_\gamma$. But $\alpha x_iy_i^{-1}=\alpha x_i$. Thus $|\alpha G_\gamma|\geq m+1$, a contradiction. hence $|\rho|\leq m$.
\end{proof}

\chapter{Transitive Groups with All Suborbits Finite III: Commentary}

\setcounter{note}{0}

\begin{note} For a faithful almost-regular permutation group $G$, define $s(G)=|G_\alpha|$ (if $G$ is not faithful define $s(G)=s(G^\Omega)$).\end{note}

\begin{thm} If $G$ is transitive on $\Omega$ and all subdegrees are at most $m$ then there is a congruence $\rho$ such that
\begin{enumerate}
\item $|\rho|\leq m$;
\item $G^{\Omega/\rho}$ is almost regular with $s(G^{\Omega/\rho})\leq m^{m^m}$.
\end{enumerate}
\end{thm}
The fact that such bounds exist can be proved by a compactness argument (Bergman--Lenstra). Then B--L prove $s(G^{\Omega/\rho})\leq m^{m^{.^{.^{.^m}}}}$ ($m$ terms in the exponent).

By $(\ast)$, we mean the following condition. Follow proof of  \ref{26.1} but choose $\Phi$ to be a non-empty finite subset of $\Omega$ minimal subject to the condition that
\[ m(\Psi)<m(\Phi)\Rightarrow|\Psi|>m(\Phi)+m+1\]
for $\Psi\neq\emptyset$, $\Psi\supseteq \Phi$.

Can prove that if $N=\gen{G_{(\Phi)}}^G$ and $\rho$ a congruence (equivalence) mod $N$, then $|\rho|\leq m$.

Also, supposing $\alpha\in\Phi$ and $N^*=\ker(G^{\Omega/\rho})$,
\[s(G^{\Omega/\rho})=|G_\alpha N^*/N^*|=|G_\alpha/G_\alpha\cap N^*|
\leq |G_\alpha/G_\alpha\cap N|\leq |G_\alpha:G_{(\Phi)}|.\]
 But \[|G_\alpha:G_{(\Phi)}|=|G_\alpha:\bigcap_{\beta\in\Phi-\{\alpha\}} G_{\alpha\beta}|\leq m^{|\Phi|-1}.\]
Let $\Phi_1=\{\alpha\}$. Choose $\Phi_2,\dots,\Phi_r$ so that
\begin{enumerate}
\item $\Phi_1\subseteq \Phi_2\subseteq\cdots\subseteq \Phi_r$,
\item $m(\Phi_{i+1})< m(\Phi_i)$,
\item $|\Phi_{i+1}|\leq m(\Phi_i)+m+1$,
\item $r$ maximal subject to (i), (ii) and (iii).
\end{enumerate}
This is possible because $r\leq m$ by (ii). The fact that $r$ is maximal tells us that $\Phi_r$ satisfies $(\ast)$ for $\Phi$ by (ii) and (iii). So
\begin{align*} \Phi&=\Phi_r
\\ |\Phi_1|&=1
\\ |\Phi_2|&\leq 2m+1
\\ |\Phi_3|&\leq 2m^2+2m+1
\\&\vdots
\\ |\Phi|&\leq 2(m^{r-1}+m^{r-2}+\cdots+m+1)-1=2\frac{m^r-1}{m-1}-1.\end{align*}
So $|\Phi|<m^m$ if $m\geq 3$. If $m=2$, easy.

\begin{note}\end{note}

\begin{observation}\label{27.2}
Suppose that $\Omega$ is an almost regular $G$-space and $\rho$ a finite congruence. Then $\Omega/\rho$ is almost regular and $s(G^{\Omega/\rho})\leq |\rho| s(G^\Omega)$.

If $H$ is the stabilizer of $\rho(\alpha)$, $G_\alpha\leq H$, so $|H:G_\alpha|=|\rho|$. Thus $s(G^{\Omega/\rho})\leq |H|\leq |\rho|\cdot |G_\alpha|$.
\end{observation}

Let's say that $\Omega$ is \emph{irreducible} if  $s(G^{\Omega/\rho})\geq s(G^\Omega)$ for all finite congruences $\rho$.

\begin{observation} \label{27.3} Suppose that $G$ is faithful and almost regular on $\Omega$. Then $\Omega$ is reducible if and only if there is a finite normal subgroup $K$ of $G$ such that $K_\alpha\neq \{1\}$ (just like  \ref{27.2}).
\end{observation}

\begin{thm} Let $\Omega$ be a faithful almost regular $G$-space. If $\Omega$ is irreducible, then $G_\alpha$ has a regular orbit. Consequently $s(G)$ is a maximal subdegree.
\end{thm}
\begin{proof} Suppose that $g\in G_\alpha$ with only finitely many $G$-conjugates. By Dietzmann's Lemma: since $g$ has finite order $K=\gen{g^G}$ is finite. By  \ref{27.3}, since $\Omega$ is irreducible, $K\cap G_\alpha=\{1\}$. So $g=1$. Apply  \ref{3.1} to $G$ acting by conjugation on the conjugacy classes represented by $G_\alpha-\{1\}$: there exists $h\in G$ such that $h^{-1}G_\alpha h\cap G_\alpha=\{1\}$. Then if $\beta=\alpha h$, then $G_{\alpha\beta}=\{1\}$. So the $G_\alpha$-orbit of $\beta$ is a regular orbit.
\end{proof}

\begin{cor} Suppose that $G$ is primitive on $\Omega$ ($|\Omega|\geq \aleph_0$) and that the maximal subdegree $m$ is finite. Then $G$ is almost regular and $s(G)=m$.
\end{cor}
\begin{ex}
Let $G$ be a Tarski monster for the prime $p$: that is, a non-abelian group all of whose proper subgroups have order $1$ or $p$. (Notice that $G$ is necessarily infinite, $2$-generator, simple, of exponent $p$,\dots.) Recently constructed by Olshanski\u i\footnote{A.\ Yu. Olshanski\u i, \lq Groups of bounded period with subgroups of prime order\rq, Algebra and Logic 21 (1982), 369–418.}. Now if $\Omega=(G:H)$ where $H$ is cyclic of order $p$ then $\Omega$ is a primitive $G$-space and $s(G)=m=p$.
\end{ex}

\chapter{Transitive Groups with All Suborbits Finite IV: More Commentary on Bergman--Lenstra}

\begin{note} Bergman and Lenstra first proved:\end{note}

\begin{thm}\label{28.1} Let $G$ be a faithful transitive group on $\Omega$ with all subdegrees $\leq 2$. Then either $|G_\alpha|\leq 2$ or there is a congruence $\rho$ such that $|\rho|\leq 2$ and $G^{\Omega/\rho}$ is regular.
\end{thm}
\begin{thm} Let $G$ be a faithful transitive permutation group on $\Omega$ and let $p$ be a prime number. If all subdegrees are $1$ or $p$, then either $G_\alpha$ is a faithful permutation group of degree $p$ or there is a congruence with $|\rho|=p$ and $G^{\Omega/\rho}$ is regular.
\end{thm}
\begin{proof}[Proof of  \ref{28.1}] 
First note note that it is impossible that the stabilizer of one point $\beta=\alpha g$ is properly contained in the stabilizer of another point $\alpha$, since if $\gamma=\beta g=\alpha g^2$ we 
get $|G_\alpha:G_\alpha\cap G_\gamma|=|G_\alpha: G_\gamma|=4$ contrary to our assumption.   From this it follows that if $|G_\alpha:G_\alpha\cap G_\beta|=2$ then $|G_\beta:G_\alpha\cap G_\beta|=2$.

Fix a point $\alpha$ in $\Omega$ and let $\Omega'$ denote the set of points $\beta$ in $\Omega$ such that $|\alpha G_\beta|=2$.   Note that if $\Omega'=\emptyset$ then the action of $G$ is regular and there is nothing to prove. 
Define $\sim$, $\approx$ on $\Omega'$ as follows:
\begin{align*}\beta\sim\gamma&\Leftrightarrow \alpha G_\beta=\alpha G_\gamma,\\
\beta\approx\gamma&\Leftrightarrow G_{\alpha\beta}=G_{\alpha\gamma}.\qquad\text{(Sort of a dual property.)}
\end{align*}
Suppose $\beta\not\sim\gamma$. Then $\alpha G_\beta=\{\alpha,\alpha'\}$ and $\alpha G_\alpha=\{\alpha, \alpha''\}$. By assumption $\alpha'\neq \alpha''$. Then $G_{\beta\gamma}$ fixes $\alpha$, so $G_{\beta\gamma}=G_{\alpha\beta\gamma}=G_{\alpha\beta}\cap G_{\alpha\gamma}$. But $G_{\beta\gamma}$ has index $2$ in $G_\beta$ and in $G_\gamma$. So $G_{\alpha\beta}=G_{\alpha\gamma}$ and thus $\beta\approx\gamma$.

Thus for all $\beta,\gamma\in\Omega'$, $\beta\sim\gamma$ or $\beta\approx \gamma$. Therefore either $\sim$ is the universal relation on $\Omega'$ or $\approx$ is the universal relation.

Suppose that $\sim$ is universal. Then $\alpha G_\beta=\alpha G_\gamma$ for all $\beta,\gamma\in\Omega'$.   Then there is an $\alpha'\neq \alpha$ such that $\alpha G_\beta=\{\alpha, \alpha'\}$ for all $\beta\in \Omega'$.   Since the set $\Omega'$ is invariant under $G_\alpha$ we see that the set $\{\alpha, \alpha'\}$ is also invariant under $G_\alpha$.  If $\beta\in \Omega-\Omega'$ then $G_\beta=G_\alpha$ and thus the set $\{\alpha, \alpha'\}$ is invariant under $G_\beta$ for all $\beta\in \Omega$.  Hence $\{\alpha, \alpha'\}$ is invariant under the normal subgroup $N$ generated by all the stabilizers of points and $|N:G_\alpha|=2$.   The orbits of $N$ give the required congruence $\rho$.

Suppose that $\approx$ is universal. Then $G_{\alpha\beta}=G_{\alpha\gamma}$ for all $\beta,\gamma\in\Omega'$. Then $G_{\alpha\beta}$ fixes all points, so $G_{\alpha\beta}=\{1\}$. Hence $|G_\alpha|\leq 2$.
\end{proof}

\begin{thm}[Praeger, unpublished]  Suppose that $G$ is transitive, all subdegrees are finite and any two distinct subdegrees are coprime. (Note: may have several say of subdegree $17$). Then there is a bound on the subdegrees.
\end{thm}

\begin{note} Bergman and Lenstra describe their theorems in terms of commensurability: subgroups $Y_1$ and $Y_2$ of a group $X$ are said to be \emph{commensurable} if $[Y_1:Y_1\cap Y_2]<\infty$ and $[Y_2:Y_1\cap Y_2]<\infty$. Their version of  \ref{28.1} is:\end{note}

\medskip
\noindent\textbf{Theorem  \ref{28.1} (Original Form).} \emph{Let $G$ be a group and $H$ a subgroup. Then $|H:H\cap gHg^{-1}|\leq 2$ for all $g\in G$ if and only if $G$ has a normal subgroup $N$ such that either (a) $H\leq N$ and $|N:H|\leq 2$; or (b) $N\leq H$ and $|H:N|\leq 2$}.

\medskip

The original version of Theorem~\ref{26.1} is:

\medskip
\noindent\textbf{Theorem  \ref{26.1} (Original Form).} \emph{Let G be a group and H a subgroup. 
The following are equivalent:
\begin{enumerate}
\item $\{|H:H\cap gHg^{-1}|\mid g\in G\}$ has a finite upper bound.
\item $H$ is commensurable with $N\normal G$.
\end{enumerate}}

\begin{note} One might hope that B--L is true in some form for infinite cardinal numbers. For example,

\begin{quotation}
Suppose $G$ is transitive on $\Omega$ and that $|\Gamma|\leq m$ for all $G_\alpha$-orbits  $\Gamma$. We'd like that there exists a congruence $\rho$ with $|\rho|\leq f_1(m)$ and $s(G^{\Omega/\rho})\leq f_2(m)$.
\end{quotation}
This is false.

\begin{example} Let $\Sigma$ be any uncountable set. Let $G=\FS(\Sigma)$ (group of finitary permutations, i.e., ones with finite support). Let $\Sigma=\bdisun \Sigma_i$ be a partition where $2\leq |\Sigma_i|<\aleph_0$ for all $i$. Let $H=\{f\in\FS(\Sigma)\mid \Sigma_i f=\Sigma_i, \forall i\}$, the weak direct product $\mathrm{Dr}_{i\in I} (\Sym(\Sigma_i))$. (P. Hall's notation.) Take $\Omega=(G:H)$, the coset space. Then $|\Omega|=|\Sigma|=|G|$. Let $g\in G$. There is a finite subset $J\subseteq I$ such that $\Sigma_i\cap \supp_\Sigma(g)\neq \emptyset\Leftrightarrow i\in J$. So $H\cap g^{-1}Hg\geq \mathrm{Dr}_{i\in I-J} \Sym(\Sigma_i)$, and this has finite index in $H$. Therefore all suborbits of $G$ in $\Omega$ have finite size. But if $\rho$ is any proper congruence, then $G^{\Omega/\rho}$ is faithful and so $|G_{\{\rho(\alpha)\}}|=|\Sigma|$.
\end{example}
\begin{observation} Suppose that $m$ is an infinite cardinal number, $G$ is transitive and $|\Gamma|\leq m$ for all $G_\alpha$-orbits $\Gamma$. If $\Delta\subseteq \Omega$ and $|\Delta|\leq m$ then there is a congruence $\rho$ such that $|\rho|\leq m$ and $\omega_1\equiv\omega_2\bmod \rho$ for all $\omega_1,\omega_2\in\Delta$. So if $|\Omega|>m$ then $G$ is totally imprimitive on $\Omega$.
\end{observation}

The idea is from Lecture 23. Consider the union of the orbital graphs for $\Gamma,\Gamma^*$ where $\Gamma$ ranges over
\[ \{\beta G_\alpha\mid \beta\in\Delta\}.\]
Connected components give the equivalence classes of a congruence.
\end{note}

\chapter*{References for Lectures 23--28 on Suborbits}

\begin{itemize}
\item[1.]\label{ref:lec23-28} George M.\ Bergmann and Hendrik W.\ Lenstra, \lq Subgroups close to normal subgroups\rq, \emph{J.\ Algebra}, 127 (1989), 80--97.
\item[2.] D.G.\ Higman, \lq Intersection matrices for finite permutation groups\rq, \emph{J.\ Algebra}, 6 (1967), 22--42.
\item[3.]W.A.\ Manning, \emph{Primitive groups (Part 1)}, Stanford University Publications, Vol 1, 1921. [This appears to be the only \lq Part\rq\ ever published.]
\item[4.] Peter M.\ Neumann, \lq  Finite permutation groups, edge-coloured graphs and matrices\rq, 82--18,  in \emph{Topics in group theory and computation} (Ed.\ M.P.J.\ Curran), Academic Press, 1977.
\item[5.] Charles C.\ Sims, \lq  Graphs and finite permutation groups\rq, \emph{Math.\ Z.}, 95 (1967), 76--86.
\item[6.] Charles C.\ Sims, \lq  Graphs and finite permutation groups, II\rq, \emph{Math.\ Z.}, 103 (1968), 276--281.
\item[7.] Helmut Wielandt, \emph{Finite permutation groups}, Academic Press, New York, 1964.
\end{itemize}

\chapter{Groups Containing Finitary Permutations}
\setcounter{note}{0}
\begin{lem}\label{29.1} Let $G$ be a primitive permutation group on $\Omega$. If $G$ contains a $3$-cycle then $\Alt(\Omega)\leq G$.
\end{lem}
\begin{proof} Let $(\omega_1,\omega_2,\omega_3)\in G$. A subset $\Gamma$ of $\Omega$ will be said to be \emph{good} if $\omega_1,\omega_2,\omega_3\in \Gamma$ and $\Alt(\Gamma)\leq G$. In particular, $\{\omega_1,\omega_2,\omega_3\}$ is good. Suppose that $\Gamma_1,\Gamma_2$ are good. Let $(\alpha,\beta,\gamma)$ be a $3$-cycle in $\Sym(\Gamma_1\cup\Gamma_2)$. If $\alpha,\beta,\gamma\in\Gamma_i$ for some $i=1,2$, then $(\alpha,\beta,\gamma)\in G$. So suppose that $\alpha,\beta\in\Gamma_1$, $\gamma\in\Gamma_2-\Gamma_1$. 
Then some conjugate of $(\alpha,\beta,\gamma)$ by an element of $\Alt(\Gamma_1)$ is $(\omega_1,\omega_2,\gamma)$, which lies in $\Alt(\Gamma_2)$. Thus $g^{-1}(\alpha,\beta,\gamma)g\in\Alt(\Gamma_2)$ for some $g\in\Alt(\Gamma_1)$. Hence $(\alpha,\beta,\gamma)\in G$. But $\Alt(\Gamma_1\cup\Gamma_2)$ is generated by $3$-cycles. Hence $\Gamma_1\cup\Gamma_2$ is good. Consequently, if
\[ \Gamma^*=\bigcup_{\Gamma\text{ good}}\Gamma,\]
then $\Gamma^*$ is good. Now suppose that $\Gamma^*\neq\Omega$. Then there exists $g\in G$ such that $\Gamma^* g\neq \Gamma^*$, but $\Gamma^* g\cap\Gamma^*\neq \emptyset$. If $\Gamma^* g\subseteq \Gamma^*$ then replace $g$ by $g^{-1}$: thus we may assume that $\Gamma^* g-\Gamma^*\neq \emptyset$. But $\Gamma^*\cup\Gamma^* g$ is good: because $\Alt(\Gamma^* g)=g^{-1}\Alt(\Gamma^*)g\leq G$, and since $\Gamma^*g\cap\Gamma^*\neq\emptyset$, $\gen{\Alt(\Gamma^*),\Alt(\Gamma^*g)}=\Alt(\Gamma^*\cup\Gamma^*g)$. This contradicts  the definition of $\Gamma^*$; hence $\Gamma^*=\Omega$.
\end{proof}

\begin{thm}[Wielandt's theorem]\label{29.2} Let $G$ be primitive on $\Omega$, and suppose that $\Omega$ is infinite. If $G$ contains a non-trivial finitary permutation then $\Alt(\Omega)\leq G$.
\end{thm}
\begin{proof} (Mainly John Dixon.) Let $x\in G\cap\FS(\Omega)$, $x\neq 1$, and let $\alpha\in\supp(x)$. Let $\{\alpha\}$, $\Gamma_1,\dots,\Gamma_k$ be the $G_\alpha$-orbits that meet $\supp(x)$. If $\Gamma$ were any other $G_\alpha$-orbit then $\Gamma$ would be fixed setwise by $\gen{G_\alpha,x}$, which is $G$ because $G_\alpha$ is a maximal proper subgroup (since $G$ is primitive) and $x\not\in G_\alpha$. This is impossible, so $\{\alpha\}$, $\Gamma_1,\dots,\Gamma_k$ are the only $G_\alpha$-orbits. Now apply Observation  \ref{24.3}: each of $\Gamma_1,\dots,\Gamma_k$ must be infinite. Therefore by Theorem  \ref{3.1} (separation lemma) there exists $g\in G_\alpha$ such that
\[ [\supp(x)]g\cap\supp(x)=\{\alpha\}.\]
Let $y=g^{-1}xg$ and $z=[x,y]=x^{-1}y^{-1}xy$. Then $z=(\alpha,\alpha y,\alpha x)$ is an element of $G$. Hence by Lemma  \ref{29.1}, $\Alt(\Omega)\leq G$.
\end{proof}
\begin{cor}\label{29.3} If $\Omega$ is infinite and $G$ is a primitive group of finitary permutations then $G=\Alt(\Omega)$ or $G=\FS(\Omega)$.
\end{cor}

\noindent
\textbf{Commentary.}

\begin{note} Proof is unpublished (in a letter from 1983).
\end{note}

\begin{note} The theorem also appears in Wielandt's \emph{Unendliche Permutationsgruppen}.
\end{note}
\begin{note} Jordan proved [1871, 1873, 1875]:
\medskip
\noindent \textbf{Jordan's Minimal Degree Theorem.} Suppose that $G$ is primitive of finite degree $n$ and contains a permutation of degree $m>1$. If $n\geq (m/4+1)(4+m\log(m/2))$ then $A_n\leq G$.
\end{note}

Suppose now that $\Omega$ is infinite, that $G\leq \FS(\Omega)$, and that $G$ is transitive.

\begin{lem}\label{29.4} If $\rho$ is a proper congruence then $|\rho|$ is finite.
\end{lem}
\begin{proof} Choose $\alpha,\beta$ with $\alpha\not\equiv \beta\bmod \rho$. Choose $g$ such that $\alpha g=\beta$. Then $\rho(\alpha)g=\rho(\beta)$; thus $\rho(\alpha)\subseteq \supp(g)$ and $\rho(\alpha)$ is  finite.
\end{proof}

\begin{terminology} If there is a maximal finite congruence then we shall say that $G$ is \emph{finite-by-primitive}. If not then $G$ is said to be \emph{totally imprimitive}.
\end{terminology}

\begin{lem}\label{29.5} If $\Omega$ is uncountable then $G$ is finite-by-primitive.
\end{lem}\begin{proof} Suppose that $G$ is totally imprimitive. Then there is a increasing chain
\[\rho_1<\rho_2<\rho_3<\cdots\]
of finite proper congruences. If $\rho_\omega=\bigcup \rho_i$, then $\rho_\omega$ is a congruence which cannot be proper (since classes are infinite). But $\rho_\omega(\alpha)=\bigcup \rho_i(\alpha)$ so 
\[ |\underbrace{\rho_\omega(\alpha)}_{=\Omega}|=\aleph_0.\]
Hence $|\Omega|=\aleph_0$.
\end{proof}

\chapter{Finitary Permutation Groups: The Classification of Finite-by-Primitive Groups}

$\Omega$ infinite, $G\leq \FS(\Omega)$, transitive. Last time: if $\rho$ proper congruence then $|\rho|$ is finite.

\begin{lem}\label{30.1} Let $\rho$ be a maximal proper congruence. Then $\rho(\alpha)$ is a union of the finite $G_\alpha$-orbits.
\end{lem}
\begin{proof}
Let $\{\Phi_i\}_{i\in I}$ be the family of finite $G_\alpha$-orbits. Let $\Gamma=\bigcup \Phi_i$. Notice that $\alpha\in\Gamma$. Certainly $\rho(\alpha)$ is a union of $G_\alpha$-orbits, hence $\rho(\alpha)\subseteq \Gamma$.

Let $A=G_{\{\rho(\alpha)\}}$. If $\omega\not\equiv \alpha\bmod \rho$ then $\omega A$ contains representatives of every $\rho$-class other than $\rho(\alpha)$ -- because $G^{\Omega/\rho}\geq \Alt(\Omega/\rho)$ which is $2$-transitive. Since $G_\alpha$ has finite index in $A$, all $G_\alpha$-orbits in $\Omega-\{\rho(\alpha)\}$ are infinite. Thus $\Gamma\subseteq \rho(\alpha)$; hence $\rho(\alpha)=\Gamma$.
\end{proof}

\begin{cor}\label{30.2} If $G$ is finite-by-primitive then the maximal proper congruence is unique.
\end{cor}

Let \begin{align*}
\rho_\mathrm{\max}&=\rho=\text{ maximal congruence}
\\ A&=G_{\{\rho(\alpha)\}}=G_{\{\Gamma\}}\qquad \text{($\Gamma=\rho(\alpha)$)}
\\ B&=G_{(\Gamma)}
\\ H&=A/B=G^\Gamma.
\end{align*}
Thus $H$ is a finite transitive permutation group (faithful) on $\Gamma$.

The general theory of imprimitive groups, Theorem  \ref{18.3}, and Theorem  \ref{19.1} give an embedding $G\leq W=H\wrr \FS(\Omega/\rho)$ and we know that $G^{(\Omega/\rho)}$ is $\Alt$ or $\FS$.

Let 
\begin{align*} N&=\{g\in G\mid \omega g\equiv \omega\bmod \rho\text{ for all }\omega\},
\\ K&=\mathop\mathrm{Dr}\limits_{\delta\in\Omega/\rho} H_\delta=\text{base group in $W$}.
\end{align*}So $N=G\cap K$.

\begin{lem}\label{30.3} Let $h\in H$, $\beta\not\equiv\alpha\bmod \rho$. Then there exists $g\in N$ such that $\supp(g)\subseteq \rho(\alpha)\cup\rho(\beta)$ and $g\ha \Gamma=h$.
\end{lem}
\begin{proof} Choose $a\in A$ such that $a^\Gamma=h$. For convenience choose $\gamma$ such that $\rho(\gamma)\subseteq \Fix(a)$. Let $\Gamma_0=\rho(\gamma)$; let $\Gamma_1=\Gamma=\rho(\alpha)$, $\Gamma_2,\dots,\Gamma_r$ be the $\rho$-classes that meet $\supp(a)$. There exists $x\in G$ such that
\[ \Gamma_0x=\Gamma_1,\;\; \Gamma_1x=\Gamma_0,\quad \Gamma_0\cup \Gamma_3\cup\cdots\cup \Gamma_r\subseteq \Fix(x).\]
[For let $x'$ be any member of $G$ mapping $\Gamma_0$ to $\Gamma_1$ and $\Gamma_1$ to $\Gamma_0$. Choose $\rho$ classes $\Gamma_2',\dots,\Gamma_r'$ contained in $\Fix(x')$. Since $G^{\Omega/\rho}\geq \Alt(\Omega/\rho)$ there exists $g'\in G$ mapping $\Gamma_0\mapsto \Gamma_0$, $\Gamma_1\mapsto\Gamma_1$, $\Gamma_i'\mapsto \Gamma_i$ for $2\leq i\leq r$. Then take $x=(g')^{-1}x'g'$.]

Now let $y=x^{-1}a^{-1}xa$. Now $\supp(x^{-1}a^{-1}x)=(\supp(a))x\subseteq (\Gamma_1\cup\cdots \cup \Gamma_r)x=\Gamma_0\cup\Gamma_2\cup\cdots\cup \Gamma_r$. So $\supp(y)\subseteq \Gamma_0\cup\Gamma_1\cup\cdots\cup \Gamma_r$. But if $\omega\in\Gamma_2\cup\Gamma_3\cup\cdots\cup\Gamma_r$, then $\omega y=\omega x^{-1}a^{-1}xa=\omega$. Thus in fact $\supp(y)\subseteq \Gamma_0\cup\Gamma_1$. If $\omega\in \Gamma_1$ then $\omega y=\omega x^{-1}a^{-1}xa=\omega a$. Thus $\Gamma_1 y=\Gamma_1$ and $\Gamma_0 y=\Gamma_0$; thus $y\in N$ and $y^{\Gamma_1}=a^{\Gamma_1}=h$.

Now choose $z\in G$ such that $\rho(\alpha)\subseteq \Fix(z)$ and $\rho(\gamma)z=\rho(\beta)$ [as choice of $x$ above]. Let $g=z^{-1}yz$. Then $g\in N$, $g\ha \Gamma=h$, and $\supp(g)\subseteq \Gamma\cup\rho(\beta)$.
\end{proof}

\begin{thm}\label{30.4} Let $G\leq W$ be the natural embedding. Then $[W:G]\leq 2|H|$.
\end{thm}
\begin{proof}
We first show that $|K:N|\leq |H|$. Fix $\delta_0=\rho(\beta)\in\Omega/\rho$. Let $f\in K$. Then $f$ has the form $f_{\delta_1}f_{\delta_2}\ldots f_{\delta_k}$ where $f_{\delta_i}\in H$. By Lemma  \ref{30.3} there exist $g_1,\dots,g_k\in N$ such that $\supp(g_i)\subseteq \delta_0\cup\delta_i$ and, furthermore, $g_i^{\delta_i}=f_{\delta_i}$ (if $\delta_i=\delta_0$ take $g_i=1$). Then $g_1g_2\ldots g_k=f_0f$ where $f_0\in G_{(\Omega-\delta_0)}$. Thus $f\in f_0^{-1}N$. Hence $K=H_0N$ where $H_0$ acts on $\delta_0$ and fixes all other blocks pointwise. Therefore $|K:N|\leq |H|$. Now $|W:GK|\leq 2$ and $|GK:G|=|K:K\cap G|=|K:N|\leq |H|$. Thus $|W:G|\leq 2|H|$.
\end{proof}

\chapter{A Classification of Finitary Permutation Groups}
\setcounter{note}{0}

Recall  \ref{30.4}: If $\Omega$ is infinite (and $n=|\Omega|$) and $G$ is finite-by-primitive then $G\leq W=H\wrr \FS(n)$ as a subgroup of finite index. Here $H$ is the finite factor of $G$, that is $H=G^\Gamma$ where $\Gamma=\rho(\alpha)$.  Let $\delta_1$ denote $\Gamma$ as an element in $\Delta$.

\begin{thm}\label{31.1} Let $H$ be a finite transitive group on a set $\Gamma$, let $\Delta$ be an infinite set, let $W=H\wrr \FS(\Delta)$ acting on $\Omega'=\Gamma\times\Delta$. Then $W/W'\cong H/H'\times C_2$ and if $G$ is any subgroup of finite index in $W$, then $W'\leq G$. Furthermore $fb\in W'$ if and only if $b\in\Alt(\Delta)$ and $\prod f_\delta\in H'$.
\end{thm}
\begin{proof} We have
\begin{align*}W=\gen{H,\FS(\Delta)\mid\ &g^{-1}hg=h\text{ if $h\in H,g\in \FS(\Delta)$ and $\delta_1g=\delta_1$},\\ & [g^{-1}h_1g,h_2]=1\text{ if $h_1,h_2\in H,g\in\FS(\Delta)$ and $\delta_1g\neq \delta_1$}}.\end{align*}
So 
\[ W/W'=\gen{H/H',C_2\mid [g,h]=1\text{ if $g\in C_2,h\in H/H'$}}=H/H'\times C_2.\]
Suppose that $G\leq W$ and $|W:G|$ finite. Then there exists $N\normal W$, $N\leq G$, $X=W/N$ finite. Let $\phi:W\to X$ be the natural surjective homomorphism. Then
\[ \phi:\Alt(\Delta)\to\{1\}.\]
Therefore $[h_1\phi,h_2\phi]=1$ since there exists $g\in \Alt(\Delta)$ such that $\delta_1g\neq \delta_1$. Thus $H\phi$ is abelian. But also $FS(\Delta)\phi$ centralises $H\phi$. Hence $W\phi$ is abelian. Then $W'\leq N\leq G$.
\end{proof}

\begin{cor}[$\Pi$MN, 1976] \label{31.2} If $G$ is a finite-by-primitive finitary permutation group on $\Omega$ of degree $n$ with finite factor $H$, then $\Omega$ can be identified with $\Omega'$ (as in  \ref{31.1}) in such a way that $W'\leq G\leq W$ where $W=H\wrr \FS(n)$.
\end{cor}
\begin{cor}[$\Pi$MN, 1976] \label{31.3} Up to equivalence, there are only countably many finite primitive finitary permutation groups of any given infinite degree $n$. In particular if $n>\aleph_0$ then up to equivalence there are only countably many transitive finitary groups of degree $n$.
\end{cor}

\begin{thm}\label{31.4} Suppose $\Omega$ is infinite and that $G$ is a transitive but totally imprimitive group of finitary permutations of $\Omega$ (so $\Omega=\aleph_0$). Let $\rho_1<\rho_2<\rho_3<\cdots$ be a properly ascending chain of proper congruences. Let $N_i=\{g\in G\mid \omega g\equiv \omega\bmod \rho_i\;\forall\omega\in\Omega\}$,
$\Gamma_i=\rho_i(\alpha)/\rho_{i-1}$ and $H_i=G^{\Gamma_i}$.
Then:
\begin{enumerate}
\item[(1)] $\bigcup \rho_i$ is the universal relation, so $\bigcup \rho_i(\alpha)=\Omega$;
\item[(2)] If $\rho$ is any proper congruence then $\rho\leq \rho_i$ for some $i$;
\item[(3)] If $\Phi$ is any finite subset of $\Omega$ then $\Phi\equiv \alpha\bmod \rho_i$ for some $i$;
\item[(4)] $G=\bigcup N_i$ where $N_1\leq N_2\leq \cdots$;
\item[(5)] $G$ is embeddable into $\mathrm{wr}^\N(H_i)$.
\end{enumerate}
\end{thm}
\begin{thm}\label{31.5} There are $2^{\aleph_0}$ inequivalent transitive finitary permutation groups of degree $\aleph_0$.
\end{thm}
\begin{proof} Certainly there are $\leq 2^{\aleph_0}$ because $|\FS(\aleph_0)|=\aleph_0$. Let $\ms P$ be a sequence $p_1,p_2,\dots$ of distinct prime numbers. Let $H_i=\Gamma_i=C_{p_i}$, $H_i$ acting regularly on $\Gamma_i$. Let
\[ G_{\ms P}=\mathrm{wr}^\N (H_i).\]
Then $G_{\ms P}$ is a transitive finitary permutation group. There are congruences $\rho_1<\rho_2<\rho_3<\cdots$ with $|\rho_i|=p_1\ldots p_i$.
These are the \emph{only} congruences. Hence $G_{\ms P}$ and $G_{\ms P'}$ are isomorphic (equivalent) if and only if $\ms P=\ms P'$.
\end{proof}

\noindent
\textbf{Commentary:}

\begin{note} If seems just possible that one might `classify' the totally imprimitive groups in terms of P.\ Hall's wreath towers (focus on 5th assertion of  \ref{31.4} showing $G$ is a very large subgroup).
\end{note}

\begin{note} Again if $N\normal G$ then either $N\leq N_i$ for some $i$ or $G'\nleq N$.
\end{note}
\begin{note} Sources are two papers by $\Pi$MN in Archiv der Math, 26, 1975 and 27, 1976\footnote{The two papers referred to are: Peter M. Neumann, \lq The lawlessness of groups of finitary permutations\rq, \emph{Arch.\ Math.\ (Basel)}, 26 (1975), 561--566 and Peter M. Neumann, \lq The structure of finitary permutation groups\rq, \emph{Arch.\ Math.\ (Basel)}, 27 (1976), 3--17.}.
\end{note}

\chapter{Bounded Permutation Groups}

Let $G$ be  a transitive $m^+$-bounded group of permutations of degree $n$ 
on $\Omega$  (so $\deg g\leq m$ for all $g\in G$) and suppose that $\aleph_0\leq n$ and $m^+<n$.

\begin{example}
$(\Omega,\leq)=(\Q\times\omega_1)$ (reverse lexicographically ordered).
\[G=\{g\in \Aut(\Omega,\leq)\mid \supp(g)\text{ bounded above}\}\]
Then $G$ is highly homogeneous and $\deg g\leq \aleph_0$ for all $g\in G-\{1\}$.
\end{example}

\begin{example}
Let $\Omega=\Q\times Y$ as a topological space where $\Q$ is as usual, $Y$ has discrete topology, $|Y|=n$, and $\Omega$ has the product topology.  Let $G=\{g\in \Aut(\Omega)\mid \deg g\leq \aleph_0\}$.  Then $G$ is highly transitive on $\Omega$.  Furthermore, $\deg g=\aleph_0$ for all $g\in G- \{1\}$.
\end{example}  

\begin{lem}   Let $G$ be as specified.  Then there is a unique maximal proper congruence $\rho$ on $\Omega$ and $|\rho|\leq m$.  [So $G$ is degree $m$-by-primitive.]
\end{lem}

\begin{proof}
If $\rho$ is any proper congruence and $\alpha\not\equiv \beta \bmod \rho$ and $g\in G$ with $\alpha g=\beta$ then $\rho(\alpha)\subseteq \supp(g)$.  So $|\rho|\leq m$.  If there is no maximal proper congruence then there is a properly ascending chain $(\rho_i)_{i<m^+}$ and $\bigcup \rho_i$ would be a proper congruence $\rho$.  If $\rho'$ is another congruence then $\sigma =\bigcup (\rho\circ \rho')^n$ is a join of $\rho$ and $\rho'$.  Clearly $\sigma$ can't be universal and $\sigma$ contains $\rho$.   Thus $\sigma=\rho$ and $\rho'\subseteq \rho$.   
\end{proof}

\begin{cor}
Let $X$ be a primitive group of degree $n$ which contains a non-identity permutation of degree $m$ where $m^+<n$.   Let $G:=\{g\in X\mid \deg g\leq m\}$.  Then $G$ is primitive on $\Omega$.
\end{cor}

\begin{proof}
For $G$ is transitive on $\Omega$.  Thus there is a unique maximal proper $G$-congruence $\rho$ on $\Omega$.  But, if $x\in X$ then $\rho x$ is another maximal proper $G$-congruence so $\rho x=\rho$.   Thus $\rho$ is an $X$-congruence  so $\rho$ is trivial.  hence $G$ is primitive.
\end{proof} 

\begin{thm}[SAA+$\Pi$MN\footnote{S,\ A.\ Adeleke and Peter M.\ Neumann, \lq Infinite bounded permutation groups\rq, \emph{J.\ London Math.\ Soc.\ (2)}, 53 (1996), 230--242.}]  If $G$ is as specified and $G$ is primitive then $G$ is highly transitive.
\end{thm}

\medskip

\noindent
\textbf{Connections with Jordan groups}

\medskip

Let $G$ be any transitive group on $\Omega$.  A subset $\Gamma$ of $\Omega$ is said to be a \emph{Jordan set} for $G$ if $|\Gamma|>1$ and $G_{(\Omega- \Gamma)}$ is transitive on $\Gamma$.  We call $\Gamma$ an \emph{improper Jordan set} if $|\Omega-\Gamma|=k<\aleph_0$ and $G$ is $k+1$-transitive.

We say that $\Gamma$ has property $\cal P$ if $\exists H\leq G_{(\Omega- \Gamma)}$ and $H$ acting on $\Gamma$ has property $\cal P$.  Typically we talk of \emph{primitive Jordan sets}, \emph{2-transitive Jordan sets}, ...

\begin{thm}
Let $G$ be as specified.  Let $\Phi\subseteq \Omega$ and $m<|\Phi|<n$.  Then there exists $\Psi\subseteq \Omega$ such that 

(i)  $\Phi\subseteq\Psi$ and $|\Psi|=|\Phi|$.

(ii)  $\Psi$ is a Jordan set for $G$.

(iii)  If $G$ is primitive then $\Psi$ may be chosen to be primitive,  if $G$ is $k$-transitive then $\Psi$ may be chosen to be $k$-transitive, etc.
\end{thm}

For $\omega_1, \omega_2\in\Phi$ choose $g_{\omega_1\omega_2}$ such that $\omega_1 g=\omega_2$.  Set $H_1=\langle g_{\omega_1\omega_2}\mid \omega_1, \omega_2\in\Phi\rangle$, $\Psi_1=\supp(H_1)$ and iterate.
Finally set $\Psi=\bigcup \Psi_i$ and $H=\bigcup H_1$.  

\newpage \blankpage
\part{Trinity Term 1989}

\blankpage

\chapter{Retrospect and Prospect}
\setcounter{prob}{0}
\begin{prob}
\begin{enumerate}
\item Does there exist a transitive group $G$ of degree $n=2^k$ such that there exist $\Gamma,\Delta\subseteq \Omega$ with $|\Gamma|=|\Delta|=k$, $|\Gamma\cap\Delta|=1$ and $\Gamma g\cap\Delta\neq\emptyset$ for all $g\in G$.
\item Can the bound $>2^{2^k}$ in Tomkinson's theorem about cofinitary groups be reduced to $>2^k$?\footnote{In the lecture $\Pi$MN mentioned that these two questions had been answered affirmatively by Zo\'e Chatzidakis and Peter Pappas.  The details can be found in Zo\'e Chatzidakis, Peter Pappas, M.\ J.\ Tomkinson, \lq Separation theorems for infinite permutation groups\rq. \emph{Bull.\ London Math.\ Soc.}, {22} (1990), 344--348.   See also the addendum to Lecture 4 above.} 
\end{enumerate}
\end{prob}

\begin{prob}[Problem 2$^*$.]
What is the least cardinal number $\ell^*(n)$ such that $S=\Sym(\Omega)$ is covered by $\ell^*$ proper subgroups?
\end{prob}

\noindent
\emph{Problem 2$^{**}$.}  What is the least cardinal number $\ell^{**}(n)$ such that $S$ is covered by $\ell^{**}$ proper cosets?

\medskip

Now known $\ell^{**}(n)=n$ and thus $\ell^*(n)\geq n$. Known that $\ell^*(\aleph_0)>\aleph_0$.\footnote{These results can be found in H.\ D.\ Macpherson and Peter M.\ Neumann, \lq Subgroups of infinite symmetric groups\rq,
\emph{J.\ London Math.\ Soc.\ (2) 42},  (1990), 64--84.} 

\begin{prob}[2nd half] Is every proper subgroup of $S$ contained in a maximal subgroup?
\end{prob}

A letter from S. Thomas announces: in ZFC plus something $\exists G<S(\aleph_0)$ such that $\{H\mid G\leq H<S\}$ is a well-ordered chain of order type of the initial ordinal of cardinality $2^{\aleph_0}$.

So the answer is `no'.\footnote{For more details see James E.\ Baumgartner, Saharon Shelah and Simon Thomas, \lq Maximal subgroups of infinite symmetric groups\rq, 
\emph{Notre Dame J.\ Formal Logic} 34 (1993), 1--11, and  Simon Thomas, \lq Aspects of infinite symmetric groups\rq, in \emph{Infinite groups and group rings (Tuscaloosa, AL, 1992)}, 139--145,
Ser. Algebra, 1, World Sci. Publ., River Edge, NJ, 1993.}

\medskip
Notation: $\Omega$ is a set, $n=|\Omega|$, $G\leq \Sym(\Omega)$. Then $\#\Orb(G,\Omega)$ is the number of orbits of $G$ in $\Omega$. Compare $\#\Orb(G,\Omega^{\{m\}})$ (writing $\{\Gamma\mid\Gamma\subseteq\Omega,\;|\Gamma|=m\}$) and $\#\Orb(G,\Omega^{\{k\}})$.

If $m,k$ are finite, $m\leq k$, then $\#\Orb(G,\Omega^{\{m\}})\leq \#\Orb(G,\Omega^{\{k\}})$. 

Let $F$ be a field of characteristic $0$. Then $F^\Omega=\{f :\Omega\to F\}$ is an $F$-vector space and has a $G$-action given by
\[ (f^g)(\omega)=f(\omega g^{-1}).\]
Then $F^\Omega$ is an $FG$-module.

If $f\in \Fix_G(F^\Omega)$ then $f(\omega g^{-1})=f(\omega)$ for all $\omega\in\Omega$ so $f$ is constant over orbits of $\Omega$.

In particular if $\#\Orb(G,\Omega)$ is finite then
\begin{align*} \#\Orb(G,\Omega)&=\dim\Fix_G(F^\Omega)
\\ &=\dim\Hom_{FG}(F,F^\Omega).
\end{align*}

\begin{lem}
There exists an injective morphism of $FG$-modules $F^{\Omega^{\{k-1\}}}\rightarrow F^{\Omega^{\{k\}}}$.  
\end{lem}
\begin{proof}
Write $V_k=F^{\Omega^{\{k\}}}$.

Define $r:V_{k-1}\to V_k$, $f\mapsto \hat f$, where $\hat f(\Gamma)=\sum_{\gamma\in\Gamma} f(\Gamma-\{\gamma\})$.   Clearly $r$ is a morphism of $FG$-modules.

Suppose first that $|\Omega|=n=2k-1$. Then $V_k$, $V_{k-1}$ are identified by $\Gamma\mapsto (\Omega-\Gamma)$. This turns $r$ into $r^*:V_{k-1}\to V_{k-1}$, $r^*:f\mapsto f^*$, where $f^*(\Delta)=\sum_{\Delta'\cap\Delta=\emptyset} f(\Delta')=\hat f(\Omega-\Delta)$ (where $\Delta,\Delta'\in\Omega^{\{k-1\}}$).

Suppose that $\hat f=0$. Then $f^*=0$. Want to get $f=0$. For $0\leq j\leq k-1$, $\Psi\in \Omega^{\{j\}}$, define
\[ F(\Psi)=\sum \{f(\Delta')\mid \Delta'\in\Omega^{\{k-1\}},\;\Delta'\cap\Psi=\emptyset\}.\]
Since $f^*=0$, we have  (assuming $\Delta,\Delta'\in\Omega^{\{k-1\}}$)
$$ 0=\sum_{\substack{\Delta\\\Delta\supseteq \Psi}}\sum_{\substack{\Delta'\\\Delta\cap\Delta'=\emptyset}}f(\Delta')
=\sum_{\substack{\Delta'\\\Delta'\cap \Psi=\emptyset}}\sum_{\substack{\Delta\\\Psi\subseteq\Delta\\\Delta\cap\Delta'=\emptyset}}f(\Delta')
=(k-j)\sum_{\substack{\Delta'\\\Delta'\cap \Psi=\emptyset}}f(\Delta')
=(k-j)F(\Psi).
$$
For $0\leq i\leq k-1$ and $i+j\leq k-1$ and $\Phi\in\Omega^{\{i\}}$ define (assuming $\Phi\cap\Psi=\emptyset$)
\[ F(\Phi,\Psi)=\sum_{\substack{\Phi\subseteq \Delta'\\\Psi\cap\Delta'=\emptyset}} f(\Delta').\]
Then if $\Phi=\emptyset$, then $F(\emptyset,\Psi)=F(\Psi)=0$.

Also, \[F(\Phi,\Psi)=F(\Phi-\{\alpha\},\Psi)-F(\Phi-\{\alpha\},\Psi\cup\{\alpha\}).\]
So induction on $i$ yields $F(\Phi,\Psi)=0$ for all $\Phi,\Psi$. In particular,
\[ 0=F(\Delta,\emptyset)=\sum_{\substack{\Delta\subseteq \Delta'\\\emptyset\cap\Delta'=\emptyset}} f(\Delta')=\sum_{\Delta=\Delta'} f(\Delta')=f(\Delta)\Rightarrow f=0.\]

Thus $r$ is injective in the case $n=2k-1$.  The general case follows by considering subsets $\Omega'\subseteq \Omega$ such that $|\Omega'|=2k-1$ and considering the \lq restriction\rq\ $r': F^{\Omega'^{\{k-1\}}}\rightarrow F^{\Omega'^{\{k\}}}$.    
\end{proof}

\noindent\emph{Problem.} Find a natural specificiation of the map
\[ \theta_{r,s} V_r\to V_s\]
in the case where $\Omega$ is finite and such that $\theta_{rs}\theta_{st}=\lambda \theta_{rt}$.

\medskip

\noindent \textbf{Candidate.} $\theta_{rs}:f\to \hat f$, where
\[ \hat f(\Sigma)=\sum_{\Gamma\in\Omega^{\{r\}}} (-1)^{|\Gamma\cap\Sigma|}f(\Gamma).\]

\begin{thm} \label{33.2} Let $0<k<\aleph_0$. If $n\geq 2k$ then
\[ \#\Orb(G,\Omega^{\{k-1\}})\leq \#\Orb(G,\Omega^{\{k\}}).\]
\end{thm}
\begin{proof} The map $r$ defined above is an injective morphism of $FG$-modules
$F^{\Omega^{\{k-1\}}}\rightarrow F^{\Omega^{\{k\}}}$.  
So \[\dim\Hom(F,F^{\Omega^{\{k-1\}}})\leq \dim\Hom(F, F^{\Omega^{\{k\}}}).\]
\end{proof}

\chapter{Commentary}
\setcounter{note}{0}

Theorem \ref{33.2} implies the following:  

\begin{thm}\label{34.1} Suppose that $|\Omega|\geq 2k$. If $G$ is $k$-homogeneous and $0\leq m\leq k$, then $G$ is $m$-homogeneous.
\end{thm}
(The analogous statement for $k$-transitivity is trivial.)

\begin{note} This theorem is proved by M.\ Brown (1959)\footnote{Morton Brown, \lq Weak $n$-homogeneity implies weak $(n-1)$-homogeneity\rq, \emph{Proc.\ Amer.\ Math.\ Soc.}, 10 (1959), 644--647.}
for topological spaces.
Theorem  \ref{34.1} was discovered by Livingstone and Wagner (1965)\footnote{Donald Livingstone and Ascher Wagner, \lq Transitivity of finite permutation groups on unordered sets\rq, \emph{Math.\ Z.} 90 (1965), 393--403.}  for finite sets. Wielandt (1967)\footnote{Helmut Wielandt, \lq Endliche $k$-homogene Permutationsgruppen\rq, \emph{Math.\ Z.}, 101 (1967), 142.}, Bercov and Hobbey (1970)\footnote{Ronald D.\ Bercov and Charles R.\ Hobby, \lq Permutation groups on unordered sets\rq, \emph{Math.\ Z.}, 115 (1970), 165--168.}, Cameron (1976) produced new proofs to infinite case.
\end{note}

\begin{note}
If $\aleph_0\leq k<n$ and $G$ is $k$-homogeneous, then $G$ is highly transitive.
\begin{proof} Certainly $G$ is transitive. There exists a subgroup $H$ that has an orbit of length $k$: let $(\alpha_i)_{i\in I}$ be $k$ distinct points and choose $h_i$ mapping $\alpha$ to $\alpha_i$. Let $H=\gen{h_i\mid i\in I}$. Then $|H|=k$ so all orbits have length $\leq k$, but $|\alpha H|\geq k$. Since $G$ is $k$-homogeneous $G_{\{\Gamma\}}$ is transitive on $\Gamma$ for all $\Gamma\in\Omega^{\{k\}}$.

Let $\Phi=\{(\Gamma,\gamma)\mid \Gamma\subseteq \Omega^{\{k\}},\gamma\in\Gamma\}$. Then $G$ acts transitively on $\Phi$. Hence $G_\alpha$ is transitive on $k$-sets containing $\alpha$. Thus $G_\alpha$ is $k$-homogeneous on $\Omega-\{\alpha\}$. Then induction gives that $H$ is highly transitive, see \hyperref[ref:lec33-48]{$\Pi$MN (1988)}.
\end{proof}
\end{note}

\begin{note} Intriguing open problem.
\setcounter{prob}{8}
\begin{prob} Is it true that if $\aleph_0\leq m\leq k<n$ and $G$ is $k$-homogeneous, then $G$ must be $m$-homogeneous?\footnote{This question is treated in Saharon Shelah and Simon Thomas, \lq Homogeneity of infinite permutation groups\rq, \emph{Arch.\ Math.\ Logic} 28 (1989), 143--147. Assuming Martin's axiom and that $2^{\aleph_0}>\aleph_2$ they construct a group acting on a set of size $n$ such that $\aleph_1<n<2^{\aleph_0}$ that is $k$-homogeneous for all $\aleph_0< m\leq k<n$ but not $\aleph_0$-homogeneous.}
\end{prob}

There are examples (piecewise automorphism groups), transitive on moieties but not on $m$-sets for $\aleph_0\leq m<n$.
\end{note}

\chapter{Highly Homogeneous Groups I}

\begin{thm}[J.P.J. McDermott]\label{35.1}
Suppose that $G$ is $3$-homogeneous but not doubly transitive on $\Omega$ (and $|\Omega|>3$). Then there is a (dense, open) $G$-invariant linear order relation on $\Omega$.
\end{thm}

\noindent
\emph{Note:} Example $\Aut(\Q,\leq)$.

\begin{proof} Note first that by 33.2 certainly $G$ is $2$-homogeneous. Therefore $G$ has exactly two orbits $\Delta$, $\Delta^*$ in $\Omega^{(2)}$. We show that $\Delta$ is a (strict linear) order.
\begin{enumerate}
\item $(\omega_1,\omega_2)\in\Delta$ implies $(\omega_2,\omega_1)\not\in\Delta$, because $G$ is not $2$-transitive.
\item If $\omega_1\neq \omega_2$ then either $(\omega_1,\omega_2)\in\Delta$ or $(\omega_2,\omega_1)\in\Delta$.

\item We can suppose that $|\Delta(\alpha)|\geq 2$. Choose $\beta,\gamma\in\Delta(\alpha)$. We may suppose $(\beta,\gamma)\in\Delta$. Now suppose $(\omega_1,\omega_2)\in \Delta$, $(\omega_2,\omega_3)\in\Delta$. There exists $g\in G$ such that
\[ \{\omega_1,\omega_2,\omega_3\}g=\{\alpha,\beta,\gamma\}\]
Then $\omega_2g=\beta$ so $\omega_3g=\gamma$ and $\omega_1g=\alpha$. Hence $(\omega_1,\omega_3)\in\Delta$.
So $(\omega_1,\omega_2)\in\Delta$ and $(\omega_2,\omega_3)\in\Delta$ implies $(\omega_1,\omega_3)\in\Delta$. 
\end{enumerate}

Thus $\Delta$ is $G$-invariant and a linear order on $\Omega$. Since $(\Omega,\Delta)$ has a $2$-homogeneous automorphism group, $\Delta$ is dense.
\end{proof}

\begin{cor}\label{35.2} If $G$ is $3$-homogeneous and not $2$-transitive and $\Omega$ is countable then we can identify $\Omega$ with $\Q$ in such a way that $G\leq \Aut(\Q,\leq)$.
\end{cor}
Let $(\Omega,\leq)$ be a linearly ordered set. Define
\[ B_<=\{(\alpha,\beta,\gamma)\in\Omega^3\mid \beta<\alpha<\gamma\text{ or }\gamma<\alpha<\beta\}.\]
$B_<$ is called the \emph{(strict) (linear) betweenness relation} on $\Omega$ derived from $\leq$.

\begin{lem}\label{35.3} If $(\Omega,<)$ is not isomorphic to $(\Omega,>)$ then $\Aut(\Omega,B_<)=\Aut(\Omega,<)$. If $(\Omega,<)$ is isomorphic to $(\Omega,>)$ then $\Aut(\Omega,<)$ is a subgroup of index $2$ in $\Aut(\Omega,B_<)$. There exists $t$ such that $t^2=1$, $t$ reverses $\Omega$. If ($\ast$) is satisfied then $\Aut(\Omega,B_<)$ is highly homogeneous  $2$-transitive but not $3$-transitive. The condition ($\ast$) is
\[ (-\infty,\alpha)\cong (-\infty,\beta),\qquad (\alpha,+\infty)\cong(\beta,+\infty)\qquad\text{for all $\alpha,\beta\in\Omega$}\]
and
\[ (\alpha,\beta)\cong(\gamma,\delta)\qquad\text{if $\alpha<\beta$ and $\gamma<\delta$}.\]
\end{lem}

\noindent
\emph{Note:} ($\ast$) is precisely the conditions that $\Aut(\Omega,<)$ is $2$-homogeneous; hence precisely the condition that $\Aut(\Omega,<)$ is highly homogeneous.

\medskip

Define 
\[C_<=\{(\alpha,\beta,\gamma)\in\Omega^3\mid \alpha<\beta<\gamma\text{ or }\beta<\gamma<\alpha\text{ or }\gamma<\alpha<\beta\},\]
the \emph{cyclic (or circular) order} derived from $<$.

\begin{lem}\label{35.4} Suppose that $\Aut(\Omega,<)$ is highly homogeneous (i.e., ($\ast$) holds). If also
\[ (\ast\ast)\qquad (\beta,+\infty)+(-\infty,\alpha)\cong (\alpha,\beta)\text{ whenever $\alpha<\beta$}\]
(sum as ordered sets). Then $\Aut(\Omega,C_<)$ is highly homogeneous $2$-transitive but not $3$-transitive.
\end{lem}

\noindent
\emph{Warning:}  \hyperref[ref:lec33-48]{Peter J.\ Cameron (1976)} defines $C$ differently, namely as the cyclic order from the linear order on $\Omega\cup\{\infty\}$. Then $\Aut(\Omega,<)$ is the stabilizer of $\infty$ in $\Aut(\Omega\cup\{\infty\},C_<)$.

\medskip
\noindent
\emph{Note:}  For $(\Q,<)$ the two constructions give isomorphic cyclic orders.

Define $D_0(\alpha,\beta;\gamma,\delta)$ to mean $\alpha<\gamma<\beta<\delta$ or $\beta<\gamma<\alpha<\delta$.

$D_1(\alpha,\beta;\gamma,\delta)=D_0(\alpha,\beta;\gamma,\delta)\vee D_0(\alpha,\beta;\delta,\gamma)$.

$D(\alpha,\beta;\gamma,\delta)=D_1(\alpha,\beta;\gamma,\delta)\vee D_1(\gamma,\delta;\alpha,\beta)$.

\[D_<=\{ (\alpha,\beta;\gamma,\delta)\in\Omega^4\mid D(\alpha,\beta;\gamma,\delta)\text{ holds}\}.\]
Then $D_<$ is known as the \emph{circular separation relation} derived from $<$. 

\medskip
\noindent
\emph{Note:}  Again PJC defines this on $\Omega\cup\{\infty\}$.

\begin{lem}\label{35.5} If $\Aut(\Omega,<)$ is highly homogeneous, and ($\ast\ast$)
holds, and $(\Omega,<)\cong(\Omega,>)$ then $\Aut(\Omega,D_<)$ is highly homogeneous and $3$-transitive but not $4$-transitive.
\end{lem}

Let $(\Omega,\leq)$ be linearly ordered, $B_<$ be linear betweenness, $C_<$ be the cyclic order, and $D_<$ be circular separation.

\begin{lem}\label{35.6} Let $R$ be a relation on $\Omega$.
\begin{enumerate}
\item $R$ is a linear order if and only if $R\mid_\Delta$ is a linear order for every $3$-set $\Delta$.
\item $R$ is a linear betweenness relation if and only if $R\mid_\Delta$ is a linear betweenness relation for every $4$-set $\Delta$.
\item $R$ is a cyclic order if and only if $R\mid_\Delta$ is a cyclic order for every $4$-set $\Delta$.
\item $R$ is a circular separation relation if and only if $R\mid_\Delta$ is a circular separation relation for every $5$-set $\Delta$.
\end{enumerate}
\end{lem}

\chapter{Highly Homogeneous Groups II}
\setcounter{note}{0}
Let $(\Omega,\leq)$ be linearly ordered., let $B_<$ be the linear betweenness, $C_<$ be the cyclic order and $D_<$ be the circular separation relation.

\begin{thm}[{{\hyperref[ref:lec33-48]{Cameron, 1976}}}]\label{36.1} Suppose that $G$ is highly homogeneous on $\Omega$ and $\Omega$ is infinite. Then either $G$ is highly transitive or there is a $G$-invariant relation $R$ on $\Omega$ which is
\begin{enumerate}
\item a linear order ($R$ is binary),
\item a linear betweenness ($R$ is ternary)
\item a cyclic order ($R$ is ternary)
\item a circular separation relation ($R$ is quaternary).
\end{enumerate}
In particular if $G$ is not highly transitive then $G$ is at most $3$-transitive.
\end{thm}

\noindent
\emph{Note:} Suppose $G$ is $r$-transitive, not $(r+1)$-transitive. Then there is a non-trivial $G$-invariant $(r+1)$-place relation $R$ which is not universal on $\Omega^{(r+1)}$.

The proof is not given here but see one of
\begin{enumerate}
\item[A.] \hyperref[ref:lec33-48]{Cameron (1976)},
\item[B.] \hyperref[ref:lec33-48]{Higman (1977)} (based on an observation of Hodges on a theorem of \hyperref[ref:lec33-48]{Hodges, Lachlan and Shelah (1977)}).
\end{enumerate}

\begin{cor}\label{36.2} If $G$ is $k$-homogeneous and $G$ is not a group of automorphisms of one of those relational systems, then $G$ is $g(k)$-fold transitive where $g(k)\to\infty$ as $k\to\infty$.
\end{cor}

\begin{note} For finite groups we have the splendid theorem of Livingstone and Wagner (1965)\footnote{Donald Livingstone and Ascher Wagner, \lq
Transitivity of finite permutation groups on unordered sets\rq, 
\emph{Math.\ Z.}, 90 (1965), 393--403.} that a $k$-homogeneous group is $k$-transitive if $k\geq 5$.
\end{note}

\begin{note} There exists $G$ which is $k$-homogeneous, $(k-1)$-transitive, not $k$-transitive.

Sketch: Let $L$ be a first-order language with a $k$-place relation $R$.  Let 
\begin{align*}\mc C=\{M\mid\ &M\text{ is finite $L$-structure satisfying: $R(a_1,\dots,a_k)\Rightarrow a_1,\dots,a_k$ distinct,}\\& \text{for all distinct $a_1,\dots,a_k$ $\exists!\sigma\in S_k$, $R(a_{1\sigma},\dots,a_{k\sigma})$}\}.\end{align*}
$\mc C$ is closed under substructures and has finite amalgamation property. By Fraiss\'e's theorem there exists a (unique) countable homogeneous $L$-structure with $\mc C$ as its class of finite substructures (which Higman called skeleton of system, and the French call the age of the system). If $G=\Aut(\Omega,R)$ then certainly $G$ is $k$-homogeneous, not $k$-transitive;: $G$ is $(k-1)$-transitive.
\end{note}

\begin{note} There exist groups $G$ which are:
\begin{enumerate}
\item[(1)] $k$-homogeneous but not $2$-transitive (linear order);
\item[(2)] $k$-homogeneous, $2$-transitive but not $3$-transitive (linear betweenness, cyclic);
\item[(3)] $k$-homogeneous, $3$-transitive but not $4$-transitive (separation);
\item[(4)] (later) $5$-homogeneous, $3$-transitive but not $4$-transitive [Cameron].
\end{enumerate}
\end{note}

Peter Cameron's conjecture: If $k\geq 6$ and $G$ is $k$-homogeneous, not a group of automorphisms of a linear order, etc., then $G$ is at least $(k-1)$-transitive.

He proved in 1976 that $G$ is $n$-transitive where $2n!\sim k$. Hodges in 1978 showed that $G$ being $k$-homogeneous implies $\sim (k+1)/2$-transitive\footnote{The editors have not been able to find a reference for this result.}. Macpherson (1986)\footnote{H.\ D.\ Macpherson, \lq Homogeneity in infinite permutation groups\rq, 
\emph{Period.\ Math.\ Hungar.} 17 (1986), 211--233.} $G$ is $(k-3)$-transitive.

\chapter{Jordan Groups I}

(After Camille Jordan, 1871)

Situation: $\Omega$ a set, $G\leq \Sym(\Omega)$.

\begin{definition}
A subset $\Gamma$ of $\Omega$ is said to be a \emph{Jordan set} if
\begin{enumerate}
\item[(1)] $|\Gamma|>1$
\item[(2)] there exists a subgroup $H$ of $G$ that is transitive on $\Gamma$ and fixes $\Omega-\Gamma$ pointwise. (Equivalently, $G_{(\Omega-\Gamma)}$ is transitive on $\Gamma$.)
\end{enumerate}
\end{definition}

We won't use the latter formulation because we want to use adjectives appropriate for permutation groups and apply them to Jordan sets -- e.g., primitive (if $H$ is primitive on $\Gamma$), $k$-homogeneous, $k$-transitive, cyclic (if $H$ ic cyclic).

The Jordan set $\Gamma$ is said to be \emph{improper} if $|\Omega-\Gamma|=k<\aleph_0$ and $G$ is $(k+1)$-fold transitive. Otherwise $\Gamma$ is a \emph{proper} Jordan set.

A \emph{Jordan group} is a primitive group $G$ that has a proper Jordan set.

\begin{exs}
(1) $\Sym(\Omega)$, $\FS(\Omega)$, $\Alt(\Omega)$.

Any set of size $\geq 2$ ($\geq 3$ for $\Alt(\Omega)$) is a Jordan set.

\medskip
\emph{Note:} By Wielandt's theorem (Theorem~\ref{29.2}) if $\Omega$ is infinite, $G$ is primitive, and $G$ has a finite Jordan set, then $\Alt(\Omega)\leq G$.

\medskip

\noindent (2) Let $\Omega$ be a path-connected $d$-manifold without boundary with $d\geq 2$, and $G=\Homeo(\Omega)$. Then $\Gamma$ is a Jordan set if and only if $\Gamma\neq\emptyset$, $\Gamma$ is open and connected.

\medskip

\noindent (3) $\Omega=\PG(d,F)$ for $d\geq 2$, the projective geometry with $F$ a skew field. $G=\PGL(d+1,F)$. Think of $\Omega$ as the set of $1$-dimensional subspaces of a $(d+1)$-dimensional $F$-space. Then $G$ is $2$-transitive and $\Gamma$ is a Jordan set if and only if $\Omega-\Gamma$ is a subspace.

\medskip

\noindent (4) $\Omega=\AG(d,F)$ for $d\geq 2$, the affine geometry, with $F$ a skew field. 
\[ G=\AGL(d,F)=\{v\mapsto vA+t\mid A\in\GL(d,F),\; t\in F^{(d)}\}.\]
Then $G$ is $2$-transitive ($3$-transitive if $F=\Z_2$). And $\Gamma$ is a Jordan set if and only if $\Omega-\Gamma$ is a subspace.

$M_{22}$, $M_{23}$, $M_{24}$ have Jordan sets of size $16$. There is a subgroup $C_2.A_7$ in $\AGL(4,2)$. These plus automorphisms of geometric figures related to the above are the only \emph{finite} Jordan groups.

In the infinite case there are many more.

\medskip
\noindent(5) $\Omega=\Q$, $G=\Aut(\Q,\leq)$. Then $G$ is primitive (in fact highly homogeneous), $\Gamma$ is a Jordan set if and only if $\Gamma\neq \emptyset$ and $\Gamma$ is convex and open.
\end{exs}

\begin{lem}\label{37.1}
\begin{enumerate}
\item[(1)] If $\Gamma$ is a Jordan set and $g\in G$ then $\Gamma g$ is a Jordan set.
\item[(2)] If $\Gamma_1,\Gamma_2$ are Jordan sets, $\Gamma_1\cap\Gamma_2\neq\emptyset$, then $\Gamma_1\cup\Gamma_2$ is a Jordan set.
\item[(3)] If $(\Gamma_i)_{i\in I}$ is a chain of Jordan sets then $\bigcup \Gamma_i$ is a Jordan set.
\item[(4)] If $(\Gamma_i)_{i\in I}$ is a connected family of Jordan sets then $\bigcup \Gamma_i$ is a Jordan set.
\end{enumerate}
\end{lem}
\begin{proof}   (1) Follows because $g^{-1}Hg\subseteq G_{(\Omega-\Gamma g)}$.

Connected family: draw a graph on $I$ by specifying an edge from $i$ to $j$ if $\Gamma_i\cap\Gamma_j\neq \emptyset$. Want graph connected.  Clearly (2) and (3) follow from (4).

To see (4) we observe that if $H=\langle H_i\rangle$ then $H$ fixes $\Omega-\bigcup\Gamma_i$.  That $G$ is transitive on $\bigcup\Gamma_i$ is proved using connectedness. 
\end{proof}

\begin{coro}
Suppose that $\Gamma_0,\Delta_0\subseteq\Omega$, $\Gamma_0\cap\Delta_0=\emptyset$, and $\Gamma_0$ is a Jordan set. Then there is a Jordan set $\Gamma$ such that $\Gamma_0\subseteq\Gamma$, $\Delta_0\cap\Gamma=\emptyset$ and $\Gamma$ is maximal subject to this.
\end{coro}
\begin{proof} Let $\Gamma=\bigcup\{\Gamma'\mid \Gamma'\text{ is Jordan set, $\Gamma'\supseteq\Gamma_0,\;\Gamma'\cap\Delta_0=\emptyset$}\}$. Then $\Gamma$ is a Jordan set by (4).
\end{proof}

\begin{lem}\label{37.2}
Let $\Gamma$ be a Jordan set. Let $\Gamma_1\subset \Gamma$, $\Gamma_1$ a Jordan set, maximal proper in $\Gamma$, and let $\Delta_1=\Gamma-\Gamma_1$. Then either $\Gamma_1$ or $\Delta_1$ is a block of imprimitivity (possibly $|\Delta_1|=1$) for $G_{\{\Gamma\}}$ in $\Gamma$.
\end{lem}
\begin{proof} Suppose that neither $\Gamma_1$ nor $\Delta_1$ is a block for $G_{\{\Gamma\}}$. Then there exists $g,h\in G_{\{\Gamma\}}$ such that $\Gamma_1g\neq \Gamma_1$ but $\Gamma_1g\cap\Gamma_1\neq\emptyset$, $\Delta_1h\neq\Delta_1$ but $\Delta_1h\cap\Delta_1\neq\emptyset$.

If $\Gamma_1g\subseteq \Gamma_1$ then $\Gamma_1g^{-1}$ would be a Jordan set such that $\Gamma_1\subseteq \Gamma_1g^{-1}\subseteq \Gamma$.  So $\Gamma_1g\not\subseteq \Gamma_1$.  Likewise $\Gamma_1\not\subseteq \Gamma_1g$. But $\Gamma_1g\cup\Gamma_1$ is a Jordan set and as it properly contains $\Gamma_1$ we have $\Gamma_1g\cup\Gamma_1=\Gamma$.

If $\Gamma_1h\subseteq \Gamma_1$ the same argument applies so $\Gamma_1h\not\subseteq \Gamma_1$. So $\Gamma_1h\cap\Gamma_1=\emptyset$ or $\Gamma_1h\cup\Gamma_1=\Gamma$. But $\Gamma_1h\cup\Gamma_1\neq \Gamma$ because $\Delta_1h\cap\Delta_1\neq\emptyset$. Thus $\Gamma_1h\cap\Gamma_1=\emptyset$. Therefore $\Gamma_1h\subseteq \Delta_1\subseteq\Gamma_1g$. Hence $\Gamma_1\subset \Gamma_1gh^{-1}$ contradicting maximality.
\end{proof}

\chapter{Cofinite Jordan Sets}

Assume that $G$ is primitive on $\Omega$, and that $G$ has cofinite Jordan sets which are proper.

\begin{thm}\label{38.1}
\begin{enumerate}
\item[(1)] $G$ is $2$-transitive.
\item[(2)] Any proper cofinite Jordan set is imprimitive.
\end{enumerate}
\end{thm}
\begin{proof} Certainly there exist maximal Jordan sets $\Gamma\neq\Omega$.  Let $\Delta=\Omega-\Gamma$.  Either $\Gamma$ or $\Delta$ is a block of imprimitivity for $G$. But $|\Gamma|\geq 2$, so $\Gamma$ cannot be a block ($G$ is primitive). Therefore $\Delta$ is, and so $|\Delta|=1$. Then $G$ is $2$-transitive.

For (2) let $\Gamma_0$ be a primitive cofinite Jordan set ($\neq \Omega$). Let $\Gamma$ be a Jordan set, minimally properly containing $\Gamma_0$. Such exists because there exists $g\in G$ such that $\Gamma_0g\neq \Gamma_0$ and $\Gamma_0g\cap\Gamma_0\neq\emptyset$. So $\Gamma_0g\cup\Gamma_0$ is a Jordan set properly containing $\Gamma_0$. Let $\Delta_0=\Gamma-\Gamma_0$. Now either $\Delta_0$ or $\Gamma_0$ is a block for $G_{\{\Gamma\}}$. But there is a Jordan set of size $<2|\Gamma_0|$ containing $\Gamma_0$ (in case $\Gamma_0$ is finite). So $|\Delta_0|<|\Gamma_0|$ and therefore $\Gamma_0$ cannot be a block. So $\Delta_0$ is a block. Let $\rho_0$ be the $G_{\{\Gamma\}}$ congruence on $\Gamma$ of which $\Delta_0$ is a class. Then $\rho_0\mid\Gamma_0$ is a $G_{\{\Gamma_0\}\{\Delta_0\}}$-congruence on $\Gamma_0$, hence a $G_{(\Omega-\Gamma_0)}$-congruence on $\Gamma_0$. Since $\Gamma_0$ is a primitive Jordan set, $\rho_0\mid \Gamma_0$ is the trivial relation. Then $|\Delta_0|=1$. Therefore $\Gamma$ is a $2$-transitive Jordan set. Now use induction: $\Gamma$ is improper so $\Gamma_0$ is improper.
\end{proof}

\begin{fact}[Wielandt]\label{38.2} Suppose that $G$ is transitive on $\Omega$. If $\Gamma_1,\Gamma_2$ are cofinite Jordan sets then there exists $g\in G$ such that $\Gamma_1g\subseteq \Gamma_2$ or $\Gamma_2\subseteq \Gamma_1g$.
\end{fact}
\begin{proof} Choose $g$ to maximize $|\Omega-(\Gamma_1g\cup\Gamma_2)|$. Suppose that $\Gamma_1g\not\subseteq \Gamma_2$ and $\Gamma_2\not\subseteq \Gamma_1g$.

Choose $\gamma_1\in \Gamma_1g-\Gamma_2$, $\gamma_2\in \Gamma_2-\Gamma_1g$. Because $G$ is transitive $\Gamma_1g\cap\Gamma_2\neq\emptyset$, so $\Gamma_1g\cup\Gamma_2$ is a Jordan set. Choose $h\in G_{(\Omega-(\Gamma_1g\cup\Gamma_2))}$ such that $\gamma_2h=\gamma_1$. Then $(\Omega-(\Gamma_1g\cup\Gamma_2))\subseteq (\Omega-(\Gamma_1gh\cup\Gamma_2))$. But also $\gamma_2h=\gamma_1\not\in \Gamma_1gh\cup\Gamma_2$. Thus
\[ |\Omega-(\Gamma_1gh\cup\Gamma_2)|>|\Omega-(\Gamma_1g\cup\Gamma_2)|,\]
contradicting maximality.
\end{proof}

\begin{fact}\label{38.3} If $\Gamma_1,\Gamma_2$ are cofinite Jordan sets which are proper, and suppose $G$ is primitive. Then $\Gamma_1\cap\Gamma_2\neq\emptyset$ unless $\Omega$ is finite and $|\Gamma_1|=|\Gamma_2|=n/2$ (and $G=\AGL(d,2)$).
\end{fact}

\noindent
Jordan's idea in the finite case.

\newpage

\begin{thm}\label{38.4} Let $\Gamma_0$ be a proper cofinite Jordan set, and let $\Delta_0=\Omega-\Gamma_0$. Let $G_0=G^{\Delta_0}$. Then
\begin{enumerate}
\item[(1)] $G_0$ is $2$-transitive;
\item[(2)] If $\Gamma$ is any Jordan set properly containing $\Gamma_0$, then $\Gamma-\Gamma_0$ is a Jordan set for $G_0$.
\end{enumerate}
\end{thm}
\begin{proof} Let $(\alpha_1,\beta_1),(\alpha_2,\beta_2)$ be pairs of distinct points of $\Delta_0$. Since $G$ is $2$-transitive, there exists $g\in G$ mapping $(\alpha_1,\beta_1)\mapsto(\alpha_2,\beta_2)$. Then $\Gamma_0g$ is a Jordan set and it meets $\Gamma_0$; hence $\Gamma_0\cup\Gamma_0g$ is a Jordan set and $\alpha_2,\beta_2\not\in\Gamma_0\cup\Gamma_0g$. Now apply  \ref{38.2} to $G_{(\Omega-(\Gamma_0\cup\Gamma_0g))}$: there exists $h$ in this group mapping $\Gamma_0g$ to $\Gamma_0$. Now $gh\in G_{\{\Gamma_0\}}$, that is $gh\in G_{\{\Delta_0\}}$. And $(\alpha_1,\beta_1)gh=(\alpha_2,\beta_2)$.

For the second part, let $\alpha,\beta\in\Gamma-\Gamma_0$. There exists $g\in G_{(\Omega-\Gamma)}$ such that $\alpha g=\beta$. Then $\Gamma_0\cap\Gamma_0g\neq\emptyset$, so if $\Gamma_1=\Gamma_0\cup\Gamma_0g$, then $\Gamma_1$ is a Jordan set. Applying  \ref{38.2} to $G_{(\Omega-\Gamma_1)}$ there exists $h\in G_{(\Omega-\Gamma_1)}$ such that $\Gamma_0gh=\Gamma_0$. Then $gh\in G_{\{\Gamma_0\}}=G_{\{\Delta_0\}}$. But also
\[\alpha gh=\underbrace{\beta}_{\not\in\Gamma_1} h=\beta.\]
Notice $gh$ fixes all points on $\Omega-\Gamma$ (which is $\Delta_0-(\Gamma-\Gamma_0)$). Thus $G_{\{\Delta_0\}(\Delta_0-(\Gamma-\Gamma_0))}$ is transitive on $\Gamma-\Gamma_0$.
\end{proof}

\chapter{Geometry of Cofinite Jordan Sets}
\setcounter{note}{0}
Situation: $G$ is primitive on $\Omega$. There are (proper) cofinite Jordan sets. Define a geometry (see below) on $\Omega$: If $\Gamma$ is finite,
\[\spa(\Gamma)=\begin{cases}\Omega&\text{all cofinite Jordan sets meet $\Gamma$},
\\ \Omega-\Sigma&\text{where $\Sigma$ is the unique maximal Jordan set avoiding $\Gamma$ otherwise}.\end{cases}\]

\begin{thm}\label{39.1} \begin{enumerate}
\item $(\Omega,\spa)$ is a geometry.
\item $G$ is transitive on ordered sequences of $k$ independent points\footnote{A set $A$ of points is \emph{independent} if $a\nin \spa(A-\{a\})$ for all $a\in A$.}.
\end{enumerate}
\end{thm}

Notion of geometry: this requires $\Gamma\subseteq \spa(\Gamma)$, and $\spa(\spa(\Gamma))=\spa(\Gamma)$. We also want $\spa(\emptyset)=\emptyset$, and $\spa(\Gamma)=\Gamma$ if $|\Gamma|=1$. (This property is implied by $G$ being doubly transitive.) We also need the following.

\medskip

\noindent\textbf{Exchange principle.} If $\gamma\in\spa(\Gamma\cup\{\beta\})$ and $\gamma\not\in \spa(\Gamma)$ then $\beta\in\spa(\Gamma\cup\{\gamma\})$ (so $\spa(\Gamma\cup\{\beta\})=\spa(\Gamma\cup\{\gamma\})$).

\begin{proof}[Proof of exchange principle] $\Gamma$ is a finite set, $\gamma\not\in \spa(\Gamma)$. Now $\gamma\in\spa(\Gamma\cup\{\beta\})$ implies $\beta\not\in\spa(\Gamma)$. Let $\Sigma=\Omega-\spa(\Gamma)$ (maximal Jordan set avoiding $\Gamma$). Let $\Sigma_1=\Omega-\spa(\Gamma\cup\{\beta\})$. So $\Sigma_1\subset \Sigma$. Now $\gamma\in\spa(\Gamma\cup\{\beta\})-\spa(\Gamma)$. So $\gamma\in\Sigma-\Sigma_1$. Choose $g\in G_{(\Omega-\Sigma)}$ such that $\beta g=\gamma$ since $\Gamma$ is a Jordan set. If $\Sigma_1g\neq\Sigma_1$, then (since $\Sigma_1g\cap\Sigma_1\neq\emptyset$) $\Sigma_1g\cup\Sigma_1$ is a Jordan set avoiding $\Gamma$ and $\gamma$. So $\Sigma_2=\Sigma_1g\cup\Sigma_1\supset \Sigma_1$. Then $\Sigma_2 g^{-1}\supset\Sigma_1$ and avoids $\beta$ and $\Gamma$, contradicting maximality of $\Sigma_1$.

Therefore $\Sigma_1g=\Sigma_1$, so $\spa(\Gamma\cup\{\gamma\})=\spa(\Gamma\cup\{\beta\})$.
\end{proof}

Let $(\alpha_1,\dots,\alpha_k)$ and $(\beta_1,\dots,\beta_k)$ be sequences of independent points.

\medskip

\noindent\textbf{Inductive hypothesis}. There exists a mapping $\alpha_i\to\beta_i$ for $1\leq i\leq k-1$. Then $\alpha_kg,\beta_k$ both lie outside $\spa(\{\beta_1,\dots,\beta_{k-1})$. But
\[ \Omega-\spa(\{\beta_1,\dots,\beta_{k-1}\})\]
is a Jordan set, so there exists $h$ fixing $\spa\{\beta_1,\dots,\beta_k\}$ pointwise and mapping $\alpha_kg\mapsto \beta_k$. Then $gh$ maps $(\alpha_1,\dots,\alpha_k)\mapsto (\beta_1,\dots,\beta_k)$.

\begin{thm}\label{39.2} Let $\Omega$ be finite and $G$ a primitive Jordan group on $\Omega$. Then we can identify $(\Omega,G)$ as follows:
\begin{enumerate}
\item[(1)] $\Omega=\PG(d,q)$ ($d\geq 2$, $q$ a prime power), $\PSL(d+1,q)\leq G\leq \PGammaL(d+1,q)$.
\item[(2)] $\Omega=\AG(d,q)$, ($d\geq 2$, $q$ a prime power) $\ASL(d,q)\leq G\leq \AGammaL(d,q)$.
\item[(3)] $\Omega=\PG(3,2)$, $G\cong A_7$ (alternating group).
\item[(4)] $\Omega=\AG(4,2)$, $n=16$, $G=(\Z_2)^4.A_7$.
\item[(5)] $\Omega =$ Witt's Steiner system on $22$ points, $G=M_{22}$.
\item[(6)] $\Omega =$ Witt's Steiner system on $23$ points, $G=M_{23}$.
\item[(7)] $\Omega =$ Witt's Steiner system on $24$ points, $G=M_{24}$.
\end{enumerate}
\end{thm}
\begin{proof} A finite primitive Jordan group is $2$-transitive (Jordan, 1877). Now if $G$ is a finite $2$-transitive group and $H$ is a minimal normal subgroup of $G$, then either
\begin{enumerate}
\item[(1)] $H$ is elementary abelian of order $p^m$, and then $G\leq \AGL(m,p)$;
\item[(2)] $H$ is simple and $G\leq \Aut(H)$ (Jordan or Burnside).
\end{enumerate}

Now apply the classification theorem for finite simple groups\footnote{Peter J.\ Cameron, \lq Finite permutation groups and finite simple groups\rq,
\emph{Bull.\ London Math.\ Soc.}, 13 (1981), 1--22.}

Case-by-case examination.
\end{proof}

\begin{thm}\label{39.3} Suppose that $\Omega$ is infinite and $G$ is primitive and there are cofinite Jordan sets with arbitrarily large complements. Then one of the following holds:
\begin{enumerate}
\item[(1)] $G$ is highly transitive;
\item[(2)] there is a prime power $q$ and a cardinal number $d$  such that $\Omega$ can be identified with $\PG(d,q)$ and $G\leq \PGammaL(d+1,q)$;
\item[(3)] there is a prime power $q$  and a cardinal number $d$ such that $\Omega$ can be identified with $\AG(d,q)$ and $G\leq \AGammaL(d,q)$.
\end{enumerate}
\end{thm}
\begin{proof}
Let $\Gamma$ be a finite set, $\Delta=\spa(\Gamma)$. Then $G^\Delta_{\{\Delta\}}$ is a finite $2$-transitive Jordan group by  \ref{38.4}. We know the finite Jordan groups. Only projective and affine geometries arise, and it is a triviality to see that they fit together.
\end{proof}

\noindent
\textbf{Commentary:}

\begin{note}
Theorem  \ref{39.2} was proved independently by Kantor (1985)\footnote{William M.\ Kantor, \lq Homogeneous designs and geometric lattices\rq, \emph{J.\ Combin.\ Theory Ser.\ A}, 38 (1985), 66--74.}, Cherlin (1985)\footnote{ G.\ Cherlin, L.\ Harrington and A.\ H.\ Lachlan, \lq $\aleph_0$-categorical, $\aleph_0$-stable structures\rq,
\emph{Ann.\ Pure Appl.\ Logic}, 28 (1985), 103--135.} , $\Pi$MN  (1985),.... Cherlin and $\Pi$MN did it to get Theorem  \ref{39.3}.
\end{note}

\begin{note}
$\Pi$MN's reason for proving  \ref{39.3} (for $|\Omega|=\aleph_0$) was to get:\end{note}

\begin{cor}\label{39.4}  S Let $G$ be a primitive group on the countable set $\Omega$. If $G$ has no countable orbits on the set of moieties of $\Omega$, then
\begin{enumerate}
\item[(1)] $G$ is highly transitive;
\item[(2)] $G\leq \PGammaL(\aleph_0,q)$; or
\item[(3)] $G\leq \AGammaL(\aleph_0,q)$.
\end{enumerate}
\end{cor}
Suppose that $G$ has no countable orbits on moieties. Let $\Gamma$ be any finite set.  Then $G_{(\Gamma)}$ has index $\aleph_0$ in $G$. Therefore $G_{(\Gamma)}$ does not stabilize any moiety. So $G_{(\Gamma)}$ has a cofinite orbit $\Sigma$. Then $\Sigma$ is an orbit of $G_{(\Omega-\Sigma)}$. Thus $\Sigma$ is a cofinite Jordan set.

\begin{note} Questions. Can it be done without the classification of finite simple groups? Answer: Yes. David Evans\footnote{David M.\ Evans, \lq Homogeneous geometries\rq, \emph{Proc.\ London Math.\ Soc.\ (3)}, 52 (1986), 305--327.}  works for finite groups $\dim\geq 23$ (possibly down to $\dim\geq 10$) Yes! Zil'ber\footnote{B.\ I.\ Zil'ber, \lq Totally categorical structures and combinatorial geometry\rq, \emph{Dokl.\ Akad.\ Nauk SSSR}, 259 (1981), 1039--1041.}  $\dim\geq 7$.
\end{note}

\chapter{Stable and Semistable Jordan Sets}

Define $\Sigma$ to be \emph{stable}\footnote{The material in this lecture was published in S.\ A.\ Adeleke and Peter M.\ Neumann, \lq Primitive permutation groups with primitive Jordan sets\rq, \emph{J.\ London Math.\ Soc.\ (2)}, 53 (1996), 209--229.  There the word \emph{catenary} is used instead of \emph{stable}.}  if $\forall g\in G$, $\Sigma\subseteq \Sigma g$ or $\Sigma g\subseteq \Sigma$ and \emph{semistable} if $\forall g\in G$, $\Sigma\cap\Sigma g\neq \emptyset$.

\begin{thm}\label{40.1} Suppose that $G$ is primitive on $\Omega$ and that $\Sigma$ is a non-empty proper stable subset of $\Omega$. Then
\begin{enumerate}
\item[(1)] there is a $G$-invariant dense linear order on $\Omega$;
\item[(2)] if $\Sigma$ is a Jordan set then $G$ is highly homogeneous.
\end{enumerate}
\end{thm}
\begin{proof} (1) Define $\alpha\leq\beta\Leftrightarrow(\forall g\in G) (\beta\in \Sigma g\Rightarrow \alpha\in\Sigma g)$. Then $\alpha\leq \alpha$, $\alpha\leq\beta\wedge\beta\leq \gamma\Rightarrow \alpha\leq \gamma$. Define $\alpha\mathrel{\rho}\beta\Leftrightarrow \alpha\leq\beta\wedge\beta\leq\alpha$. Then $\rho$ is a $G$-congruence and proper. Therefore $\rho$ is equality. Thus $\alpha\leq \beta\wedge \beta\leq\alpha\Rightarrow\alpha=\beta$. Suppose that $\alpha\not\leq\beta$. Then $(\exists h)(\beta\in\Sigma h \wedge \alpha\not\in\Sigma h)$. Suppose that $\alpha\in\Sigma g$. Either $\Sigma g\subseteq \Sigma h$ or $\Sigma h\subseteq \Sigma g$. But $\Sigma g\not\subseteq \Sigma h$ since $\alpha\not\in\Sigma h$. Therefore $\Sigma h\subseteq \Sigma g$, so $\beta\in\Sigma g$.
Thus $\alpha\not\leq \beta\Rightarrow \beta\leq\alpha$ and therefore $\leq$ is a linear order (dense).

(2) Suppose $\alpha_1<\cdots<\alpha_k$ and $\beta_1<\cdots<\beta_k$. Inductive hypothesis: there exists $g_1$ mapping $\alpha_i\mapsto \beta_i$ for $2\leq i\leq k$. Then let $\gamma=\alpha_1g$, so that $\gamma<\beta_2$ and $\beta_1<\beta_2$. Thus there are translates $\Sigma_1,\Sigma_2$ of $\Sigma$ such that $\gamma\in\Sigma_1$, $\beta_2\not\in\Sigma_1$, and $\beta_1\in\Sigma_2$, $\beta_2\not\in\Sigma_2$. One of these contains the other so we get $\Sigma'$ a translate of $\Sigma$ such that $\gamma,\beta_1\in\Sigma'$, $\beta_2,\dots,\beta_k\not\in\Sigma'$. Since $\Sigma'$ is a Jordan set we can map $\gamma\mapsto\beta_1$ while fixing $\beta_2,\dots,\beta_k$. By induction the result now follows.
\end{proof}

\begin{thm}\label{40.2} Suppose that $G$ is primitive on $\Omega$ and that there is a Jordan set $\Sigma\neq\Omega$ which is semistable. Then either $G$ is $2$-transitive, or there is a $G$-invariant dense linear order on $\Omega$ and $G$ is highly homogeneous.
\end{thm}

Typical example: $G=\Aut(\Q,\leq)$ and $\Sigma$ is an open left section (i.e., all elements less than a given rational number).

\begin{cor}\label{40.3} Suppose that $G$ is a primitive Jordan subgroup of a `geometric group' ($\PGammaL(d,F)$, $\AGammaL(d,F)$, \dots ). Then $G$ is $2$-transitive.
\end{cor}

\begin{proof}[Proof of 40.2]
Let $\Sigma_0$ be our semistable Jordan set. Let $\alpha\in\Omega-\Sigma_0$. Let $\Sigma$ be a maximal Jordan set such that $\Sigma_0\subseteq \Sigma$ and $\alpha\not\in\Sigma$. Then $\Sigma$ is also semistable. Note that $\Sigma$ is a $G_\alpha$-orbit. Let $g\in G$ and suppose that $\alpha g\not\in\Sigma$. If $\alpha g=\alpha$, then $\Sigma g=\Sigma$. If $\alpha g\neq \alpha$, then $\Sigma\subset \Sigma g$.

\medskip

\noindent\textbf{Case 1.} For all $g\in G$ either $\alpha g\not\in\Sigma$ or $\alpha g^{-1}\not\in \Sigma$. But then for all $g\in G$ either $\Sigma\subseteq \Sigma g$ or $\Sigma\subset \Sigma g^{-1}$, so $\Sigma$ is stable and  \ref{40.1} applies.

\medskip

\noindent\textbf{Case 2.} $(\exists g_0\in G)(\alpha g_0\in\Sigma\text{ and }\alpha g_0^{-1}\in\Sigma)$. We prove that $\{\alpha\}=\Omega-\Sigma$. Suppose there exists $\beta\neq \alpha$, $\beta\not\in\Sigma$. Choose $g_1:\alpha\mapsto \beta$. Thus $\Sigma\subseteq \Sigma g_1$ so $\alpha\in\Sigma g_1$. So $\alpha=\gamma g_1$, $\gamma\in\Sigma$. Choose $h\in G_\alpha$ mapping $\alpha g_0$ to $\gamma$. Let $g_2=g_0h$, $\delta=\alpha g_2^{-1}$. Then $\delta\in\Sigma$ (since $\delta=\alpha g_0^{-1}$). And $\delta g_2g_1=\alpha g_1=\beta$, $\alpha g_2g_1=\alpha g_0hg_1=\gamma g_1=\alpha$. This contradicts $g_2g_1\in\ G_\alpha$, $\delta\in\Sigma$ which is a $G_\alpha$-orbit, $\beta\not\in\Sigma$.
\end{proof}

\newpage 

\chapter{Semilinear Order Relations}

M. Droste (1985)\footnote{Manfred Droste,  \lq Structure of partially ordered sets with transitive automorphism groups\rq, \emph{Mem.\ Amer.\ Math.\ Soc.}, 57 (1985), no. 334}.

\begin{defn}
A partially ordered set $(\Lambda,\leq)$ is said to be \emph{upper semilinearly ordered} if
\begin{enumerate}
\item $(\alpha\leq \beta\wedge \alpha\leq \gamma)\Rightarrow (\beta\leq \gamma\vee \gamma\leq\beta)$,
\item $(\forall \alpha,\beta)(\exists \gamma)(\alpha\leq\gamma\wedge \beta\leq \gamma)$.
\end{enumerate}
\end{defn}

Look up and get linear order though may get branching down below.

\medskip
\noindent\emph{Note:} Definition of lower semi-linear order is dual.

\medskip

A subset $\Sigma$ of $\Lambda$ is an \emph{upper section} if it is a filter and has a lower bound. The upper sections are semilinearly ordered by inclusion and form a Dedekind complete semilinearly ordered set.

Define $\Lambda_-(\Sigma)=\{\gamma\mid\gamma<\sigma\text{ for all }\sigma\in\Sigma\}$.

Define a binary relation $\rho$ on $\Lambda_-(\Sigma)$ by
\[ (\gamma_1,\gamma_2)\in\rho\Leftrightarrow (\exists \gamma\in\Lambda_-(\Sigma))(\gamma_1\leq \gamma\wedge \gamma_2\leq \gamma).\]

Then $\rho$ is an equivalence relation. The equivalence classes are themselves semilinearly ordered: the \emph{cones} below $\Sigma$. Call $\Sigma$ a \emph{ramification point} of ramification index $s$ if there are $s\geq 2$ cones below $\Sigma$. Write $\alpha||\beta$ if $\alpha\not\leq \beta$ and $\beta\not\leq \alpha$. Note that if $\alpha||\beta$ and $\Sigma=\{\gamma\mid \alpha\leq\gamma\wedge \beta\leq\gamma\}$, then $\Sigma$ is an upper section and it is a ramification point. Conversely, all ramification points are of this form.

We say that $\Lambda$ is \emph{dense} if $\alpha<\beta\Rightarrow (\exists\gamma)(\alpha<\gamma<\beta)$ (and also $(\forall \alpha)(\exists\gamma)(\alpha<\gamma)$). We say that $\Lambda$ is of \emph{positive type} if every pair of elements has a least upper bound, so that $\Lambda$ forms a semilattice. We say that $\Lambda$ is of \emph{negative type} if no two incomparable elements has a least upper bound.

Let us (temporarily) say that $\Lambda$ is \emph{structurally homogeneous} if:
\begin{enumerate}
\item[(1)] $\Lambda$ is dense;
\item[(2)] $\Lambda$ is of positive or negative type;
\item[(3)] all ramification indices are the same.
\end{enumerate}

Let us say that $(\Omega,R)$ is \emph{relatively $k$-homogeneous} if, whenever $\Delta_1,\Delta_2$ are isomorphic $k$-subsets of $\Omega$ then there exists $f\in\Aut(\Omega,R)$ such that $f:\Delta_1\to\Delta_2$; \emph{relatively $k$-transitive} if, whenever $f_0:\Delta_1\to\Delta_2$ is an isomorphism of $k$-element subsets there exists $f\in\Aut(\Omega,R)$ with $f\ha \Delta_1=f_0$.

\medskip

\noindent\emph{Warning:} These are often used the other way round.

\begin{thm}\label{41.1} If $(\Lambda,\leq)$ is a relatively $2$-homogeneous semilinearly ordered set then $\Lambda$ is structurally homogeneous. If $\Lambda$ is countable and structurally homogeneous then $(\Lambda,\leq)$ is relatively $2$-homogeneous.
\end{thm}

\begin{wrapfigure}{r}{4cm}   
\begin{tikzpicture}
\node(x) at (-0.3,0) {};
\draw[->] (3,3) -> (0,0);
\draw[->] (2,2) -> (2,0.25);
\draw[->] (2,1) -> (3.5,0.125);
\draw[->] (3.2,0.3) -> (2.7,-0.3);
\node(q) at (2.4, 2.8) {$\Q$};
\node(q1) at (1.6, 2.0) {$q_1$};
\node(q2) at (1.75, 1.0) {$q_2$};
\node(q3) at (3.47, 0.37) {$q_3$};
\node(u1) at (2.25, 1.5) {$u_1$};
\node(u2) at (2.9, 0.75) {$u_2$};
\node(u3) at (3.2, -0.1) {$u_3$};
\end{tikzpicture}
\end{wrapfigure}
\noindent
{\bf Example 41.3.}
Construction of a countable dense semilinear order of positive type and of ramification index $s$, where $2\leq s\leq \aleph_0$. We denote this by $\Lambda(\Q,s,+)$, the set of non-empty words $q_1u_1q_2u_2\ldots q_{k-1}u_{k-1}q_k$ where $u_i\in U$, $U$ a set such that $|U|=s-1$, $q_i\in\Q$, $q_1>q_2>\cdots >q_k$. Define
\[ q_1u_1\ldots q_{k-1}u_{k-1}q_k\leq q_1'u_1'\ldots q_{\ell-1}'u_{\ell-1}'q_\ell'\]
if and only if $\ell \leq k$ and $q_1u_1\ldots q_{\ell-1}u_{\ell-1}q_\ell=q_1'u_1'\ldots q_{\ell-1}'u_{\ell-1}$, with $q_\ell\leq q_{\ell}'$.

\medskip

\noindent\textbf{Exercise.} $(\Lambda,\leq)$ is semilinearly ordered of positive type, dense, ramification index $s$.

\chapter{Automorphism Groups of Semilinear Relations}

\begin{fact}\label{42.1} The theory of dense upper semilinear orders of positive type and branching number $s$ (for a given value of $s$, $2\leq s\leq\aleph_0$) is $\aleph_0$-categorical.
\end{fact}
\begin{cor}\label{42.2}
If $G=\Aut(\Lambda)$ where $\Lambda$ is as above, then $G$ has finitely many orbits on $\Lambda^k$ for every $k$ (Ryll--Nardzewski).
\end{cor}

Take $\Lambda$ to be a countable dense upper semilinear order of positive type and ramification number $s\geq 2$ at each ramification point. Let $G=\Aut(\Lambda)$.

\begin{thm}\label{42.3} $G$ is a primitive permutation group of rank $4$.
\end{thm}
\begin{proof} Let $\alpha\in\Lambda$. Define $\Sigma_\alpha=\{\gamma\mid\alpha<\gamma\}$, $\Lambda_\alpha=\{\gamma\mid \gamma<\alpha\}$, $\Gamma_\alpha=\{\gamma\mid\gamma||\alpha\}$. The $G_\alpha$-orbits are $\{\alpha\}$, $\Sigma_\alpha$, $\Lambda_\alpha$, $\Gamma_\alpha$. Since cones are isomorphic to $\Lambda$, $G$ is transitive on $\Lambda$ because
\begin{enumerate}
\item[(1)] there are $s$ cones below $\alpha$ and $s$ cones below $\beta$ all isomorphic, so there exists $\phi:\Lambda_\alpha\mathop{\to}\limits^\sim \Lambda_\beta$,
\item[(2)] $\Sigma_\alpha$, $\Sigma_\beta$, are dense open linearly ordered hence isomorphic to $\Q$, so there is an isomorphism $\psi:\Sigma_\alpha\to\Sigma_\beta$,
\item[(3)] for each $\gamma\in\Sigma_\alpha$ there are $s-1$ cones below $\gamma\phi$ not containing $\beta$. So there is an isomorphism \[\theta_\gamma:\Lambda_\gamma-(\text{cone that contains } \alpha)\to\Lambda_{\gamma\phi}-(\text{cone that contains } \beta)\]
\end{enumerate}
Putting this all together forms an isomorphism $\Lambda\to\Lambda$ mapping $\alpha\mapsto\beta$, so $G$ is transitive. Also $G_\alpha$ is transitive on $\Lambda_\alpha$, $\Sigma_\alpha$ and $\Gamma_\alpha$.

Let $\rho$ be any $G$-congruence, $\rho\neq =$. Suppose $\alpha\neq\beta$, $\alpha\equiv\beta\bmod \rho$, so $\rho$ is universal.
\end{proof}

\begin{thm}\label{42.4} Let $H=G_\alpha$. Then $H=H_-\times H_+$ where $H_-=\Aut(\Lambda_\alpha)\cong  G\Wrr \Sym(s)$, and $H_+\cong G\Wrr \Sym(s-1) \Wrr \Aut(Q,\leq)$.
\end{thm}
\begin{proof} Define $H_+=G_{(\Delta_\alpha\cup\{\alpha\})}$, $H_-=G_{(\Sigma_\alpha\cup\Gamma_\alpha\cup\{\alpha\})}\cong G\Wrr \Sym(s)$.
\end{proof}
\noindent \emph{Problem.} What is the rate of growth of the number of orbits of $G$ on $\Lambda^k$ or $\Lambda^{\{k\}}$?\footnote{This problem is discussed in Peter J.\ Cameron, \lq Some treelike objects\rq, \emph{Quart.\ J.\ Math.\ Oxford Ser.\  (2)} 38 (1987), 155–183 and  H.\ D.\ Macpherson, \lq Orbits of infinite permutation groups\rq, \emph{Proc.\ London Math.\ Soc.\ (3)} 51 (1985), 246--284 (in particular Example 2.4).    For more recent work on this topic see Samuel Braunfeld, \lq 
Monadic stability and growth rates of $\omega$-categorical structures\rq,
\emph{Proc.\ Lond.\ Math.\ Soc.\ (3)} 124 (2022), 373--386.} 

\begin{thm} Let $\Sigma$ be an upper section not of the form $\Sigma_\alpha$.
\begin{enumerate}
\item If $\Gamma$ is a union of cones below $\Sigma$, then $\Gamma$ is a Jordan set for $G$.
\item All Jordan sets for $G$ are of this form.
\item If $\Gamma$ is a cone below $\Sigma$ then $\Gamma$ is a primitive Jordan set.
\item All primitive Jordan sets are of this form.
\end{enumerate}
\end{thm}

\chapter{\texorpdfstring{$C$}{C}-relations}

A ternary relation $G$ is a \emph{$C$-relation} if
\begin{enumerate}
\item[(C1)] $C(\alpha;\beta,\gamma)\Rightarrow C(\alpha;\gamma,\beta)$,
\item[(C2)] $C(\alpha;\beta,\gamma)\Rightarrow \lnot C(\beta;\alpha,\gamma)$,
\item[(C3)] $C(\alpha;\beta,\gamma)\Rightarrow C(\alpha;\gamma,\delta)\vee C(\delta;\beta,\gamma)$,
\item[(C4)] $\alpha\neq\beta\Rightarrow C(\alpha;\beta,\beta)$.
\item[(C5)] $(\forall \beta,\gamma)$ $(\exists \alpha)$ $C(\alpha;\beta,\gamma)$,
\item[(C6)] $\alpha\neq \beta\Rightarrow (\exists\delta\neq \beta) C(\alpha;\beta,\delta)$.
\end{enumerate}
It is \emph{dense} if also
\begin{enumerate}
\item[(C7)] $C(\alpha;\beta,\gamma)\Rightarrow\exists\delta(C(\alpha;\beta,\delta)\wedge C(\delta;\beta,\gamma))$.
\end{enumerate}

\begin{wrapfigure}{r}{3cm}  \begin{tikzpicture}[node distance=1cm,line width=1pt]
\node(A) at (0,-0.2) {$\alpha$};
\node(L)  at (1,-0.2){$\beta$};
\node(K) at (1.5,-0.2){$\gamma$};
\draw [<-] (0,0) -- (1,2); 
\draw [->] (0.75,1.5) -- (1,0); 
\draw [->] (0.875, 0.75) -- (1.5,0); 
\end{tikzpicture}\end{wrapfigure}

\ \newline
\emph{Example.}  Let $(\Lambda,\leq)$ be a semilinear order. Let $\Omega$ be the set of maximal chains in $\Lambda$.
Define $C(\alpha;\beta,\gamma)\Leftrightarrow \alpha\cap\beta=\alpha\cap\gamma\neq \beta\cap\gamma$.   To guarantee (C5) and (C6) we have to impose some weak extra conditions on $(\Lambda, \leq)$.
\ \newline

\begin{lem}\label{43.1}
\begin{enumerate}
\item $\lnot C(\alpha;\alpha,\beta)$.
\item $\lnot C(\beta;\gamma,\alpha)\wedge \lnot C(\gamma;\delta,\alpha)\Rightarrow \lnot C(\beta;\delta,\alpha)$.
\item $C(\alpha;\beta,\gamma)\wedge C(\alpha;\gamma,\delta)\Rightarrow C(\alpha;\beta,\delta)$.
\end{enumerate}
\end{lem}

Now fix attention on $\alpha$, define relations $R_\alpha$, $S_\alpha$ on $\Omega-\{\alpha\}$ by
\begin{align*} R_\alpha(\beta,\gamma)&\Leftrightarrow \lnot C(\beta;\gamma,\alpha)\wedge \lnot C(\gamma;\beta,\alpha),
\\ S_\alpha(\beta,\gamma)&\Leftrightarrow C(\alpha;\beta,\gamma).
\end{align*}

\begin{cor}\label{43.2} $R_\alpha$, $S_\alpha$ are equivalence relations on $\Omega-\{\alpha\}$. Also $S_\alpha$ refines $R_\alpha$.
\end{cor}

\begin{wrapfigure}{r}{3cm}\begin{tikzpicture}[node distance=1cm,line width=1pt]
\node(A) at (0,-0.2) {$\alpha$};

\node(L)  at (1,-0.2){$\beta$};
\node(K) at (1.5,-0.2){$\gamma$};
\node(K) at (2,-0.2){$\delta$};
\node(M)   at (1.5, 3){};
\draw [<-] (0,0) -- (1,2); 
\draw [->] (0.75,1.5) -- (1,0); 
\draw [->] (0.75,1.5) -- (1.5,0); 
\draw [->] (1.1225, 0.775) -- (2,0); 
\end{tikzpicture}\end{wrapfigure}
\ \newline\ \newline \emph{Example (continuation).}   In the case of chains $R_\alpha(\beta, \gamma)$ says that  $\beta,\gamma$ with $\alpha$ determine the same ramification point. And $S_\alpha(\gamma, \delta)$ indicates that $\gamma, \delta$ lie in the same cone of the ramification point with $\alpha$

\bigskip
The $R_\alpha$-classes in $\Omega-\{\alpha\}$ are linearly ordered
\[ R_\alpha(\beta)\leq R_\alpha(\gamma)\Leftrightarrow \lnot C(\beta;\gamma,\alpha).\]

Let $(\Omega,C)$ be a set with $C$-relation. Define binary relations $\preceq$, $\equiv$ on $\Omega^{\{2\}}$ by $\{\alpha,\beta\}\preceq\{\gamma,\delta\}\Leftrightarrow \lnot C(\alpha;\gamma,\delta)\wedge \lnot C(\beta;\gamma,\delta)$.
Define
\[\{\alpha,\beta\}\equiv\{\gamma,\delta\}\Leftrightarrow \{\alpha,\beta\}\preceq \{\gamma,\delta\}\text{ and }\{\gamma,\delta\}\preceq \{\alpha,\beta\}.\]

\begin{thm}\label{43.3} If $\Lambda=\Omega^{\{2\}}/\equiv$ and $\leq$ is induced by $\preceq$ then $(\Lambda,\leq)$ is a semilinearly ordered set. The map
\[ \alpha\mapsto \{\{\alpha,\beta\}\mid\beta\in\Omega-\{\alpha\}\}/\equiv\]
is an identification of $\Omega$ with a covering set of maximal chains. Moreover $(\Lambda,\leq)$ is of positive type.
\end{thm}

Define the ramification order of a pair $\{\alpha,\beta\}$ to be $ 1+\#\big(S_\alpha\text{-classes in }R_\alpha(\beta)\big)$.

\begin{thm}\label{43.4} The theory of dense $C$-relations with ramification order $s$ at each `branch point' is $\aleph_0$-categorical (for any $s$, $2\leq s\leq \aleph_0$).
\end{thm}
\begin{thm}\label{43.5} Let $(\Omega,C)$ be the countable dense $C$-relation with ramification order $s$ at each pair $\{\alpha,\beta\}$, $\alpha\neq\beta$. Let $G=\Aut(\Omega,C)$. Then $G$ is $2$-transitive. Also $G_\alpha$ acting on $\Omega-\{\alpha\}$ is isomorphic to $G\Wrr \Sym(s-1)\Wrr\Aut(\Q,\leq)$ acting on $\Omega\times(s-1)\times\Q$ in the natural way.
\end{thm}

\noindent \textbf{Observation.}  $G\leq \Aut(\Lambda,\leq)$.

\smallskip

There are four types of Jordan sets.
\begin{enumerate}
\item[I]\begin{enumerate}
\item[(i)] Set of all $\gamma$ in a cone ($S_\alpha(\gamma)$ for some $\alpha$)
\item[(ii)] Union of two or more such cones at the same branch point $\bigcup_\gamma S_\alpha(\gamma)$.
\end{enumerate}
\item[II] There is a proper open upper section $\Phi$ in $\Lambda$ and $\Gamma=\{\gamma\mid\Phi\subseteq\gamma\}$.
\item[III] There is an open convex linear bounded subset $\Phi$ of $\Lambda$ and $\Gamma=\{\gamma\mid R_\alpha(\gamma)\subseteq\Phi\}$.
\end{enumerate}

\begin{thm}\label{43.6} If $(\Omega,C)$ is relatively $3$-homogeneous then these are Jordan sets and all are of this form. Those of type I(i), II are $2$-transitive. The others are imprimitive.
\end{thm}

\chapter{More on \texorpdfstring{$C$}{C}-relations}

\setcounter{note}{0}
Let $\Omega$ be as usual, $\{\rho_i\}_{i\in I}$ a chain of proper equivalence relations (all classes non-trivial) on $\Omega$ such that $\bigcap\rho_i$ is trivial and $\bigcup\rho_i$ is universal. Define a ternary relation:
\[ C(\alpha;\beta,\gamma)\Leftrightarrow (\exists i) \alpha\not\equiv\beta\equiv\gamma\bmod \rho_i.\]

\begin{claim}\label{44.1} This is a $C$-relation.
\end{claim}

\begin{proof}
It is clear the $C(\alpha;\beta, \gamma)\Leftrightarrow C(\alpha; \gamma, \beta)$ and that 
$C(\alpha;\beta, \gamma)\Leftrightarrow \neg C(\beta; \alpha, \beta)$.

Suppose $C(\alpha; \beta, \gamma)$ and that $\alpha\not\equiv\beta\equiv\gamma\bmod \rho_i$.  If $\delta\equiv\gamma \bmod \rho_i$ then $C\alpha; \delta, \gamma)$ and if $\delta\not\equiv\gamma \bmod \rho_i$ then $C(\delta; \beta, \gamma)$.  

Suppose $\alpha\neq \beta$.  There exists $i$ such that $\alpha\not\equiv \beta \bmod \rho_i$.  Chose $\gamma\neq \beta$ in $\rho_i(\beta)$.  Then $C(\alpha; \beta, \gamma)$.

Let $\beta, \gamma$ be in $\Omega$.   Then there exists $i$ such that $\beta\equiv \gamma \bmod \rho_i$.  Then there exists $\alpha$ such that $\alpha\not\equiv \beta \bmod \rho_i$ and then $C(\alpha;\beta, \gamma)$.
\end{proof}

\begin{note}
$\Lambda=\{\rho_i(\gamma)\mid i\in I,\gamma\in\Omega\}$. Then $\Lambda$ is a semilinearly ordered subset of $P(\Omega)$. This semilinear order is usually the semilinear order determined by $(\Omega,C)$.
\end{note}

\begin{note} We can use  \ref{44.1} to construct the relatively homogeneous dense $s$-branched countable $C$-relation as follows. Let $S$ be a set with $(s+1)$ members, one of which is $\ast$. Let
\[ \Omega=\{f:\Q\to S\mid \supp(f)\text{ is finite}\},\]
where $\supp(f)$ is
\[ \{q\in\Q\mid f(q)\neq \ast\}.\]
Define $\rho_q$ by
\[ f\equiv g\bmod\rho_q\Leftrightarrow f(r)=g(r)\;\;\forall r>q.\]
The family $\{\rho_q\}_{q\in\Q}$ of equivalence relations determines a homogeneous $C$-relation.
\end{note}

Return to permutation groups $G$ on $\Omega$. If $\Gamma\subseteq \Omega$, then we say $\Gamma$ is \emph{highly atypical} if
\begin{enumerate}
\item $|\Gamma|>1$, $\Gamma\neq \Omega$,
\item $(\forall g\in G)(\Gamma g\cap\Gamma=\emptyset \vee \Gamma g\subseteq \Gamma \vee \Gamma\subseteq\Gamma g)$.
\end{enumerate}

\begin{lem}\label{44.2} Suppose that $G$ is primitive on $\Omega$ and that $\Gamma$ is a highly atypical subset of $\Omega$. Then, for all $\alpha,\beta\in\Omega$, there exists $g\in G$ such that $\alpha,\beta\in\Gamma g$.
\end{lem}
\begin{proof} Define a relation $\rho$ on $\Omega$ by
\[ (\alpha,\beta)\in\rho\Leftrightarrow (\exists g)(\alpha,\beta\in\Gamma g).\]
Then $\rho$ is a $G$-invariant relation.  

Since $G$ is transitive we see that $(\alpha, \alpha)\in \rho$.

Symmetric since definition symmetric.

Suppose that $\alpha, \beta\in \Gamma g$, $\beta, \gamma\in \Gamma h$.  Then $\Gamma g\cap \Gamma h\neq \emptyset$ so $\Gamma g\subseteq \Gamma h$ or $\Gamma h\subseteq \Gamma g$.  Thus $(\alpha, \gamma)\in \rho$.  

Thus $\rho$ is a $G$-congruence.  Non-trivial, so it must be universal.
\end{proof}


\begin{thm}\label{44.3} Suppose that $G$ is primitive and has a highly atypical subset $\Sigma_0$. Then either there is a $G$-invariant dense semilinear order on $\Omega$ or there is a $G$-invariant dense $C$-relation on $\Omega$. (Note that the semilinear order may in fact be linear.)
\end{thm}
\begin{proof} Recall from Lecture 22 we say that $\Sigma_0$ is idealistic if it separates pairs but not ordered pairs. Suppose $\Sigma_0$ is idealistic.  (So $G$ is not strongly primitive.) Define
\[ \alpha\leq \beta\Leftrightarrow (\beta\in\Sigma_0g\Rightarrow \alpha\in\Sigma_0g).\]
This is a partial order. Suppose that $\alpha\leq \beta$ and $\alpha\leq \gamma$. Suppose also that $\beta\not\leq \gamma$. There is a translate $\Sigma_1$ of $\Sigma_0$ such that $\gamma\in\Sigma_1$, $\beta\not\in\Sigma_1$. But $\alpha\leq\gamma$ so $\alpha\in\Sigma_1$. Let $\Sigma$ be any translate of $\Sigma_0$ containing $\beta$. Then $\alpha\in\Sigma$. Now $\Sigma\cap\Sigma_1\neq\emptyset$. So $\Sigma\subseteq\Sigma_1$ or $\Sigma_1\subseteq \Sigma$. But $\Sigma\not\subseteq\Sigma_1$. Therefore $\Sigma_1\subseteq\Sigma$ so $\gamma\in\Sigma$ and $\gamma\leq \beta$. Prove also the existence of common upper bounds.

Now suppose that $\Sigma_0$ is not idealistic. Then for all $\alpha,\beta\in\Omega$ if $\alpha\neq\beta$ then $(\exists g\in G)(\alpha\in\Sigma_0g\wedge \beta\not\in\Sigma_0g)$. Define $C(\alpha;\beta,\gamma)\Leftrightarrow$ there exists a translate $\Sigma$ of $\Sigma_0$ such that $\beta,\gamma\in\Sigma$, $\alpha\not\in\Sigma$. This is a $C$-relation. Clearly condition (C1) holds.  If both $C(\alpha; \beta, \gamma)$ and $(C(\beta;\alpha, \gamma)$ then there would be translates of $\Sigma_1$ and $\Sigma_2$ of $\Sigma_0$ such that $\beta, \gamma\in \Sigma_1$, $\alpha\nin \Sigma_1$ and $\alpha, \gamma\in \Sigma_2$, $\beta\nin\Sigma_2$.  This contradicts the assumption that $\Sigma_0$ is highly atypical.  For (C3) we assume that $C(\alpha;\beta, \gamma)$.  Let $\Sigma$ be a translate of $\Sigma_0$ such that $\beta, \gamma\in \Sigma$ but $\alpha\nin \Sigma$,  If $\delta in \Sigma$ then $C(\alpha; \gamma, \delta)$ and if $\delta \nin \Sigma$ then $C(\delta; \beta, \gamma)$.   Conditions (C4), (C5) and (C6) are obvious.  

Denseness comes from primitivity. Suppose $\Lambda$ is not dense. Then there exists $\alpha,\alpha^+$, with $\alpha^+$ immediately above $\alpha$. Since $G$ is transitive every point has an immediate successor.
Define $\alpha,\beta$ to be `at finite distance' if there exists $\gamma$ above $\alpha,\beta$  such that $[\alpha,\beta],[\beta,\gamma]$ are finite. This is a $G$-invariant non-trivial equivalence relation so is universal. Now define $\rho$ by
\[ \alpha\equiv\beta\bmod \rho\Leftrightarrow [\alpha,\gamma]+[\beta,\gamma]\text{ is even}.\]
Thus $\rho$ is a non-trivial, proper $G$-congruence, contradicting the assumption that $G$ is primitive.
\end{proof}

\chapter{\texorpdfstring{$B$}{B}- and \texorpdfstring{$D$}{D}-relations}

A ternary relation $B$ on $\Omega$ is a \emph{$B$-relation} if it satisfies
\begin{enumerate}
\item[(B1)] $B(\alpha;\beta,\gamma)$ ($\alpha$ is between $\beta$ and $\gamma$) $\Rightarrow$ $B(\alpha;\gamma,\beta)$,
\item[(B2)] $B(\alpha;\beta,\gamma)\wedge B(\beta;\alpha,\gamma)\Leftrightarrow \alpha=\beta$,
\item[(B3)] $B(\alpha;\beta,\gamma)\Rightarrow B(\alpha;\beta,\delta)\vee B(\alpha;\gamma,\delta)$,
\item[(B4)] $B(\alpha;\gamma,\delta)\wedge B(\beta;\alpha,\gamma)\Rightarrow B(\beta;\gamma,\delta)$,
\item[(B5)] $B(\alpha;\gamma,\delta)\wedge B(\beta;\gamma,\delta)\Rightarrow B(\alpha;\beta,\gamma)\vee B(\beta;\alpha,\gamma)$.
\end{enumerate}
If also
\begin{enumerate}
\item[(B6)] $\lnot B(\alpha;\beta,\gamma)\Rightarrow (\exists \delta\neq\alpha)(B(\delta;\alpha,\beta)\wedge B(\delta;\alpha,\gamma))$,
\end{enumerate}
then $B$ is a \emph{general betweenness relation}.

Denseness, being infinite in all directions, etc., can be similarly formulated as before. Likewise discreteness also (combinatorial trees).

\begin{example}\label{45.1}
Let $\Lambda$ be a semilinearly ordered set. Define $B(\alpha;\beta,\gamma)$ if
\begin{enumerate}
\item[(1)] $\beta\leq\alpha\leq\gamma$ or $\gamma\leq\alpha\leq\beta$, or
\item[(2)] $\alpha=\sup(\beta,\gamma)$, or
\item[(3)] $(\beta\leq\alpha$ and $\alpha||\gamma$) or ($\gamma\leq\alpha$ and $\alpha||\beta$).
\end{enumerate}
Check that this is a betweenness relation. If $\Lambda$ is countable $\aleph_0$-categorical then $(\Lambda,B)$ is also $\aleph_0$-categorical. And $\Aut(\Lambda,B)$ is $2$-transitive.
\end{example}

Call $\alpha,\beta,\gamma$ \emph{collinear} if $B(\alpha;\beta,\gamma)\vee B(\beta;\alpha,\gamma)\vee B(\gamma;\alpha,\beta)$. Call a set $\Gamma$ \emph{linear} if any three points are collinear and \emph{convex} if $\beta,\gamma\in\Gamma\wedge B(\alpha;\beta,\gamma)\Rightarrow \alpha\in\Gamma$. If $B(\delta;\alpha,\beta)\wedge B(\delta,\alpha,\gamma)\wedge B(\delta;\beta,\gamma)$ then call $\delta$ the \emph{centroid} of $\alpha,\beta,\gamma$. Say that $(\Omega,B)$ is of \emph{positive} type if any three points have a centroid, negative if no three non-collinear points have a centroid.

Define a subset $\Delta$ to be \emph{closed} if it is convex and co-convex (i.e., $\Omega-\Delta$ is convex), a \emph{component} if it is not the disjoint union of two non-empty closed sets. Three non-collinear points determine a branch point and branch point determines components.

So we get branching numbers.

\begin{thm} \label{45.2} The theory of dense $B$-relations (or betweenness) of positive type and given branching number $s$, $3\leq s\leq \aleph_0$ at every point is $\aleph_0$-categorical. The automorphism group of the countable structure is $2$-transitive and has primitive Jordan sets. (Similarly for negative type.)
\end{thm}
Define $\Gamma$ be be a \emph{line} in $(\Omega,B)$ if $\Gamma$ is a maximal linear set. Then $\Gamma$ is convex, carries a linear betweenness relation (restrict $B$), hence a pair of linear orderings (reverse of each other). Then define a \emph{half-line} to be a Dedekind section in one of these orderings on a line in $(\Omega,B)$. Define half-lines $\Gamma_1,\Gamma_2$ to be \emph{equivalent} if $\Gamma_1\wedge \Gamma_2$ is a half-line. Call equivalence classes of half-lines \emph{directions} or \emph{points at infinity} or `leaves' or `ends'. (Two different directions determine a unique line -- not easy.)

A \emph{$D$-relation} is a four-place relation $D$ satisfying
\begin{enumerate}
\item[(D1)] $D(\alpha,\beta;\gamma,\delta)\Rightarrow D(\beta,\alpha;\gamma,\delta)\wedge D(\alpha,\beta;\delta,\gamma)\wedge D(\gamma,\delta;\alpha,\beta)$,
\item[(D2)] $D(\alpha,\beta;\gamma,\delta)\Rightarrow \lnot D(\alpha,\gamma;\beta,\delta)$,
\item[(D3)] $D(\alpha,\beta;\gamma,\delta)\Rightarrow D(\alpha,\ep;\gamma,\delta)\vee D(\alpha,\beta;\gamma,\ep)$,
\item[(D4)] $\alpha,\beta,\gamma$ distinct $\Rightarrow$ $(\exists\delta\neq\gamma)(D(\alpha,\beta;\gamma,\delta)$.
\end{enumerate}
Axiom for denseness, etc.

\begin{wrapfigure}[5]{r}{4cm}
 \begin{tikzpicture}
      \draw [<->,domain=-45:45] plot ({cos(\x)},{1+sin(\x)});
   \draw [<->,domain=135:225] plot ({4+cos(\x)},{1+sin(\x)});
   \draw (1,1) -- (3,1);
   \node(A) at ({-0.2+cos(-45)}, {1+sin(-45)}) {$\alpha$};
    \node(B) at ({-0.2+cos(45)}, {1+sin(45)}) {$\beta$};
   \node(C) at ({4.2+cos(135)}, {1+sin(135)}) {$\gamma$};
    \node(D) at ({4.2+cos(225)}, {1+sin(225)}) {$\delta$};
\end{tikzpicture} 
\end{wrapfigure}  \ \newline
The figure shows a schematic way to think about the meaning of $D(\alpha, \beta;\gamma, \delta)$.

If $(\Omega,D)$ is a $D$-set, then $(\Omega-\{\alpha\},C)$ is a $C$-set, where $C(\beta;\gamma,\delta)\Leftrightarrow D(\alpha,\beta;\gamma,\delta)$.

If $(\Omega,D)$ is a $D$-set then there is a $B$-relation $(\Lambda,B)$ of positive type determined by $D$ (a quotient of $\Omega^{\{3\}}$) such that $\Omega$ is embedded into the set of directions of $(\Lambda,B)$.  

\begin{thm}\label{45.3} If $\Omega$ is a countable $D$-relation such that each `branch point' has `branching number' $s$, where $3\leq s\leq \aleph_0$. Then $\Omega$ is $\aleph_0$-categorical and relatively $k$-transitive for all finite $k$. Hence $\Aut(\Omega,D)$ is $3$-transitive. It is a Jordan group with $2$-transitive Jordan sets. In fact the stabilizer of a point $\alpha$ is $\Aut(\Omega-\{\alpha\},C)$ where $C$ is the corresponding homogeneous $C$-relation.
\end{thm}

\chapter{Jordan Groups II}

Recall: $G$ primitive on $\Omega$, $H\leq G$, $H\leq G_{(\Omega-\Gamma)}$, $H$ transitive on $\Gamma$, $|\Gamma|\geq 2$.

\medskip
\noindent\textbf{Lemma 37.1.} \emph{$\Gamma g$ is a Jordan set if $g\in G$. The union of a connected family of Jordan sets is a Jordan set.}

\begin{lem}\label{46.1} Let $\Gamma_i$ be Jordan sets with associated groups $H_i$, and suppose the family $(\Gamma_i)$ is connected.
\begin{enumerate}
\item[(1)] If the $\Gamma_i$ are all primitive then $\bigcup \Gamma_i$ is primitive.
\item[(2)] If the $\Gamma_i$ are (properly) $k$-homogeneous then $\bigcup \Gamma_i$ is $k$-homogeneous.
\item[(3)] If the $\Gamma_i$ are $k$-transitive then $\bigcup \Gamma_i$ is $k$-transitive.
\end{enumerate}
[Properly $k$-homogeneous: $k$-homogeneous of degree $\geq 2k$.]
\end{lem}
\begin{proof}(1) Suppose that the $\Gamma_i$ are all primitive. Let $H=\gen{H_i}_{i\in I}$, $\Gamma=\bigcup \Gamma_i$. let $\rho$ be a non-trivial $H$-congruence on $\Gamma$. Choose $\alpha,\beta$ such that $\alpha\neq\beta$ and $\alpha\equiv\beta\bmod \rho$. Suppose $\alpha\in\Gamma_i$, $\beta\not\in\Gamma_i$. Operating with elements of $H_i$ we see that $\alpha\equiv\gamma\bmod \rho$ for all $\gamma\in\Gamma_i$. So without loss of generality we can suppose that $\beta\in\Gamma_i$. But $\rho|\Gamma_i$ is an $H_i$-congruence. By primitivity $\rho|\Gamma_i$ is universal. Now use connectedness.
\end{proof}

\begin{cor}\label{46.2} Given that $G$ is primitive on $\Omega$, if $G$ has a $k$-transitive ($k$-homogeneous) Jordan set, then $G$ is $k$-transitive ($k$-homogeneous).
\end{cor}
\begin{proof} Let $\Gamma_0$ be a $k$-transitive Jordan set. Let $\Gamma^*=\bigcup \{\Gamma\mid \Gamma_0\subseteq\Gamma\text{ and $\Gamma$ is $k$-transitive Jordan set}\}$. Then by  \ref{46.1} $\Gamma^*$ is a $k$-transitive Jordan set. Also $\Gamma^*$ is maximal amongst such. Suppose that $\Gamma^*\neq \Omega$. Since $G$ is primitive there exist $g\in G$ such that $\Gamma^*g\cap\Gamma^*\neq\emptyset$ and   $\Gamma^* g\neq \Gamma^*$. If $\Gamma^* g\subseteq \Gamma^*$ then replace $g$ by $g^{-1}$. So we may suppose $\Gamma^* g\not\subseteq \Gamma^*$. Then $\Gamma^*g\cup\Gamma^*$ is a $k$-transitive Jordan set properly containing $\Gamma^*$, a contradiction.
\end{proof}

\begin{lem}\label{46.3} Let $\Gamma_1,\Gamma_2$ be Jordan sets such that $\Gamma_1\cap\Gamma_2\neq\emptyset$, $\Gamma_1-\Gamma_2\neq\emptyset$, $\Gamma_2-\Gamma_1\neq \emptyset$, $\Gamma=\Gamma_1\cup\Gamma_2$.
\begin{enumerate}
\item If $\Gamma_1,\Gamma_2$ are primitive then $\Gamma_1\cup\Gamma_2$ is $2$-homogeneous.
\item If $\Gamma_1,\Gamma_2$ are $k$-transitive then $\Gamma$ is at least $(k+1)$-transitive.
\end{enumerate}
\end{lem}
\begin{proof} (i) Choose $\alpha\in\Gamma_1-\Gamma_2$, $\beta\in\Gamma_2-\Gamma_1$. $H$ is the group on $\Gamma$, $H_i$ on $\Gamma_i$. Let $\gamma,\delta$ be any two distinct points of $\Gamma$. Since $H$ is primitive on $\Gamma$, there exists $h_0\in H$ such that one of $\gamma h_0,\delta h_0$ lies in $\Gamma_1$ and the other not. Say $\gamma h_0\in\Gamma_1$, $\delta h_0\not\in\Gamma_1$. Choose $h_1\in H_1$, mapping $\gamma h_0$ to $\alpha$. Then $\delta h_0h_1=\delta h_0$. Now choose $h_2\in H_2$ mapping $\delta h_0$ to $\beta$. Then $\alpha h_2=\alpha$ so $\gamma h_0h_1h_2=\alpha$; and $\delta h_0h_1h_2=\delta h_0h_2=\beta$. Then $h_0h_1h_2$ maps $\{\gamma,\delta\}$ to $\{\alpha,\beta\}$.

\medskip
\noindent (ii) Take $k=2$ and get $3$-transitivity. If $|\Gamma_1-\Gamma_2|\leq 2$ then $\Gamma$ is at least $(k+|\Gamma_1-\Gamma_2|)$-fold transitive. So suppose that $|\Gamma_1-\Gamma_2|\geq 3$. Choose $\alpha_1,\alpha_2$ distinct in $\Gamma_1-\Gamma_2$, $\beta\in\Gamma_2-\Gamma_1$. Let $\delta_1,\delta_2,\delta_3$ be distinct points of $\Gamma$. Since $H$ is $2$-transitive we can find $h_0\in H$ such that $\delta_1h_0=\alpha_1$, $\delta_3h_0=\beta$. If $\delta_2 h_0\in\Gamma_2$ choose $h_2\in H_2$ mapping $\delta_2 h_0$, $\beta$ to $\delta'$, $\beta$, where $\delta'\in\Gamma_1\cap\Gamma_2$. Replace $h_0$ by $h_0h_2$; we can assume that $\delta_2h_0\in\Gamma_1$. Now by $2$-transitivity of $\Gamma_1$, we have $h_1\in H_1$, mapping $\alpha_1$, $\delta_2 h_0$ to $\alpha_1$, $\alpha_2$. Then $h_0h_1$ maps $\delta_1,\delta_2,\delta_3$ to $\alpha_1,\alpha_2,\beta$. So $H$ is $3$-transitive.
\end{proof}

\begin{lem}\label{46.4} Let $H$ be a transitive group on a set $\Gamma$ and let $\Gamma_0$ be a Jordan set for $H$. Let
\begin{align*} \rho^*&=\text{congruence generated by $\Gamma_0\times \Gamma_0$}
\\ \rho_0&=\bigcap_{h\in H} \{(\Gamma_0\times \Gamma_0\cup (\Gamma-\Gamma_0)\times (\Gamma-\Gamma_0))h\mid h\in H\}.
\end{align*}
If $\rho$ is any $H$-congruence on $\Gamma$ then either $\rho\leq \rho_0$ or $\rho^*\leq \rho$. If, furthermore, there is a Jordan set $\Gamma_1$ such that $\Gamma=\Gamma_0\cup\Gamma_1$, $\Gamma\neq\Gamma_0$, $\Gamma\neq\Gamma_1$, $\Gamma_0\cap\Gamma_1\neq\emptyset$, then $\rho^*$ is universal, so $\rho_0$ is the unique maximal proper congruence on $\Gamma$ and $H$ acting on $\Gamma/\rho_0$ is $2$-homogeneous.
\end{lem}
\begin{proof} Suppose that $\rho\not\leq\rho_0$. Then there exist $\alpha\in\Gamma_0$, $\beta\not\in\Gamma_0$, $\alpha\equiv \beta\bmod \rho$. Applying elements of $H_0$ we see that $\beta\equiv\alpha$ for all $\alpha\in\Gamma_0$. So $\Gamma_0\times \Gamma_0\subseteq \rho$, hence $\rho^*\leq \rho$.

If `furthermore' holds then all elements of $\Gamma_0$ are equivalent modulo $\rho^*$: If $\alpha\in\Gamma_0\cap\Gamma_1$, $\beta\in\Gamma_0-\Gamma_1$ then operating with elements of $H_1$ we see that $\gamma\equiv\beta\bmod \rho^*$ for all $\gamma\in\Gamma_1$. Hence $\rho^*$ is universal.

Now $\Gamma_0,\Gamma_1$ are unions of $\rho_0$-classes and $\Gamma_0/\rho_0$, $\Gamma_1/\rho_0$ are Jordan sets for $\Gamma/\rho_0$ which is primitive, and Lemma  \ref{46.3} required only that $H$ was primitive on $\Gamma$. Hence $H$ acting on $\Gamma/\rho_0$ is $2$-homogeneous.
\end{proof}

\chapter{Jordan Groups III}

\begin{thm}\label{47.1} Let $G$ be a primitive Jordan group. Then either $G$ is $2$-transitive or there is a $G$-invariant relation $R$ such that one of the following holds:
\begin{enumerate}
\item $R$ is a linear order (and $G$ is highly homogeneous);
\item $R$ is a semilinear order (and $G$ has primitive proper Jordan sets);
\item $R$ is a $C$-relation.
\end{enumerate}
\end{thm}

\begin{wrapfigure}{r}{7.5cm}
\begin{tikzpicture}\node(x) at (0,-0.3) {};
\draw [draw=black, rounded corners] (0,0) rectangle (7,4);
\node(O) at (6.7, 3.7) {$\Omega$};
\draw [draw=black, rounded corners] (0.5,0.8) rectangle (2.0,2.3);
\node(G) at (0.7, 2.05) {$\Gamma$};
\node(A) at (2.2, 2.275) {\textbullet};
\node(A) at (2.4, 2.3) {$\alpha$};
\draw [draw=black, rounded corners] (2.8,1.8) rectangle (6.0,3.0);
\node(G) at (5.7, 2.7) {$\Gamma_0$};
\draw [draw=black, rounded corners] (1.5,1.8) rectangle (4.7,3.0);
\node(G) at (1.75, 2.7) {$\Gamma_1$};

\end{tikzpicture}
\end{wrapfigure}
\noindent
\emph{Proof.} Suppose that $G$ does not leave any relation as in (i)--(iii) invariant. By  \ref{44.3} if $\Gamma\subseteq\Omega$, $|\Gamma|\geq 2$, then $(\exists g\in G)(\Gamma g-\Gamma\neq\emptyset$, $\Gamma\not\subseteq \Gamma g$ and $\Gamma\cap\Gamma g\neq\emptyset)$. Let $\Gamma_0$ be a Jordan set, maximal avoiding $\alpha$. Choose $g$ so that if $\Gamma_1=\Gamma_0g$, then $\Gamma_1-\Gamma_0,\Gamma_0-\Gamma_1,\Gamma_0\cap\Gamma_1$ are all non-empty. Write $\Gamma_2=\Gamma_0\cup\Gamma_1$. So $\Gamma_2$ is a Jordan set, strictly containing $\Gamma_0$, so $\alpha\in\Gamma_2-\Gamma_0$, that is $\alpha\in \Gamma_1-\Gamma_0$. Let $H_i=G_{(\Omega-\Gamma_i)}$. Let $\beta=\alpha g\in\Gamma_0-\Gamma_1$; then $\alpha g\not\in\Gamma_0g=\Gamma_1$ but $\Gamma_0g$ is a maximal Jordan set avoiding $\alpha g$. Choose $h\in H_1$ mapping $\beta g$ to $\alpha$. Let $x=gh$. Then $x$ interchanges $\alpha,\beta$. By Lemma  \ref{46.4} there is a unique maximal proper $H_2$-congruence $\rho$ on $\Gamma_2$ and $H_2$ is $2$-homogeneous on $\Gamma_2/\rho$.  Suppose that $\Gamma_2\neq \Omega$. Then there is a translate $\Gamma$ of $\Gamma_2$ such that $\Gamma-\Gamma_2,\Gamma_2-\Gamma_1,\Gamma_1\cap\Gamma_2$ are all non-empty. If $\Gamma$ meets all $\rho$-classes choose $h_2\in H_2$ such that $\Gamma_3=\Gamma h_2$ avoids $\alpha$. Then $\Gamma_3\cap\Gamma_0\neq\emptyset$ and $\Gamma_3-\Gamma_0\neq\emptyset$. So $\Gamma_3\cup\Gamma_0$ is a Jordan set strictly containing $\Gamma_0$ and avoiding $\alpha$, which is not possible. So $\Gamma$ does not meet all $\rho$-classes. Thus choose $h_2\in H_2$ such that one of $\rho(\alpha)$, $\rho(\beta)$ avoids $\Gamma_3=\Gamma h_2$, the other meets it. Suppose that $\Gamma_3$ avoids $\rho(\alpha)$, meets $\rho(\beta)$. Then again $\Gamma_3\cup\Gamma_0$ is a Jordan set avoiding $\alpha$; not possible. If $\Gamma_3$ avoids $\rho(\beta)$ not $\rho(\alpha)$ then $\Gamma_3\cup\Gamma_1$ is a Jordan set avoiding $\beta$ and containing $\Gamma_1$ properly; not possible. Thus $\Gamma_2=\Omega$, so $\rho$ is equality, $G$ is $2$-homogeneous. But it has a permutation interchanging $\alpha,\beta$ so $G$ is $2$-transitive.
\hfill
$\Box$

\begin{thm}\label{47.2} Suppose that $G$ is $3$-primitive (i.e., $3$-transitive and $G_{\alpha\beta}$ is primitive on $\Omega-\{\alpha,\beta\}$). If $G$ has a primitive proper Jordan set, then $G$ is highly transitive.
\end{thm}

\noindent
\emph{Proof.} Note that this Jordan set must be co-infinite.  Suppose that $G$ is not highly transitive. By Peter Cameron's theorem $G$ is not highly homogeneous. By Corollary  \ref{46.2}, there exists $m$ such that $G$ has $m$-homogeneous (primitive) proper Jordan sets but not $m'$-homogeneous proper Jordan sets for $m'>m$. Define 
\[\ms G=\{\Gamma\mid\Gamma\text{ is $m$-homogeneous (primitive) proper Jordan set}\}.\]. 
\begin{wrapfigure}{r}{7.5cm}
\begin{tikzpicture}\node(x) at (0,-0.3) {};
\draw [draw=black, rounded corners] (0,0) rectangle (7,5);
\node(O) at (6.7, 4.7) {$\Omega$};
\node(A) at (2.2, 3.275) {\textbullet};
\node(A) at (2.4, 3.3) {$\alpha$};
\node(B) at (2.2, 2.275) {\textbullet};
\node(B) at (2.4, 2.3) {$\beta$};
\node(G) at (5.2, 3.275) {\textbullet};
\node(G) at (5.4, 3.3) {$\gamma$};
\node(D) at (5.2, 2.275) {\textbullet};
\node(D) at (5.4, 2.3) {$\delta$};
\draw [draw=black, rounded corners] (1.8,1.0) rectangle (2.8,4.6);
\node(G) at (2.1, 4.35) {$\Gamma_1$};
\draw [draw=black, rounded corners] (4.8,1.0) rectangle (5.8,4.6);
\node(G) at (5.1, 4.35) {$\Gamma_2$};
\draw [draw=black, rounded corners] (1.0,2.8) rectangle (6.0,3.6);
\node(G) at (1.3, 3.3) {$\Gamma_3$};
\end{tikzpicture}
\end{wrapfigure}
Suppose first that $\exists \Gamma_1,\Gamma_2\in\ms G$, $\Gamma_1\cap\Gamma_2=\emptyset$; choose distinct $\alpha,\beta\in\Gamma_1$, distinct $\gamma,\delta\in\Gamma_2$. There exists $g_1\in G$ such that $\Gamma_1g_1$ contains $\alpha$ but not $\beta$. But $G_{\alpha\beta}$ is primitive, so there exists $g_2\in G_{\alpha\beta}$ such that $\Gamma_1g_1g_2$ contains one but not both of $\gamma,\delta$, say $\gamma$, not $\delta$. Let $\Gamma_3=\Gamma_1g_1g_2$. So $\alpha,\gamma\in\Gamma_3$, $\beta,\delta\not\in\Gamma_3$. Now $\Gamma_3\cup\Gamma_2$ is an $(m+1)$-homogeneous Jordan set (Lemma  \ref{46.3}). Therefore $\Gamma_3\cup\Gamma_2$ is improper, hence cofinite.  So $\Gamma_1-\Gamma_3$ is finite. But then $\Gamma_1\cup\Gamma_3$ is an $(m+1)$-homogeneous proper Jordan set which is impossible. Alternatively, for all $\Gamma_,\Gamma_2\in\ms G$, $\Gamma_1\cap\Gamma_2\neq \emptyset$, so all $\Gamma\in\ms G$ are semistable. Now apply:

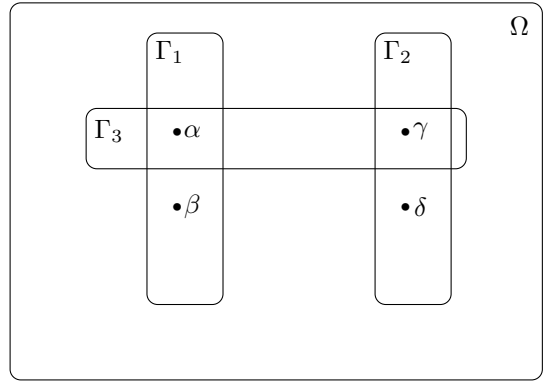
\begin{lem}\label{47.3} (This should have been in Lecture 40.) Suppose that $G$ is primitive on $\Omega$ and has a primitive semistable Jordan set. Then $G$ is highly homogeneous.
\end{lem}
The proof of this is similar to Theorem 40.4.  Exercise.
\hfill $\Box$

\chapter{Jordan Groups IV}

\begin{thm}\label{48.1} Suppose that $G$ is primitive and has primitive Jordan sets. If $G$ is not highly homogeneous, then there is a $G$-invariant relation $R$ which is a semilinear order, a $C$-relation, a $B$-relation or a $D$-relation (dense in each case and various other nice things).
\end{thm}

Let us assume that $G$ is not highly homogeneous. Can assume that $G$ is $2$-transitive (else $G$ is contained in $\Aut(\Omega,R)$ where $R$ is a dense semilinear order or a $C$-relation (\ref{47.1})). By Corollary  \ref{46.2} there exists $m$ such that there are no $(>m)$-homogeneous Jordan sets and there are $m$-homogeneous (primitive) proper Jordan sets.

Let $\ms G=\{\Gamma\mid \Gamma\text{ primitive, $m$-homogeneous proper Jordan set}\}$. By  \ref{44.3} if $\Gamma\in\ms G$ then $\exists g\in G$ such that $\Gamma g-\Gamma\neq\emptyset$, $\Gamma-\Gamma g\neq \emptyset$, $\Gamma\cap\Gamma g\neq\emptyset$. By  \ref{47.3} if $\Gamma\in\ms G$ then $\exists g\in G$ such that $\Gamma \cap\Gamma g= \emptyset$.

\begin{lem}\label{48.2} If $\Gamma_0,\Gamma_1\in\ms G$ then there exists $g\in G$ such that $\Gamma_1 g\subset \Gamma_0$.
\end{lem}

\begin{wrapfigure}[8]{r}{7.5cm}
\begin{tikzpicture}
\node(x) at (0,-0.3) {};
\draw [draw=black, rounded corners] (0,0) rectangle (7,5);
\node(O) at (6.7, 4.7) {$\Omega$};
\node(A) at (3.2, 3.675) {\textbullet};
\node(A) at (3.45, 3.7) {$\gamma_1$};
\node(G) at (5.2, 3.675) {\textbullet};
\node(G) at (5.45, 3.7) {$\gamma_2$};
\draw [draw=black, rounded corners] (2.8,1.0) rectangle (3.8,4.6);
\node(G) at (3.1, 1.2) {$\Gamma_1$};
\draw [draw=black, rounded corners] (4.8,1.0) rectangle (5.8,4.6);
\node(G) at (5.1, 1.2) {$\Gamma_2$};
\draw [draw=black, rounded corners] (1.8,3.3) rectangle (6.0,4.1);
\node(G) at (2.2, 3.5) {$\Gamma_0x$};
\draw [draw=black, rounded corners] (0.3,1.6) rectangle (4.5,2.4);
\node(G) at (0.7, 1.8) {$\Gamma_0$};
\end{tikzpicture}
\end{wrapfigure}
\noindent
\emph{Proof.}  Choose $\Gamma_2=\Gamma_1y$ such that $\Gamma_2\cap\Gamma_1=\emptyset$. Let $\gamma_1\in\Gamma_1$, $\gamma_2\in\Gamma_2$ and $x\in G$ such that $\Gamma_0x$ contains $\gamma_1,\gamma_2$. Suppose, if possible, that $\Gamma_1\not\subseteq \Gamma_0x$, $\Gamma_2\not\subseteq \Gamma_0x$. Then $\Gamma_1\cup\Gamma_0x$ is a Jordan set that is at least $(m+1)$-homogeneous. Therefore it is improper. So $\Gamma_2-\Gamma_0x$ is finite. Then $\Gamma_2\cup\Gamma_0x$ is $(m+1)$-homogeneous so $\Gamma_1-\Gamma_0x$ is finite. Hence $\Gamma_0x$ is cofinite, which is false. Hence $\Gamma_1\subseteq \Gamma_0x$ or $\Gamma_2\subseteq \Gamma_0x$, so $(\exists g)(\Gamma_1g\subset \Gamma_0)$.
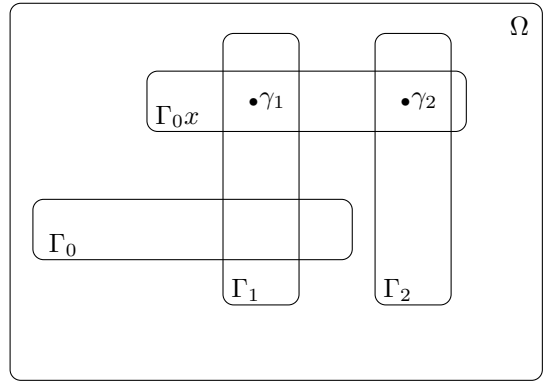
\hfill $\Box$

\vspace{2.0cm}
\begin{lem}\label{48.3} Temporarily assume that $G$ is primitive on $\Omega$ and that $\Gamma$ is a Jordan set. If $\alpha,\beta\in\Omega$ then $\exists g\in G$ such that $\Gamma g$ contains $\alpha,\beta$.
\end{lem}
\begin{wrapfigure}[5]{r}{3.5cm}
\begin{tikzpicture}
\draw [draw=black, rounded corners] (0,0) rectangle (3.6,3);
\node(O) at (3.3, 2.7) {$\Omega$};
\node(A) at (0.6, 0.875) {\textbullet};
\node(A) at (0.9, 0.9) {$\omega_1$};
\node(B) at (1.7, 0.875) {\textbullet};
\node(B) at (2.0, 0.9) {$\omega_2$};
\node(C) at (1.7, 1.775) {\textbullet};
\node(C) at (2.0, 1.8) {$\omega_3$};
\draw [draw=black, rounded corners] (1.4,0.3) rectangle (2.5,2.5);
\node(G) at (1.7, 2.2) {$\Gamma_2$};
\draw [draw=black, rounded corners] (0.3,0.3) rectangle (2.5,1.4);
\node(G) at (0.6, 0.5) {$\Gamma_1$};
\end{tikzpicture}\end{wrapfigure}
\noindent
\emph{Proof.}  Define $\omega_1\equiv\omega_2$ if $(\exists g)(\omega_1,\omega_2\in\Gamma g)$. The relation is reflexive and symmetric. For transitivity suppose that $\omega_1,\omega_2\in\Gamma_1=\Gamma g_1$, $\omega_2,\omega_3\in\Gamma_2=\Gamma g_2$.

Either $\omega_3\in\Gamma_1$ , in which case $\Gamma_1$ shows $\omega_1\equiv\omega_3$; or $\omega_1\in\Gamma_1-\Gamma_2$, $\omega_3\in\Gamma_2-\Gamma_1$. In this case $\exists h\in G_{(\Omega-\Gamma_2)}$ such that $\omega_2h=\omega_3$. Then $\Gamma_1h$ witnesses $\omega_1\equiv\omega_3$.
\hfill $\Box$

\begin{proof}[Proof of Theorem  \ref{48.1}] \textbf{Case 1.} $G$ is not $2$-primitive, that is, $G_\alpha$ is not primitive on $\Omega-\{\alpha\}$. Then we get a $B$-relation.

Define $B(\alpha;\beta,\gamma)\Leftrightarrow (\Gamma\in\ms G\wedge \beta,\gamma\in\Gamma\Rightarrow \alpha\in\Gamma)$. We check the axioms.
\begin{enumerate}
\item[(B1)] $B(\alpha;\beta,\gamma)\Rightarrow B(\alpha;\gamma,\beta)$.
\item[(B2)] Certainly $\alpha=\beta\Rightarrow B(\alpha;\beta,\gamma)\wedge B(\beta;\alpha,\gamma)$. Conversely, suppose $\alpha\neq\beta$ and $B(\alpha;\beta,\gamma)$. Choose $\Gamma_1,\Gamma_2\in\ms G$, $\Gamma_1\cap\Gamma_2\neq\emptyset$, $\alpha\in\Gamma_1-\Gamma_2$, $\beta\in\Gamma_2-\Gamma_1$. Then $\Gamma_1\cup\Gamma_2$ is $(m+1)$-homogeneous Jordan set so cannot be proper, hence is cofinite. Since $G$ is not $2$-primitive, $\Gamma_1\cup\Gamma_2=\Omega$. Then $\gamma\in\Gamma_1$ or $\gamma\in\Gamma_2$. But $\gamma\not\in\Gamma_2$ because $\alpha\not\in\Gamma_2$ and $B(\alpha;\beta,\gamma)$. Therefore $\gamma\in\Gamma_1-\Gamma_2$. Hence $\Gamma_1$ shows $\lnot B(\beta;\alpha,\gamma)$.
\item[(B3)] Want to show $B(\alpha;\beta,\gamma)\Rightarrow B(\alpha;\beta,\delta)\vee B(\alpha;\gamma,\delta)$. Suppose $\lnot B(\alpha;\beta,\delta)$ and $\lnot B(\alpha;\gamma,\delta)$. Then there exists $\Gamma_1,\Gamma_2\in\ms G$ such that $\alpha\not\in\Gamma_1$, $\beta,\delta\in\Gamma_1$, $\alpha\not\in\Gamma_2$, $\gamma,\delta\in\Gamma_2$. Then $\Gamma_1\cap\Gamma_2\neq \emptyset$. If $\Gamma_1-\Gamma_2\neq\emptyset$ and $\Gamma_2-\Gamma_1\neq\emptyset$, then $\Gamma_1\cup\Gamma_2=\Omega$ which is not so since $\alpha\not\in\Gamma_1\cup\Gamma_2$. Thus $\Gamma_1\subseteq \Gamma_2$ or $\Gamma_2\subseteq\Gamma_1$. Then $\Gamma_2$ shows $\lnot B(\alpha;\beta,\gamma)$ or $\Gamma_1$ shows $\lnot B(\alpha;\beta,\gamma)$.
\end{enumerate}
The other betweenness axioms are similarly verified.  Infact we have a genuine betweenness relation.

\medskip

\noindent\textbf{Case 2.} $G$ is $2$-primitive. Get $D$-relation. Define \[D(\alpha,\beta;\gamma,\delta)\Leftrightarrow (\exists \Gamma_1,\Gamma_2\in\ms G)(\Gamma_1\cap\Gamma_2= \emptyset\wedge \alpha,\beta\in\Gamma_1\wedge \gamma,\delta\in\Gamma_2).\]
\begin{enumerate}
\item[(D1)] Obvious.
\item[(D2)] Suppose $D(\alpha,\beta;\gamma,\delta)$ and $D(\alpha,\gamma;\beta,\delta)$. Not possible.
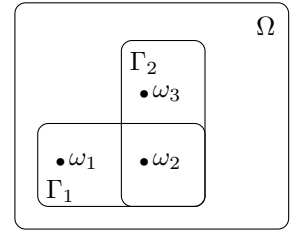
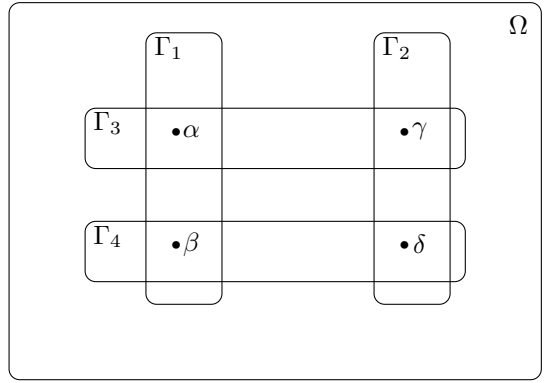
\begin{center}
\begin{tikzpicture}
\draw [draw=black, rounded corners] (0,0) rectangle (7,5);
\node(O) at (6.7, 4.7) {$\Omega$};
\node(A) at (2.2, 3.275) {\textbullet};
\node(A) at (2.4, 3.3) {$\alpha$};
\node(B) at (2.2, 1.775) {\textbullet};
\node(B) at (2.4, 1.8) {$\beta$};
\node(G) at (5.2, 3.275) {\textbullet};
\node(G) at (5.4, 3.3) {$\gamma$};
\node(D) at (5.2, 1.775) {\textbullet};
\node(D) at (5.4, 1.8) {$\delta$};
\draw [draw=black, rounded corners] (1.8,1.0) rectangle (2.8,4.6);
\node(G) at (2.1, 4.4) {$\Gamma_1$};
\draw [draw=black, rounded corners] (4.8,1.0) rectangle (5.8,4.6);
\node(G) at (5.1, 4.4) {$\Gamma_2$};
\draw [draw=black, rounded corners] (1.0,2.8) rectangle (6.0,3.6);
\node(G) at (1.3, 3.4) {$\Gamma_3$};
\draw [draw=black, rounded corners] (1.0,1.3) rectangle (6.0,2.1);
\node(G) at (1.3, 1.9) {$\Gamma_4$};
\end{tikzpicture}
\end{center}
\item[(D3)] Suppose $D(\alpha,\beta;\gamma,\delta)$ and $D(\alpha,\ep;\gamma,\delta)\vee D(\alpha,\beta;\gamma,\ep)$.  

Suppose $\Gamma_1$ is maximal avoiding $\gamma$ and $\delta$.

Can suppose $\ep\not\in \Gamma_1\cup\Gamma_2$. Since $G_\ep$ is primitive om $\Omega-\{\ep\}$, there is a $G_\ep$-translate of $\Gamma'$ containing $\alpha,\delta$. Now $\Gamma'\supseteq \Gamma_1$, $\Gamma'\supseteq \Gamma_2$ else in either case $\Gamma'$ is $(n+1)$-homogeneous with infinite $\Gamma$ in complement. Thus $\Gamma\cap\Gamma_1=\Gamma\cap\Gamma_2=\emptyset$.

Now suppose $\Gamma$ is maximal avoiding $\alpha,\delta$. Now expand $\Gamma_2$ to be maximal avoiding $\alpha,\ep$. Hence $\Gamma_2$ avoids $\Gamma_1,\Gamma$. Take $\Gamma_1$ maximal avoiding $\delta,\ep$. Now use Lemma  \ref{48.2} and work inside $G_\delta$. There exists $x\in G_\delta$ such that $\Gamma_1 x\subset \Gamma_1$ and by  \ref{48.3} applied to $G_{(\Omega-\Gamma_1)}$ there exists $y\in G_{(\Omega-\Gamma_1)}$ such that $\alpha,\beta\in\Gamma_1xy$. So $\alpha,\beta\in\Gamma_1xy\subset \Gamma_1$, $xy\in G_\delta$. Then $\Gamma_2xy\supset \Gamma_2$. Then $\ep\in\Gamma_2 xy$. So $\Gamma_2xy$, $\Gamma_1xy$ show $D(\alpha,\beta;\ep,\delta)$.
\end{enumerate}
\end{proof}

\begin{prob} Classify all Jordan groups.\footnote{Infinite primitive Jordan groups are classified in S.\ A.\ Adeleke and Dugald Macpherson, Dugald, \lq Classification of infinite primitive Jordan permutation groups\rq, \emph{Proc.\ London Math.\ Soc.\ (3)} 72 (1996), 63--123.}
\end{prob}

\noindent\textbf{Conjecture.} If $G$ is not highly transitive ($k$- but not $(k+1)$-transitive for $k\geq 4$) there exists $G$-invariant Steiner system on $\Omega$.

\chapter*{References for Lectures 33--48\\ (Trinity Term)}

\begin{itemize}
\item[1.]\label{ref:lec33-48} Peter J.\ Cameron, \lq Transitivity of permutation groups on unordered sets\rq, \emph{Math.\ Z.}, 148 (1976), 127--139.
\item[2.] Graham Higman \lq Homogeneous relations\rq, \emph{Quart.\ J.\ Math. Oxford (2)}, 28 (1977), 31--39.
\item[3.] Wilfrid Hodges, A.H.\ Lachlan and Saharon Shelah, \lq Possible orderings of an indiscernible sequence\rq, \emph{Bull.\ London Math.\ Soc.}, 9 (1977), 212--215.
\item[4.] Peter M.\ Neumann, \lq Homogeneity of infinite permutation groups\rq, \emph{Bull.\ London Math.\ Soc.}, 20 (1988), 305--312.
\item[5.] G.\ Cherlin, L.\ Harrington and A.H.\ Lachlan, \lq $\aleph_0$-categorical, $\aleph_0$-stable structures\rq, \emph{Ann.\ Pure Appl.\ Logic}, 28 (1985), 103--135.
\item[6.] Peter M.\ Neumann, \lq Some primitive permutation groups\rq, \emph{Proc.\ London Math.\ Soc.\ (3)}, 50 (1985), 265--281. 
\item[7.] David M.\ Evans, \lq  Homogeneous geometries\rq, \emph{Proc.\ London Math.\ Soc.\ (3)}, 52 (1986), 305--327.
\item[8.]   B.I.\ Zil'ber, \lq Strongly minimal countably categorical theories, I, II, III\rq, \emph{Siberian Math.\ J.}, 21 (1980), 219--230,  25 (1984), no.\ 3, 71--88, no.\ 4, 63--77.
\item[9.]  Manfred Droste, \lq Structure of partially ordered sets with transitive automorphism group\rq, \emph{Mem.\ Amer.\ Math.\ Soc.} 57 (1985), no.\ 334, pp.\ vi+100.
\item[10.] S.A.\ Adeleke, M.A.E.\ Dummett and Peter M.\ Neumann, \lq On a question of Frege about right-ordered groups\rq, \emph{Bull.\ London Math.\ Soc.}, 19 (1987), 513--521.
\item[11.] S.A.\ Adeleke and Peter M.\ Neumann, \lq Relations related to betweeness: their structure and automorphisms\rq \emph{Mem.\ Amer.\ Math.\ Soc.} 131 (1998), no.\ 623, pp.\ viii+125.
\item[12.] S.A.\ Adeleke and Peter M.\ Neumann, \lq On infinite Jordan groups\rq, Primitive permutation groups with primitive Jordan sets.
\emph{J. London Math. Soc. (2)}, 53 (1996), 209--229. 
\item[12.] S.A.\ Adeleke and Peter M.\ Neumann, \lq Infinite bounded permutation groups\rq, \emph{J.\ London Math.\ Soc.\ (2)}, 53 (1996), 230-242.
\end{itemize}

\end{document}